\documentclass[12pt, a4paper, oneside]{article}
\usepackage{amsmath, amsthm, amssymb, appendix, bm, graphicx, hyperref, mathrsfs,geometry,color,cite}
\usepackage{booktabs,microtype,enumitem}
\usepackage{array}
\usepackage{caption}
\usepackage{subfigure}
\newtheorem{theorem}{Theorem}[section]
\numberwithin{equation}{section}

\newtheorem{lemma}{Lemma}[section]

\newtheorem{example}{Example}[section]
\newtheorem{proposition}{Proposition}[section]
\newtheorem{remark}{Remark}[section]

\usepackage{sectsty}
\sectionfont{\fontsize{14}{16}\selectfont}

\usepackage{authblk}
\makeatletter
\renewcommand{\maketitle}{
	\begin{center}
		{\LARGE\textbf{ \@title} \par}
		\vskip 1em
		{ \@author \par}
		\vskip 1em
		{\small \@date \par}
	\end{center}
	\vskip 1em
	\@thanks
}
\makeatother

\title{Global dynamics of nonautonomous parabolic systems with spatiotemporal heterogeneity and anisotropy}

\author[a]{Hailong Ye}
\author[b]{Xiao-Qiang Zhao\thanks{Corresponding author: zhao@mun.ca}}

\affil[a]{School of Mathematical Sciences, Shenzhen University, Shenzhen 518061, China}
\affil[b]{Department of Mathematics and Statistics, Memorial University of Newfoundland,
	St. John's, NL A1C 5S7, Canada}

\date{}

\begin{document}

\maketitle

\begin{abstract}

\smallskip

In this paper, we investigate nonautonomous parabolic systems on a bounded domain, driven by differential operators with spatiotemporally heterogeneous and anisotropic diffusion as well as advection, subject to Dirichlet boundary conditions. We establish the global dynamics under asymptotically time-periodic, asymptotically partially and fully degenerate regimes, respectively.
It turns out that in the asymptotically degenerate case, the long-term behavior of such a system is no longer determined by its limiting system alone.  We  further address  potential connections between our results and the dichotomy/trichotomy phenomena for monostable evolution systems with moving boundary conditions. As an illustration, we also apply the developed theory
to a large class of reaction-diffusion-advection systems on time-varying domains.
\end{abstract}

\section{Introduction}

A wide variety of biological processes and population dynamics exhibit spatiotemporal heterogeneity and anisotropy.
These features not only emerge naturally from environmental variations but also reflect underlying structural or
mechanical anisotropy \cite{Jin-Hilker-Steffler-Lewis2014,Timm-Okubo1992}.
For instance, these effects have been incorporated into models to explain stripe orientations in Turing instability \cite{Shoji-Iwasa-Mochizuki-Kondo2002,Shoji-Mochizuki-Iwasa-Hirata-Watanabe-Hioki-Kondo2003}, to model and predict
tumor growth \cite{Jbabdi-Mandonnet-Duffau,Painter-Hillen2013}, and to determine evolutionarily stable strategies
for multiple competing species in heterogeneous habitats \cite{Bouin-Legendre-Lou-Slover2021,Yu-Bu-Wu2025}.
Biologically, spatiotemporal dependence and anisotropy correspond to distinct mechanisms in movement processes:
the former is typically driven by environmental fluctuations across space and time, whereas the latter reflects
directional biases in diffusion or advection caused by intrinsic tissue architecture or extrinsic cues.

Another motivation comes from the increasing recognition of the importance of incorporating domain growth into
mathematical models for a variety of biological processes, including pattern formation
\cite{Crampin-Gaffney-Maini,MR2643885,Hetzer-Madzvamuse-Shen2012}, cell migration during embryonic development
\cite{Baker-Yates-Erban2010,Landman-Pettet-Newgreen2003,Simpson2015}, and various population dynamics
\cite{FanZhao2026, Fang-Li-Su2024,Lin-Ye-Zhao2026,MR4784487,Meng-Lin-Pedersen2021}.
Mathematically, under reasonable assumptions, reaction-diffusion systems on evolving domains can be rigorously recast
via a diffeomorphic transformation as nonautonomous parabolic systems on a fixed reference domain (see, e.g.,
\cite{MR2643885,Hetzer-Madzvamuse-Shen2012,Lin-Ye-Zhao2026} on evolving domains and \cite{Gorder-Klika-Krause2021,Krause-Ellis-Gorder2019} on
evolving manifolds). This reformulation provides a starting point for theoretical analysis and numerical simulations,
which differs from free boundary models \cite{Du-Lin2010,Feng-Lewis-Wang-Wang2022}. Despite considerable research
in this area, most progress has been restricted to isotropic and spatially uniform domain evolution.
In developing organisms or tissues, however, domain evolution is spatially heterogeneous and anisotropic,
shaped by complex cellular or subcellular structures as well as external environmental factors
\cite{Crampin-Hackborn-Maini,Page-Maini-Monk2005,Bittig-Wartlick2008}. As a result, when the governing systems are
pulled back to a fixed reference domain, one naturally obtains nonautonomous parabolic systems with spatiotemporally
dependent anisotropic differential operators. For example, using such a pullback transformation from general growing domains, Hetzer et al. \cite{Hetzer-Madzvamuse-Shen2012} derived the following nonautonomous system on a fixed reference domain:
\[
\left\{
\begin{aligned}
\frac{\partial \hat{u}}{\partial t}
&=
\sum_{j,k=1}^{n} a_{j,k}(y,t)\frac{\partial^2 \hat{u}}{\partial y_j \partial y_k}
+\sum_{j=1}^{n}\bigl(b_j(y,t)-c_j(y,t)\bigr)\frac{\partial \hat{u}}{\partial y_j}
+c_0(y,t)\hat{u}
+\gamma f(\hat{u},\hat{v}),
\\[1mm]
\frac{\partial \hat{v}}{\partial t}
&=
d_c\sum_{j,k=1}^{n} a_{j,k}(y,t)\frac{\partial^2 \hat{v}}{\partial y_j \partial y_k}
+\sum_{j=1}^{n}\bigl(d_c\,b_j(y,t)-c_j(y,t)\bigr)\frac{\partial \hat{v}}{\partial y_j}
+c_0(y,t)\hat{v}
+\gamma g(\hat{u},\hat{v}).
\end{aligned}
\right.
\]
The evolution family of the linearized  system was subsequently utilized to characterize Turing diffusion-driven
instability. Nevertheless, the threshold dynamics of such systems remains far from well understood, as the intricate
interplay of spatiotemporal heterogeneity and anisotropy in both diffusion and advection poses substantial
analytical challenges.

\vspace{1mm}
Inspired by the aforementioned works, we consider the following  cooperative parabolic systems:
\begin{equation}\label{intro-1}
\left\{
\begin{aligned}
& \frac{\partial v_i}{\partial t}   +A_i(t)v_i=f_i(x,t,v) \ \ \ &{\rm in}&\ \ \ \Omega\times(0,\infty),\\
& v_i =0 \ \ \ \ &{\rm on}&\ \ \ \partial\Omega\times(0,\infty),
\end{aligned}
\right.
\end{equation}
for $i=1,\dots,m$, where $v=(v_1,\dots,v_m)^\top$ represents the density vector of $m$ interacting species,
$\Omega\subset\mathbb{R}^n$ is a bounded domain, and
$$
A_i(t)v_i=-\sum_{j,k=1}^{n}a^{(i)}_{jk}(x,t)\partial_{jk}v_i +\sum_{j=1}^{n}b^{(i)}_{j}(x,t)\partial_jv_i.
$$
Our purpose is  to systematically explore the global dynamics of nonautonomous system \eqref{intro-1}.
In particular, we focus on its threshold dynamics under three distinct asymptotic regimes for the
operator families $\{A_i(t)\}_{t\ge0}$ ($i=1,\dots,m$):
asymptotically time-periodic regime, asymptotically partially degenerate regime, and
 asymptotically fully degenerate regime.
We also address  potential connections between our results and the dichotomy/trichotomy phenomena in some free boundary
problems, see Examples \ref{asymptotically-periodic-example-1}--\ref{asymptotically-periodic-example-2} and
Examples \ref{degenerate-example-1}--\ref{degenerate-example-2} for details.
Next, we present our main results according to these three regimes, respectively.

\medskip

{\it Asymptotically time-periodic regime.}
In this regime, the operator family $\{A_i(t)\}_{t\ge0}$ converges coefficientwise to an $\omega$-periodic operator
family $\{\tilde A_i(t)\}$ as $t\to\infty$ for each $i=1,\dots,m$. The corresponding limiting system of \eqref{intro-1}
can be written as
\begin{equation}\label{limiting system-1}
\left\{
\begin{aligned}
& \frac{\partial v_i}{\partial t}+\tilde A_i(t)v_i=\tilde f_i(x,t,v) \ \ \ &{\rm in}&\ \ \ \Omega\times(0,\infty),\\
& v_i = 0 \ \ \ \ &{\rm on}&\ \ \ \partial\Omega\times(0,\infty),
\end{aligned}
\right.
\end{equation}
where $\tilde f_i$ denotes the corresponding limiting nonlinearity.
Let $\lambda^*_\omega$ be the principal eigenvalue associated with the linearization of system  \eqref{limiting system-1} at the zero solution.

For isotropic diffusion with asymptotically periodic coefficients, Lam et al. \cite{MR4784487} used a time transformation
to reduce the problem to an equivalent system with time-independent diffusion coefficients, thereby establishing the
threshold dynamics under an integral convergence condition. This technical
condition  was subsequently removed in
\cite{FanZhao2026} by exploiting the connection between evolution families and analytic semigroups. These methods,
however, are not applicable to system \eqref{intro-1}, due to the spatiotemporal heterogeneity and anisotropy in
both diffusion and advection. Such features make the spectral structure of the associated nonautonomous operators
substantially more involved, and standard techniques based on energy or Schauder estimates are insufficient
to characterize the desired threshold dynamics.   Instead, we will use the theories of asymptotically periodic semiflows
and chain transitive sets (see \cite{MR3643081}) to lift the global dynamics of the limiting systems to the original nonautonomous 
systems.

In order to apply the dynamical systems approach,  we first show that system \eqref{intro-1} generates an asymptotically periodic semiflow, and then
construct a brand new subsolution for system  \eqref{intro-1} to ensure that every nontrivial solution
remains strongly separated from zero as $t\to\infty$ in the case where  $\lambda_\omega^{\ast}<0$.  Accordingly, we have the following result.
\begin{theorem}\label{main results-case1}
	Let \ref{ass:a1}--\ref{ass:a7} hold. Then the following statements are valid:
	\begin{itemize}
		\item [(a)] If $\lambda_\omega^{\ast}\ge0$, then the zero solution is globally asymptotically stable for \eqref{intro-1} in $\mathcal{M}_{\bar u}^{\beta}$.
		\item [(b)] If $\lambda_\omega^{\ast}<0$, then every nontrivial solution $\phi$
		of \eqref{intro-1}  with $\phi(\cdot,0)\in\mathcal{M}_{\bar u}^{\beta}$
                    satisfies $\lim_{t\to+ \infty} \|\phi(\cdot,t)-v_\omega^*(\cdot,t)\|_{\beta}=0$,
                    where $v_\omega^*$ is the unique positive $\omega$-periodic solution of the
                    limiting system  \eqref{limiting system-1}.
	\end{itemize}
\end{theorem}

Theorem \ref{main results-case1} nontrivially extends a broad class of known results including those in  \cite{MR4784487,FanZhao2026,Meng-Lin-Pedersen2021,Jiang-Wang2018},
and provides a general framework for understanding the critical habitat size
governing the vanishing (or vanishing-balancing) phenomena in some free boundary problems, see
Examples \ref{asymptotically-periodic-example-1}--\ref{asymptotically-periodic-example-2}.

Note that Theorem \ref{main results-case1} is also helpful  for the analysis of  reaction--diffusion systems on evolving domains that preserve total area or volume, see Remark \ref{reaction-diffusion-advection system-7}.
This setting seems to have received little attention in the existing literature \cite{Gorder-Klika-Krause2021}.
Such geometric constraints arise in various biological and physical contexts, including the buckling of intestinal
crypts during development \cite{Drasdo-Loeffler2001}, and cell shape changes \cite{Hunding1985} or stationary
shapes of vesicles \cite{Seifert-Berndl-Lipowsky1991}.

\medskip

{\it Asymptotically partially degenerate regime.}
Here partial degeneracy refers to the situation in which, as time increases, the diffusion and advective
effects gradually diminish and eventually disappear in certain spatial directions.
We focus on the product-domain setting $\Omega=\Omega_1\times\Omega_2\subset\mathbb{R}^{n}$,
where $\Omega_k\subset\mathbb{R}^{n_k}$($k=1,2$), and $n_1+n_2=n$.
Let $x=(x^{(1)},x^{(2)})\in\mathbb{R}^{n}$
with $x^{(1)}=(x_1,\cdots,x_{n_1})$ and $x^{(2)}=(x_{n_1+1},\cdots,x_{n})$.
For each $i=1,\dots,m$, we decompose the operator $A_i(t)$ as
$$
A_i(t)=A^{(1)}_i(t)+A^{(2)}_i(t),
$$
where
\begin{align*}
&A^{(1)}_i(t)v_i=-\sum_{j,k=1}^{n_1}a^{(i)}_{jk}(x,t)\partial_{jk}v_i +\sum_{j=1}^{n_1}b^{(i)}_{j}(x,t)\partial_jv_i,
\\
&A^{(2)}_i(t)v_i=-\sum_{j,k=n_1+1}^{n}a^{(i)}_{jk}(x^{(2)},t)\partial_{jk}v_i +\sum_{j=n_1+1}^{n}b^{(i)}_{j}(x^{(2)},t)\partial_jv_i,
\end{align*}
the coefficients of $A^{(1)}_i(t)$ decay to zero as $t\to\infty$,
while ${A^{(2)}_i(t)}$ converges coefficientwise to an $\omega$-periodic operator family ${\tilde A^{(2)}_i(t)}$
in $\Omega_2\times(0,\infty)$. Thus, the limiting system of \eqref{intro-1} takes the following form:
\begin{equation}\label{limiting system-2}
\left\{
\begin{aligned}
& \frac{\partial w_i}{\partial t}+\tilde A^{(2)}_i(t)w_i=\tilde f_i(x^{(2)},t,w) \ \ \ &{\rm in}&\ \ \ \Omega_2\times(0,\infty),\\
& w_i = 0 \ \ \ \ &{\rm on}&\ \ \ \partial\Omega_2\times(0,\infty).
\end{aligned}
\right.
\end{equation}

This partial asymptotic degeneracy may occur in multi-species competition models, offering a natural
mechanism for studying the anisotropic evolutionarily stable strategies in heterogeneous habitats \cite{Bouin-Legendre-Lou-Slover2021,Yu-Bu-Wu2025}.
It also arises in reaction--diffusion systems on anisotropically evolving domains,
where disproportionate expansion in different directions progressively stretches the corresponding spatial scales \cite{Onder-Gorder2023,Gorder-Klika-Krause2021}.
A typical example is a tubular evolving domain whose length becomes asymptotically unbounded, while its cross-section
remains asymptotically bounded. After transforming the problem into a fixed reference domain, such strongly
anisotropic geometric evolution manifests itself as the asymptotic degeneracy of the diffusion and advective
effects in the stretched directions.

Let $\lambda^*_\omega$ be the principal eigenvalue associated with  the linearization of system \eqref{limiting system-2} at the zero solution.
By Theorem \ref{main results-case1},  system \eqref{limiting system-2} admits a unique
positive $\omega$-periodic solution $w_\omega^*(x^{(2)},t)$
under the condition that $\lambda_\omega^{\ast}<0$.  By  the comparison arguments and
the construction of two novel subsolutions, we are able to prove the following result.
\begin{theorem}\label{main results-case2-1}
Assume that \ref{ass:b1}--\ref{ass:b5} hold, and let $\alpha(t)$ and $\beta(t)$ be defined as in \ref{ass:b1} and \ref{ass:b3}, respectively.  Then the following statements are valid:
\begin{itemize}
\item [(a)] If $\lambda_\omega^{\ast}\ge0$, then the zero solution is globally asymptotically stable for \eqref{intro-1} in $C_0(\overline\Omega,\mathbb{M}_{\bar{u}})$.
\item [(b)] If $\lambda_\omega^{\ast}<0$ and there exists $T>0$ such that
            \begin{equation}\label{main results-case2-1-1}
            \limsup_{t\to\infty}\frac{1}{T}\int_t^{t+T}\frac{\beta(s)^2}{\alpha(s)}\,\mathrm{d}s<4\lambda_0(-\lambda_\omega^*),
            \end{equation}
            then every nontrivial solution $u$ of \eqref{intro-1} with $u(\cdot,0)\in C_0(\overline\Omega,\mathbb{M}_{\bar{u}})$ satisfies
                    $$
                    \lim_{t\to+ \infty}\sup_{(x^{(1)},x^{(2)})\in K\times\overline\Omega_2}|u(x,t)-w_\omega^*(x^{(2)},t)|=0
                    $$
                    for every compact subset $K$ of $\Omega_1$.
\end{itemize}
\end{theorem}

In population biology, the homogeneous Dirichlet boundary condition represents a hostile exterior environment.
Along the degenerate directions, even though both diffusion and
advection vanish asymptotically, their cumulative effects may still influence the long-time dynamics through
persistent boundary loss. Condition \eqref{main results-case2-1-1} effectively regulates the long-time averaged
strength of the advective drift relative to the random diffusion. It ensures that the advective effect in the
degenerate directions is not strong enough to destroy persistence in the habitat interior.
Remarkably, this type of condition for persistence reduces to
the asymptotic spreading speed in some free boundary problems, see Examples \ref{degenerate-example-1}--\ref{degenerate-example-2} for details.
In the case where \eqref{main results-case2-1-1} does not hold, we have the following result.
\begin{theorem}\label{main results-case2-2}
Assume that \ref{ass:b1}--\ref{ass:b3}, \ref{ass:b4'}--\ref{ass:b5'}, and \ref{ass:b6} hold, and that
$\lambda_\omega^{\ast}<0$. If $p:=\frac{\chi}{\alpha}$ is locally absolutely continuous on $[t_0,\infty)$ and there exists $T>0$ such that
\begin{equation}\label{main results-case2-2-1}
\liminf_{t\to\infty}\frac{1}{T}\int_{t}^{t+T}\left(\frac{\chi(s)^2}{\alpha(s)}-2d_e\bigl(p'(s)\bigr)_+\right)\,\mathrm{d}s>4\Lambda_0(-\lambda_\omega^{\ast}),
\end{equation}
where $(x)_+=\max\{x,0\}$, $d_e=\max_{x^{(1)}\in\overline{\Omega}_1} x^{(1)}\cdot e-\min_{x^{(1)}\in\overline{\Omega}_1} x^{(1)}\cdot e$, and $\chi$ is given in \ref{ass:b6},
then every solution of \eqref{intro-1} with initial data in $C_0(\overline{\Omega}, \mathbb{R}_+^m)$ converges to zero uniformly on $\overline\Omega$ as $t\to\infty$.
\end{theorem}

Combining Theorem \ref{main results-case2-1} with Theorem \ref{main results-case2-2}, we see that vanishing diffusion and advection along the
degenerate directions can still exert a critical impact on the long-time behavior. We will also  provide
Example \ref{partial-degenerate-example} to illustrate the sharpness of both the upper bound
$4\lambda_0(-\lambda_\omega^*)$ in \eqref{main results-case2-1-1} and the lower bound
$4\Lambda_0(-\lambda_\omega^*)$ in \eqref{main results-case2-2-1}.

\medskip

{\it  Asymptotically fully degenerate regime. }
In this regime, all coefficients in $A_i(t)$ vanish asymptotically over time. This situation arises naturally in
reaction--diffusion systems posed on asymptotically unbounded domains, which can be used to describe species migrating
in constantly expanding habitats \cite{Tang-Zhang-Lin2011,MR4784487}, such as those of {\it Aedes} mosquitoes.
These mosquitoes are vectors of dengue fever, yellow fever, Zika virus, and other infectious diseases, and global
warming and frequent human activities have substantially enlarged the areas suitable for their survival and reproduction.
When the governing systems are pulled back to a fixed reference domain, the unbounded growth of the original
domain is reflected in the coefficientwise decay of the induced operators $A_i(t)$.

Assume here that $f$ is independent of $x$. Then the limiting system of~\eqref{intro-1} in this regime reduces to
\begin{equation}\label{limiting system-3}
	\frac{\mathrm{d}w}{\mathrm{d}t}=\tilde f(t,w), \qquad t>0.
\end{equation}
Let $\lambda^*$ be the principal eigenvalue of the linearization of~\eqref{limiting system-3} at the zero solution. As in the partially degenerate regime, we have the following two results.
\begin{theorem}\label{main results-case3-1}
Assume that \ref{ass:c1}--\ref{ass:c5} hold, and let $\alpha(t)$ and $\beta(t)$ be defined as in \ref{ass:c1} and \ref{ass:c3}, respectively.  Then the following statements are valid:
	\begin{itemize}
		\item [(a)] If $\lambda^*\ge0$, then the zero solution is globally asymptotically stable
		for \eqref{intro-1}
		in $C_0(\overline\Omega,\mathbb{M}_{\bar{u}})$.
		\item [(b)] If $\lambda^*<0$ and there exists $T>0$ such that
                        \begin{equation}\label{main results-case3-1-1}
                        \limsup_{t\to\infty}\frac{1}{T}\int_t^{t+T}\frac{\beta(s)^2}{\alpha(s)}\,\mathrm{d}s<4\lambda_0(-\lambda^*),
                        \end{equation}
                    then every nontrivial solution $u(x,t)$ of \eqref{intro-1}
                   with $u(\cdot,0)\in C_0(\overline\Omega,\mathbb{M}_{\bar{u}})$
                    satisfies $\lim_{t\to+ \infty}|u(\cdot,t)-w_e(t)|=0$ locally uniformly in $\Omega$, where $w_e$ is the unique positive $\omega$-periodic solution of system  \eqref{limiting system-3}.
	\end{itemize}
\end{theorem}

In the autonomous Stefan models, \eqref{main results-case3-1-1} reduces to the known upper bound $k_0<2\sqrt{ad}$ on the asymptotic spreading speed (see Example \ref{degenerate-example-1}).

\begin{theorem}\label{main results-case3-2}
Assume that \ref{ass:c1}--\ref{ass:c3}, \ref{ass:c4'}--\ref{ass:c5'}, and \ref{ass:c6} hold, and that
$\lambda^*<0$. If $p:=\frac{\chi}{\alpha}$ is locally absolutely continuous on $[t_0,\infty)$ and there exists $T>0$ such that
\begin{equation}\label{main results-case3-2-1}
\liminf_{t\to\infty}\frac{1}{T}\int_{t}^{t+T}\left(\frac{\chi(s)^2}{\alpha(s)}-2d_e\bigl(p'(s)\bigr)_+\right)\,\mathrm{d}s>4\Lambda_0(-\lambda^*),
\end{equation}
where $d_e=\max_{x\in\overline{\Omega}} x\cdot e-\min_{x\in\overline{\Omega}} x\cdot e$, and $\chi$ is given in \ref{ass:c6},
then every solution of \eqref{intro-1} with initial data in $C_0(\overline{\Omega}, \mathbb{R}_+^m)$ converges to zero uniformly on $\overline\Omega$ as $t\to\infty$.
\end{theorem}

\medskip

The rest of this paper is organized as follows.  In Sections \ref{sec:asymptotically time-periodic}--\ref{sec:asymptotically fully degenerate}, we establish the
threshold  dynamics for system \eqref{intro-1} in three different asymptotic regimes, respectively.
In Section \ref{sec:applications}, we apply the developed theory  to a large class of reaction-diffusion-advection systems on time-varying domains.

\section{Asymptotically time-periodic regime}
\label{sec:asymptotically time-periodic}
In this section, we study the threshold dynamics in the fractional power space for a class of parabolic  systems
with asymptotically periodic coefficients.

\subsection{Asymptotically periodic semiflows}
Let $\Omega\subset\mathbb{R}^{n}$ be a bounded domain with boundary of class $C^{2+\theta}$ for some $0<\theta<1$.
Consider the following two nonautonomous parabolic systems:
\begin{equation}\label{nonautonomous asymptotically periodic equations}
\left\{
\begin{aligned}
& \frac{\partial v_i}{\partial t}   +A_i(t)v_i=f_i(x,t,v) \ \ \ &{\rm in}&\ \ \ \Omega\times(0,\infty),\\
& v_i =0 \ \ \ \ &{\rm on}&\ \ \ \partial\Omega\times(0,\infty),
\end{aligned}
\right.
\end{equation}
and
\begin{equation}\label{nonautonomous periodic equations}
\left\{
\begin{aligned}
& \frac{\partial v_i}{\partial t}+\tilde A_i(t)v_i=\tilde f_i(x,t,v) \ \ \ &{\rm in}&\ \ \ \Omega\times(0,\infty),\\
& v_i = 0 \ \ \ \ &{\rm on}&\ \ \ \partial\Omega\times(0,\infty),
\end{aligned}
\right.
\end{equation}
for $1\le i\le m$, where $v(x,t)=(v_1,\dots,v_m)^\top$, and
\begin{align*}
A_i(t)v_i=-\sum_{j,k=1}^{n}a^{(i)}_{jk}(x,t)\partial_{jk}v_i +\sum_{j=1}^{n}b^{(i)}_{j}(x,t)\partial_jv_i,
\\
\tilde A_i(t)v_i=-\sum_{j,k=1}^{n}\tilde a^{(i)}_{jk}(x,t)\partial_{jk}v_i +\sum_{j=1}^{n}\tilde b^{(i)}_{j}(x,t)\partial_jv_i.
\end{align*}
Throughout this section, we assume that
\begin{enumerate}[label=(A\arabic*),ref=(A\arabic*),leftmargin=*]
\item\label{ass:a1}
$a^{(i)}_{jk}$, $\tilde a^{(i)}_{jk}$, $b^{(i)}_{j}$, $\tilde b^{(i)}_{j}\in C^{\theta, \theta/2}(\overline\Omega\times\mathbb{R}_+, \mathbb{R})$,
$a^{(i)}_{jk}=a^{(i)}_{kj}$, $\tilde a^{(i)}_{jk}=\tilde a^{(i)}_{kj}$, 
and there exist positive constants $\lambda_i, \Lambda_i$ such that
\begin{equation}\label{uniform parabolicity condition}
\lambda_i|\xi|^2\le\sum_{j,k=1}^{n}a^{(i)}_{jk}(x,t)\xi_j\xi_k\le\Lambda_i|\xi|^2,
\end{equation}
for all $(x,t)\in \Omega\times\mathbb{R}_+$ and $\xi=(\xi_1,\dots,\xi_n)^\top$, $1\le i\le m$, $1\le j, k\le n$.

\item\label{ass:a2} 
For $1\le i,j\le m$,
$f_i,\tilde f_i,\frac{\partial f_i}{\partial u_j},
\frac{\partial\tilde f_i}{\partial u_j}
\in C(\overline{\Omega}\times\mathbb{R}_+\times\mathbb{R}^m,\mathbb{R})$,
and for any fixed $T>0$,
$f_i(\cdot,\cdot,u),\tilde f_i(\cdot,\cdot,u)
\in C^{\theta,\theta/2}(\overline{\Omega}\times[0,T],\mathbb{R})$
uniformly for $u$ in any bounded subset of $\mathbb{R}^m$.
Moreover, $\frac{\partial\tilde f_i}{\partial u_j}(\cdot,\cdot,0)
\in C^{\theta,\theta/2}(\overline{\Omega}\times[0,T],\mathbb{R})$
for every $T>0$.

\item\label{ass:a3} 
$\tilde A_i(t)$ and $\tilde f_i(\cdot,t,\cdot)$ are $\omega$-periodic in $\mathbb{R}$,
\begin{align*}
\lim_{t \to \infty}\left(\left\|a^{(i)}_{jk}(\cdot,t) - \tilde a^{(i)}_{jk}(\cdot,t) \right\|_{C(\overline\Omega)}
        +\left\| b^{(i)}_{j}(\cdot,t) - \tilde b^{(i)}_{j}(\cdot,t) \right\|_{C(\overline\Omega)}\right)= 0,
\end{align*}
and
\begin{align*}
    & \lim_{t \to \infty}\left\|f(\cdot,t,u) -\tilde f(\cdot,t,u) \right\|_{C(\overline\Omega)}
        = 0
\end{align*}
uniformly for all $u$ in any bounded subset of $\mathbb{R}^m$, $1\le i\le m$ and $1\le j,k\le n$, where
$\tilde f_u=\left(\frac{\partial\tilde f_i}{\partial u_{j}}\right)_{m\times m}$
and $f_u=\left(\frac{\partial f_i}{\partial u_{j}}\right)_{m\times m}$.
\end{enumerate}

Here and in what follows, we will occasionally use the notation
$\|\cdot\|_{S}$ to denote the $S$-norm of a matrix-valued function $M(x,t)=(m_{ij}(x,t))_{p\times q}$, $p,q\ge1$, defined by
$$
\|M\|_{S}=\max_{1\le i\le p}\sum_{j=1}^{q}\left\|m_{ij}\right\|_{S}.
$$
For $u,v\in\mathbb{R}^m$, we set $\|u\|_{\mathbb{R}^m}:=\max_{1\le i\le m}|u_i|$, and write $u\ge v$ if $u-v\in\mathbb{R}^m_+$, $u>v$ if $u-v\in\mathbb{R}^m_+\setminus\{0\}$, and $u\gg v$ if $u-v\in\mathrm{Int}(\mathbb{R}^m_+)$.

Let $X_{i}^{\beta}$ with $\beta\in[0,1]$ and $1 \le i \le m$ be the fractional power space of $L^p(\Omega)$ ($p>1$)
with respect to $(A_i(0), B_D)$ (e.g., see \cite{MR1100011}), and denote $X^{\beta}=\prod_{i=1}^m X_{i}^{\beta}$,
where $B_D$ means the homogeneous Dirichlet boundary condition on $\partial\Omega$.
Throughout this section, we choose
$p>n$ and $\beta\in(\frac{1}{2} + \frac{n}{2p}, 1)$
such that $X_{i}^{\beta}\hookrightarrow C^1_B(\overline\Omega)$, where $C_B^1(\overline{\Omega}) = \{u \in C^1(\overline{\Omega}) : u|_{\partial\Omega} = 0\}$. Thus,  $X_{i}^{\beta}$ is an ordered Banach space
with positive cone $P_i$ having nonempty interior $\mathrm{int}(P_i)$. Let $X^{\beta}=\prod_{i=1}^m X_{i}^{\beta}$
and $P = \prod_{i=1}^m P_i$. Then $(X^{\beta}, P)$ is an ordered Banach space.
We write $\|u\|_{\beta}=\max_{1\le i\le m}\|u_i\|_{X_i^\beta}$ for $u\in X^\beta$.
Moreover, let $\mathcal{L}(S_1, S_2)$ denote the space of bounded linear operators from $S_1$ to $S_2$.
When $S_1=S_2$, we simply write $\mathcal{L}(S_1)$ instead of $\mathcal{L}(S_1, S_2)$.
For brevity, the operator norm on $\mathcal{L}(X^{\alpha_1},X^{\alpha_2})$ is denoted by $\|\cdot\|_{\alpha_1,\alpha_2}$,
where $0\le \alpha_1,\alpha_2 \le 1$. In particular, $\|\cdot\|_{\mathcal{L}(L^p(\Omega)^m)}=\|\cdot\|_{0,0}$.
Denote
\begin{align*}
A(t)&=\mathrm{diag}(A_1(t),\dots,A_m(t)), \quad \tilde A(t)=\mathrm{diag}(\tilde A_1(t),\dots, \tilde A_m(t)),
\\
(A(t)-\tilde A(t))&\tilde A^{-1}(t)=\mathrm{diag}\left((A_1(t)-\tilde A_1(t))\tilde A^{-1}_1(t), \dots, (A_m(t)-\tilde A_m(t))\tilde A^{-1}_m(t)\right).
\end{align*}
Unless otherwise stated, we identify a diagonal operator matrix with its componentwise action on vector-valued functions.

\begin{lemma}\label{asymptotically periodic semiflows-lemma-1}
For every $t\ge 0$, $(A(t)-\tilde A(t))\tilde A^{-1}(t)\in\mathcal{L}(L^p(\Omega)^m)$. Moreover, we have
$\lim\limits_{t\to\infty}\left\|(A(t)-\tilde A(t))\tilde A^{-1}(t)\right\|_{0,0}=0$.
\end{lemma}

\begin{proof}
By the regularity theory of elliptic equations in Sobolev spaces, it follows that for any fixed $\varphi\in L^p(\Omega)^m$
and any $t\ge0$, there exists a unique function $u^{(t)}\in {\rm Dom}(\tilde A(t), B_D)\subset W^{2,p}(\Omega)^m$
such that $\tilde A(t)u^{(t)}=\varphi$, where $(\tilde A(t),B_D)=((\tilde A_1(t),B_D),\dots,(\tilde A_m(t),B_D))^\top$.
By \ref{ass:a3}, taking the limit in \eqref{uniform parabolicity condition} along $t+k\omega$ gives the uniform ellipticity for $\tilde A(t)$. By periodicity and \ref{ass:a1}, its coefficients have uniform H\"older bounds. Thus, we obtain (see, e.g., \cite[Chapter 9.1]{WuYinW2006})
$$
\|\tilde A^{-1}(t)\varphi\|_{W^{2,p}(\Omega)^m}=\|u^{(t)}\|_{W^{2,p}(\Omega)^m}\le c\|\varphi\|_{L^p(\Omega)^m},
$$
where $\tilde A^{-1}(t)\varphi=(\tilde A_1^{-1}(t)\varphi_1,\dots,\tilde A_m^{-1}(t)\varphi_m)^\top$, and $c>0$
is independent of both $t$ and $\varphi$. Throughout this section, $c$ denotes a generic positive constant,
which may vary from line to line. Clearly, for any fixed $t\ge0$,
$\tilde A^{-1}(t)\in\mathcal{L}(L^p(\Omega)^m, {\rm Dom}(\tilde A(t), B_D))$.
We then  have
\begin{align*}
    \|(A(t)-\tilde A(t))\tilde A^{-1}(t)\varphi\|_{L^p(\Omega)^m}
    & = \|(A(t)-\tilde A(t))u^{(t)}\|_{L^p(\Omega)^m}
    \\
	& \le \sum_{i=1}^{m}\sum_{j,k=1}^{n}\|a^{(i)}_{jk}(\cdot,t)-\tilde a^{(i)}_{jk}(\cdot,t)\|_{C(\overline\Omega)}\|\partial_{jk}u_i^{(t)}\|_{L^p(\Omega)}
\\
&  {\hspace{17pt}} +\sum_{i=1}^{m}\sum_{j=1}^{n}\|b^{(i)}_{j}(\cdot,t)-\tilde b^{(i)}_{j}(\cdot,t)\|_{C(\overline\Omega)}\|\partial_ju_i^{(t)}\|_{L^p(\Omega)}
\\
	&  \le  c\eta(t)\|\varphi\|_{L^p(\Omega)^m},
\end{align*}
where $c>0$ is independent of both $t$ and $\varphi$, and
\begin{align}\label{asymptotically periodic semiflows-lemma-1-eta}
\eta(t)=\sum_{i=1}^{m}\left(\sum_{j,k=1}^{n}\|a^{(i)}_{jk}(\cdot,t)-\tilde a^{(i)}_{jk}(\cdot,t)\|_{C(\overline\Omega)}+\sum_{j=1}^{n}\|b^{(i)}_{j}(\cdot,t)-\tilde b^{(i)}_{j}(\cdot,t)\|_{C(\overline\Omega)}\right).
\end{align}
Now the conclusion follows from the convergence assumptions on coefficients.
\end{proof}

We note that the domains ${\rm Dom}(A_i(t), B_D)$ and ${\rm Dom}(\tilde A_i(t), B_D)$ coincide and are independent of $t$ for each fixed $1 \le i \le m$.
Furthermore, by combining the uniform ellipticity of $A_i(t)$ and $\tilde A_i(t)$ with the regularity and asymptotic properties of their coefficients, one readily verifies that both $(A_i(t), B_D)$ and $(\tilde A_i(t), B_D)$ fulfill the structural assumptions given in \cite[Chapter I.2]{MR1204883} (see also \cite[Chapter 3.4.1]{MR2573296}) for all $1 \le i \le m$ and $t \ge 0$.
These assumptions ensure the existence of corresponding evolution operators on $\{(t,s); t\ge s\ge0\}$, denoted by $U_i(t,s)$ and $\widetilde U_i(t,s)$, respectively.
Denote
\begin{align*}
U(t,s)=\mathrm{diag}(U_1(t,s),\dots,U_m(t,s)), &\quad \widetilde U(t,s)=\mathrm{diag}(\widetilde U_1(t,s),\dots,\widetilde U_m(t,s)).
\end{align*}

The following lemma implies that the evolution operators depend continuously on the coefficients of the uniformly elliptic operator with respect to the
$\mathcal{L}(X^{\beta})$-norm.

\begin{lemma}\label{asymptotically periodic semiflows-lemma-2}
For any given $s_0>0$,  $\lim\limits_{s\to\infty}\left\|U(t+s,s)-\widetilde U(t+s,s)\right\|_{\beta,\beta}=0$
uniformly for all $t\in[0,s_0]$.
\end{lemma}

\begin{proof}
Set $\tilde f \equiv 0$ in \eqref{nonautonomous periodic equations} and take $\varphi\in X^{\beta}$.
For $s>0$ and $t\in[0,s_0]$, let $v(t+s,s,\varphi)$ with $v(s, s, \varphi) = \varphi$
be the unique solution of \eqref{nonautonomous periodic equations}.
Then \eqref{nonautonomous periodic equations} can be equivalently rewritten in the form
$\frac{\partial v}{\partial t}+A(t)v=A(t)v-\tilde A(t)v$.
By the variation of constants formula (see, e.g., \cite[Chapter I.2]{MR1204883}), we have
\begin{align*}
v(t+s,s,\varphi) & = \widetilde U(t+s,s)\varphi
    \\
    & = U(t+s,s)\varphi+\int_s^{t+s}U(t+s,\tau)\left(A(\tau)-\tilde A(\tau)\right)v\,\mathrm{d}\tau
    \\
    & = U(t+s,s)\varphi+\int_s^{t+s}U(t+s,\tau)\left(A(\tau)-\tilde A(\tau)\right)\widetilde U(\tau,s)\varphi\,\mathrm{d}\tau, \ t\in(0,s_0].
\end{align*}
By \cite[Lemma 5.2, Proposition 6.8]{MR1204883} (see also \cite[Chapter 3.8.1]{MR2573296}), we have the following two estimates:
\begin{align}\label{estimates of evolution operator-1}
\left\|U(t+s,\tau)\right\|_{0,\beta}&\le c_1(\beta)(t+s-\tau)^{-\beta},  {\hspace{-16mm}}&  s<\tau<t+s\le s_0+s,
\\[2mm]\label{estimates of evolution operator-2}
\left\|\widetilde U(\tau,s)\right\|_{\beta,1}&\le c_2(\beta)(\tau-s)^{\beta-1},  {\hspace{-6mm}}&s<\tau<t+s\le s_0+s,
\end{align}
where the positive constants $c_1(\beta), c_2(\beta)$ are independent of $s$ and $t$.
Then for $t\in(0, s_0]$, using $\beta\in(1/2 + n/(2p), 1)$, we deduce
\begin{align*}
    &\left\|U(t+s,s)\varphi-\widetilde U(t+s,s)\varphi\right\|_{\beta}
    \\
    & \le \int_s^{t+s}\left\|U(t+s,\tau)\right\|_{0,\beta}\left\|\left(A(\tau)-\tilde A(\tau)\right)\tilde A^{-1}(\tau)\right\|_{0,0}\left\|\tilde A(\tau)\widetilde U(\tau,s)\varphi\right\|_{0}\mathrm{d}\tau
    \\[2mm]
    & \le c\int_s^{t+s}\left\|U(t+s,\tau)\right\|_{0,\beta}\left\|\left(A(\tau)-\tilde A(\tau)\right)\tilde A^{-1}(\tau)\right\|_{0,0}\left\|\widetilde U(\tau,s)\right\|_{\beta,1}\|\varphi\|_{\beta}\,\mathrm{d}\tau
    \\
    & \le c\|\varphi\|_{\beta}\sup_{\tau\in[0,s_0]}\eta(s+\tau)\int_s^{t+s}(t+s-\tau)^{-\beta}(\tau-s)^{\beta-1}\,\mathrm{d}\tau
    \\
    &  =  c\|\varphi\|_{\beta}\ \frac{\pi}{\sin (\beta\pi)}\,\sup_{\tau\in[0,s_0]}\eta(s+\tau), \quad \forall t\in(0, s_0],
\end{align*}
where $\eta(\tau)$ is given by \eqref{asymptotically periodic semiflows-lemma-1-eta}, and $c>0$ is independent of $\varphi$, $s$ and $t$.
Combining with the case $t=0$, we have
\begin{align}\label{asymptotically periodic semiflows-lemma-2-1}
\left\|U(t+s,s)\varphi-\widetilde U(t+s,s)\varphi\right\|_{\beta}\le c\,\frac{\pi \|\varphi\|_{\beta}}{\sin (\beta\pi)}\sup_{\tau\in[0,s_0]}\eta(s+\tau)
\end{align}
uniformly for $t\in[0,s_0]$, where $c>0$ is independent of $\varphi$, $s$ and $t$. Letting $s\to\infty$ in the above inequality, we see that the desired conclusion
holds true.
\end{proof}

For $t\ge s\ge0$ and any $\varphi\in X^{\beta}$, let $\phi(t, s, \varphi)$ and $\tilde\phi(t, s, \varphi)$ denote the solutions of
systems \eqref{nonautonomous asymptotically periodic equations} and \eqref{nonautonomous periodic equations} respectively, with
$\phi(s, s, \varphi) =\tilde\phi(s, s, \varphi)= \varphi$. By the results in \cite[Chapter~III.20]{MR1100011},
these solutions are unique classical solutions defined on their maximal intervals of existence.
Moreover, they exist globally on $[s,\infty)$ provided that they remain bounded in $L^\infty(\Omega)^m$ on finite time intervals.

\begin{proposition}\label{asymptotically periodic semiflows-lemma-3}
Assume that $\phi(t, s, \varphi)$ and $\tilde\phi(t, s, \varphi)$ are uniformly bounded in $X^{\beta}$ for all $0\le s\le t<\infty$
and $\varphi$ in any bounded subset of $X^{\beta}$.
Then for any positive integer $k_0$, we have
\begin{equation}\label{asymptotically periodic semiflows-lemma-3-1}
\lim_{k\to\infty}\sup_{t\in[k\omega,(k+k_0)\omega]} \left\|\phi(t, k\omega, \varphi)-\tilde\phi(t, k\omega, \varphi)\right\|_{\beta}=0
\end{equation}
uniformly for all $\varphi$ in any bounded set of $X^{\beta}$.

In particular, for any $\varphi\in X^{\beta}$, $\gamma^+(\varphi)=\{\phi(k\omega,0,\varphi); k\ge0\}$ is precompact in $X^{\beta}$,
and $\phi: \Delta_0\times X^{\beta}\to X^{\beta}$ is an asymptotically periodic semiflow with limit periodic semiflow $T(t): X^{\beta}\to X^{\beta}$, $t\ge0$,
where $\Delta_0:=\{(t,s): 0\le s\le t<\infty\}$.
\end{proposition}

In  the case where $\tilde a^{(i)}_{jk}, \tilde b^{(i)}_{j}, \tilde f$ are independent of $t$ for $1\le i\le m$, $1\le j, k\le n$,
Proposition \ref{asymptotically periodic semiflows-lemma-3} implies that  $\phi(t,s,\varphi)$ is an asymptotically autonomous semiflow.
Unless otherwise specified, $X^\beta$-norms appearing here are understood to be taken in the spatial variable.
For simplicity, we suppress the spatial variable when no confusion is possible.

\begin{proof}
Our proof is motivated by  \cite[Proposition 3.2.1]{MR3643081}.
Fix a bounded set $B\subset X^\beta$ and a positive integer $k_0$.
Let $M$ be a common $X^\beta$-bound for the two solutions
for all $0\le s\le t$ and $\varphi\in B$.
By the variation of constants formula, for $t\in[k\omega,(k+k_0)\omega]$ with $k\in\mathbb{N}$, we have
\begin{align*}
\phi(t):=\phi(t, k\omega, \varphi)=U(t,k\omega)\varphi+\int_{k\omega}^{t}U(t,\tau)f(x,\tau,\phi(\tau))\,\mathrm{d}\tau,
\\
\tilde\phi(t):=\tilde\phi(t, k\omega, \varphi)=\widetilde U(t,k\omega)\varphi+\int_{k\omega}^{t}\widetilde U(t,\tau)\tilde f(x,\tau,\tilde\phi(\tau))\,\mathrm{d}\tau,
\end{align*}
where $f=(f_1,\dots,f_m)^\top$ and $\tilde f=(\tilde f_1,\dots,\tilde f_m)^\top$.
It then follows that
\begin{align*}
    \left\|\phi(t)-\tilde\phi(t)\right\|_{\beta} & \le \left\|U(t,k\omega)\varphi-\widetilde U(t,k\omega)\varphi\right\|_{\beta}+\int_{k\omega}^{t}\left\|\left(U(t,\tau)-\widetilde U(t,\tau)\right)\tilde f(\cdot,\tau,\tilde\phi(\tau))\right\|_{\beta}\,\mathrm{d}\tau
    \\
    & {\hspace{17pt}}  +\int_{k\omega}^{t}\left\|U(t,\tau)\left(f(\cdot,\tau,\phi(\tau))-\tilde f(\cdot,\tau,\tilde\phi(\tau))\right)\right\|_{\beta}\,\mathrm{d}\tau
    \\
    & = I_1+I_2+I_3.
\end{align*}
In view of  \eqref{asymptotically periodic semiflows-lemma-2-1}, there holds
\begin{align*}
    I_1 & =\left\|U(t,k\omega)\varphi-\widetilde U(t,k\omega)\varphi\right\|_{\beta}
    \le c\,\frac{\pi \|\varphi\|_{\beta}}{\sin (\beta\pi)}\sup_{\tau\in[k\omega,(k+k_0)\omega]}\eta(\tau)
\end{align*}
for all $t\in[k\omega,(k+k_0)\omega]$, where $c>0$ is independent of $t,k$.

For $I_2$, the compact embedding $X^\beta\hookrightarrow C(\bar\Omega)^m$ and
continuity of $\tilde f$ imply that
\[
\mathcal G=\bigl\{x\mapsto\tilde f(x,\sigma,z(x)):
0\le\sigma\le\omega,\ z\in X^\beta,\ \|z\|_\beta\le M\bigr\}
\]
is relatively compact in $L^p(\Omega)^m$. For any $\varepsilon>0$, choose
$h_1,\ldots,h_N\in C_c^\infty(\Omega)^m$ such that
\[
\sup_{g\in\mathcal G}\min_{1\le j\le N}\|g-h_j\|_0<\varepsilon.
\]
Since $\tilde f(\cdot,\tau,\tilde\phi(\tau))\in\mathcal G$, choose $j$ with
$\|\tilde f(\cdot,\tau,\tilde\phi(\tau))-h_j\|_0<\varepsilon$.
By \eqref{estimates of evolution operator-1} and \eqref{asymptotically periodic semiflows-lemma-2-1}, we have
\begin{align*}
	&\left\|\big(U(t,\tau)-\widetilde U(t,\tau)\big)\tilde f(\cdot,\tau,\tilde\phi(\tau))\right\|_\beta
	\\
	&\le c\,\frac{\pi\|h_j\|_\beta}{\sin (\beta\pi)}\sup_{\tau\in[k\omega,(k+k_0)\omega]}\eta(\tau)+C(t-\tau)^{-\beta}
	\left\|\tilde f(\cdot,\tau,\tilde\phi(\tau))-h_j\right\|_0.
\end{align*}
Then
\begin{align*}
	I_2\le C_\varepsilon\sup_{\tau\in[k\omega,(k+k_0)\omega]}\eta(\tau)+c\varepsilon
\end{align*}
for all $t\in[k\omega,(k+k_0)\omega]$, where $C_\varepsilon>0$ depends on the choice of $h_1,\ldots,h_N$, and $c>0$ is independent of $t,k,\varepsilon$.

For $I_3$, choose $R>0$ such that $\|z\|_{C(\bar\Omega)^m}\le R$
whenever $\|z\|_\beta\le M$.
By \ref{ass:a2}  and the $\omega$-periodicity of $\tilde f$, $\tilde f_u$ is bounded on $\bar\Omega\times[0,\infty)\times[-R,R]^m$.
Thus $\tilde f(x,\tau,\cdot)$ is Lipschitz on $[-R,R]^m$ with a
constant independent of $x$ and $\tau$. By the estimate \eqref{estimates of evolution operator-1}, we obtain
\begin{align*}
	I_3 & \le \int_{k\omega}^{t}\left\| U(t,\tau)\right\|_{0,\beta}\left\|\left(f(\cdot,\tau,\phi(\tau))-\tilde f(\cdot,\tau,\phi(\tau))\right)\right\|_{0}\,\mathrm{d}\tau
	\\
	& {\hspace{17pt}} +\int_{k\omega}^{t}\left\| U(t,\tau)\right\|_{0,\beta}\left\|\left(\tilde f(\cdot,\tau,\phi(\tau))-\tilde f(\cdot,\tau,\tilde\phi(\tau))\right)\right\|_{0}\,\mathrm{d}\tau
	\\
	& \le
	c\frac{(k_0\omega)^{1-\beta}}{1-\beta}\|f-\tilde f\|_{C(D_k)}+ c\int_{k\omega}^{t}(t-\tau)^{-\beta}\left\|\phi(\tau)-\tilde\phi(\tau)\right\|_{\beta}\,\mathrm{d}\tau,
\end{align*}
where $t\in[k\omega,(k+k_0)\omega]$ and $D_k:=\overline\Omega\times[k\omega, (k+k_0)\omega]\times[-R,R]^m\subset\mathbb{R}^{n}\times\mathbb{R}_+\times\mathbb{R}^m$.

Summarizing the estimates of $I_1,I_2,I_3$, we derive that
\begin{align*}
    &\left\|\phi(t)-\tilde\phi(t)\right\|_{\beta}
    \\
    & \le c\varepsilon+C_\varepsilon\sup_{\tau\in[k\omega,(k+k_0)\omega]}\eta(\tau)+c\|f-\tilde f\|_{C(D_k)}+c\int_{k\omega}^{t}(t-\tau)^{-\beta}\left\|\phi(\tau)-\tilde\phi(\tau)\right\|_{\beta}\,\mathrm{d}\tau,
\end{align*}
where $c>0$ is independent of $t,k,\varepsilon$ and $\varphi\in B$.
By a generalized Gronwall inequality (\cite[Lemma 19.4]{MR1100011}), we obtain
\begin{align*}
\left\|\phi(t)-\tilde\phi(t)\right\|_{\beta}\le c\varepsilon+C_\varepsilon\sup_{\tau\in[k\omega,(k+k_0)\omega]}\eta(\tau)+c\|f-\tilde f\|_{C(D_k)}, \quad \forall t\in[k\omega,(k+k_0)\omega].
\end{align*}
Passing to the limit first as $k\to\infty$ (with $\varepsilon$ fixed) and then as $\varepsilon\to 0$ yields \eqref{asymptotically periodic semiflows-lemma-3-1}.

Denote by $S: X^{\beta}\to X^{\beta}$ the Poincar\'e map associated with the periodic system \eqref{nonautonomous periodic equations}.
For each $\varphi\in X^{\beta}$ and $k\in\mathbb{N}$, define
\begin{align*}
T_k(\varphi):= \phi(k\omega,0,\varphi), \quad S_k(\varphi):= \phi((k+1)\omega,k\omega,\varphi).
\end{align*}
Then \eqref{asymptotically periodic semiflows-lemma-3-1} yields $\lim_{k\to\infty}\|S_k(\varphi)-S(\varphi)\|_{\beta}=0$
uniformly for $\varphi$ in any bounded set of $X^{\beta}$. By the uniform boundedness of $T_k(\varphi)$ in $X^{\beta}$
for fixed $\varphi\in X^{\beta}$, we see that
$$
\lim_{k\to\infty}\left\|T_{k+1}(\varphi)-S(T_k(\varphi))\right\|_{\beta}=\lim_{k\to\infty}\left\|S_{k}(T_k(\varphi))-S(T_k(\varphi))\right\|_{\beta}=0.
$$
Since $S$ is continuous and compact, $S(\gamma^+(\varphi))$ is precompact in $X^{\beta}$.
Therefore, $\gamma^+(\varphi)$ is precompact in $X^{\beta}$.

For any fixed $(t_0,\varphi_0)\in \mathbb{R}_+\times X^{\beta}$, we have
\begin{align*}
    &\left\|\phi(t+k\omega,k\omega,\varphi)-T(t_0)\varphi_0\right\|_{\beta}
    \\
    & \le \left\|\phi(t+k\omega,k\omega,\varphi)-T(t)\varphi\right\|_{\beta}+\left\|T(t)\varphi-T(t_0)\varphi_0\right\|_{\beta}
    \\
    &  =  \left\|\phi(t+k\omega,k\omega,\varphi)-\tilde\phi(t+k\omega, k\omega, \varphi)\right\|_{\beta}+\left\|T(t)\varphi-T(t_0)\varphi_0\right\|_{\beta}.
\end{align*}
By \eqref{asymptotically periodic semiflows-lemma-3-1} and the continuity of $T(t)\varphi$ for $(t,\varphi)\in\mathbb{R}_+\times X^{\beta}$,
it then follows that
$$
\lim_{(t,\varphi,k)\to(t_0,\varphi_0,\infty)}\left\|\phi(t+k\omega,k\omega,\varphi)-T(t_0)\varphi_0\right\|_{\beta}=0.
$$
Thus, $\phi(t,s,\varphi): \Delta_0\times X^{\beta}\to X^{\beta}$ is asymptotic to the $\omega$-periodic semiflow $T(t): X^{\beta}\to X^{\beta}$, $t\ge0$.
\end{proof}

We remark that
Proposition \ref{asymptotically periodic semiflows-lemma-3} includes \cite[Theorem A.1]{FanZhao2026} as a special case within a more general framework of uniformly parabolic systems.
The argument in \cite[Theorem A.1]{FanZhao2026} relies on the connection between two-parameter evolution operators and analytic semigroups,
and cannot be applied to systems involving advection terms and anisotropic diffusion.

\subsection{Threshold dynamics}

For a given $\bar{u} \gg 0$, let $\mathbb{M}_{\bar{u}}:= [0, \bar{u}]_{\mathbb{R}^m}= \{\, u \in \mathbb{R}^m : 0 \le u \le \bar{u} \,\}$.
To study the global dynamics of system \eqref{nonautonomous asymptotically periodic equations}, we make the following additional assumptions on $\tilde f$ for all $(x,t)\in\overline\Omega\times\mathbb{R}$:
\begin{enumerate}[label=(A\arabic*),ref=(A\arabic*),leftmargin=*,start=4]
\item\label{ass:a4} 
$\tilde f_u(x, t,u)$ is cooperative and irreducible for $(x,t,u)\in\overline\Omega\times\mathbb{R}\times \mathbb{M}_{\bar{u}}$.
\item\label{ass:a5} 
$\tilde f(\cdot,\cdot,u)$ is strictly subhomogeneous on $\mathbb{M}_{\bar{u}}$ in the sense that
$\alpha\tilde f(\cdot,\cdot,u)<\tilde f(\cdot,\cdot,\alpha u)$ for all  $u\in\mathbb{M}_{\bar{u}}\cap\mathrm{Int}(\mathbb{R}^{m}_{+})$ and $\alpha\in (0,1)$.
\end{enumerate}

We first introduce an auxiliary vector-valued barrier function $\Phi$, which will be used in our construction of subsolutions for system \eqref{nonautonomous asymptotically periodic equations}.

\begin{lemma}\label{Sec2 auxiliary functions}
There exist $\mu>0$ and a vector-valued function $\Phi=(\Phi_1,\dots,\Phi_m)^\top$ with $\Phi_i\in C^{2}(\overline\Omega)$ for all  $1\leq i\leq m$  such that
\begin{equation}\label{Sec2 auxiliary functions-1}
\left\{
\begin{aligned}
& (A(t)+\mu)\Phi\ge (1,\dots,1)^\top \ \ \ &{\rm in}&\ \ \ \Omega\times(0,\infty),\\
& \Phi\gg0 \ \ \ \ &{\rm in}&\ \ \ \Omega,\\
& \Phi=0 \ \ \ \ &{\rm on}&\ \ \ \partial\Omega,
\end{aligned}
\right.
\end{equation}
where $(A(t)+\mu)\Phi:=((A_1(t)+\mu)\Phi_1,\dots,(A_m(t)+\mu)\Phi_m)^\top$.
\end{lemma}

\begin{proof}
For any given symmetric $n\times n$ matrix $M$ and  each $1\le i\le m$,
we introduce Pucci's extremal operators with parameters $0<\lambda_i<\Lambda_i$
as follows:
\begin{align*}
\mathcal{M}^+_{\lambda_i,\Lambda_i}(M)=\sup_{A\in\mathfrak{A}_i}{\rm tr}(AM), \quad \mathcal{M}^-_{\lambda_i,\Lambda_i}(M)=\inf_{A\in\mathfrak{A}_i}{\rm tr}(AM),
\end{align*}
where $\mathfrak{A}_i$ denotes the set of all symmetric matrices whose eigenvalues lie in the interval $[\lambda_i,\Lambda_i]$.
For basic properties of Pucci's extremal operators, we refer to \cite[Section 2.2]{Caffarelli1995}.
By the uniform ellipticity of $A_i(t)$, we have
$$
\sum_{j,k=1}^{n}a^{(i)}_{jk}(x,t)\partial_{jk}v_i={\rm tr}(\mathcal{A}_i\mathcal{D}^2v_i)\le \mathcal{M}^+_{\lambda_i,\Lambda_i}(\mathcal{D}^2v_i),
$$
where $v_i\in C^{2}(\Omega, \mathbb{R})$, $\mathcal{A}_i=(a^{(i)}_{jk})_{n\times n}$ and $\mathcal{D}^2$ denotes
Hessian matrix. It follows that if $v_i\in C^{2}(\Omega,\mathbb{R}_+)$ for some $1\le i\le m$, then
\begin{equation}\label{Sec2 auxiliary functions-2}
-A_i(t)v_i\le\mathcal{M}^+_{\lambda_i,\Lambda_i}(\mathcal{D}^2v_i)+\bar b^{(i)}|\nabla v_i|,
\end{equation}
where $\bar b^{(i)}=\sqrt{n}\sup_{1\le j\le n}\|b^{(i)}_{j}\|_{C(\overline\Omega\times[0,\infty))}$.
For each $i\in\{1,\dots,m\}$, we consider the following elliptic problem:
\begin{equation}\label{Sec2 auxiliary functions-3}
\left\{
\begin{aligned}
& -\mathcal{M}^+_{\lambda_i,\Lambda_i}(\mathcal{D}^2\Phi_i)-\bar b^{(i)}|\nabla \Phi_i|+\mu\Phi_i=1 \ \ \ &{\rm in}&\ \ \ \Omega,\\
& \Phi_i=0 \ \ \ \ &{\rm on}&\ \ \ \partial\Omega,
\end{aligned}
\right.
\end{equation}
where $\mu>0$.
Denote
$$
H_i(\mathcal{D}^2u,\mathcal{D}u,u,x)=-\mathcal{M}^+_{\lambda_i,\Lambda_i}(\mathcal{D}^2u)-\bar b^{(i)}|\nabla u|+\mu u
$$
for all $u\in C^2(\overline\Omega, \mathbb{R})$ and $x\in\overline\Omega$.
Note that the operator $H_i$ is proper, positively homogeneous of order 1, and uniformly elliptic, where {\it proper} means that $H_i$ is nondecreasing in $u$.
By \cite[Theorem 2.3(i)]{Armstrong2009} and \cite[Theorem 77]{MR3488816},
there exists sufficiently large $\mu$ such that the fully nonlinear elliptic problem~\eqref{Sec2 auxiliary functions-3} admits a unique nonnegative
viscosity solution $\Phi_i\in C^{1}(\overline\Omega)$ for each $1\le i\le m$.
By the strong maximum principle and Hopf boundary lemma~\cite[Theorem~A.1]{Armstrong2009}, it follows that $\Phi_i>0$ in $\Omega$ for each $1\le i\le m$.
Moreover, the regularity results~\cite[Theorem 1, Theorem 78]{MR3488816} and \cite[Theorem 1.2]{Krylov1984}
imply that
$\Phi_i\in C^{2}(\overline\Omega)$ for each $1\le i\le m$.
Combining with \eqref{Sec2 auxiliary functions-2}--\eqref{Sec2 auxiliary functions-3}, we complete the proof.
\end{proof}

Let $(\lambda^*_\omega, \phi^*_\omega)$ be the principal eigen-pair of the following $\omega$-periodic parabolic eigenvalue problem:
\begin{equation}\label{eigenvalue-nonautonomous periodic equations}
\left\{
\begin{aligned}
& \frac{\partial v}{\partial t}+\tilde A(t)v=\tilde f_u(x,t,0)v+\lambda v \ \ \ &{\rm in}&\ \ \ \Omega\times\mathbb{R},\\
& v=0 \ \ \ \ &{\rm on}&\ \ \ \partial\Omega\times\mathbb{R},\\
& v(x,\omega+t)=v(x,t) \ \ \ &{\rm in}&\ \ \ \Omega\times\mathbb{R}.
\end{aligned}
\right.
\end{equation}
Then $\phi^*_\omega\in C^{2,1}(\overline\Omega\times\mathbb{R},\mathbb{R}^m)$ and $\phi^*_\omega\gg0$ in $\Omega\times\mathbb{R}$. By Hopf boundary lemma, there exists $c_*\in(0,\mu)$ such that
\begin{equation}\label{phi-Phi}
\phi^*_\omega\ge c_*\Phi \ \ \ {\rm on}\ \ \ \overline\Omega\times[0,\omega].
\end{equation}

In order to construct an appropriate  subsolution of system \eqref{nonautonomous asymptotically periodic equations}  in the case where $\lambda_\omega^{\ast}<0$, we  need the following analogous assumptions on $f$ for all $(x,t)\in\overline\Omega\times\mathbb{R}_+$:
\begin{enumerate}[label=(A\arabic*),ref=(A\arabic*),leftmargin=*,start=6]
\item\label{ass:a6}
$f_u(x, t,u)$ is cooperative and irreducible for $(x,t,u)\in\overline\Omega\times\mathbb{R}_+\times \mathbb{M}_{\bar{u}}$.

\item\label{ass:a7} 
For each $1 \le i \le m$, $ f_i(\cdot,\cdot,u)\ge0$ whenever $u \in \mathbb{M}_{\bar{u}}$ with $u_i = 0$,
and $f_i(\cdot,\cdot,u)\le 0$ whenever $u\in\mathbb{M}_{\bar{u}}$ with $u_i = \bar{u}_i$;
$ f(\cdot,\cdot,0) = 0$;  and $f(\cdot,\cdot,u)$ is subhomogeneous on $\mathbb{M}_{\bar{u}}$ in the sense that
$\alpha f(\cdot,\cdot,u)\le f(\cdot,\cdot,\alpha u)$ for all  $u\in\mathbb{M}_{\bar{u}}\cap\mathrm{Int}(\mathbb{R}^{m}_{+})$ and $\alpha\in (0,1)$.
\end{enumerate}

The following proposition shows that whenever $\lambda_\omega^{\ast}<0$,
there exists a strongly positive periodic function $\underline v_{\omega}$ such that every nontrivial solution of system \eqref{nonautonomous asymptotically periodic equations}
is eventually bounded below by $\underline v_{\omega}$, and hence remains strongly separated from zero as $t\to\infty$.

\begin{proposition}\label{lowersolution-for-general-threshold-result}
Assume that $\lambda_\omega^{\ast}<0$. Then there exist positive constants $\eta, \varepsilon, t_0$ such that the
$\omega$-periodic function $\underline v_{\omega}(x,t):=\eta\big(\phi^*_\omega(x,t)-\varepsilon\Phi(x)\big)$
satisfies $\underline v_{\omega}\gg0$ in $\Omega\times[0,\omega]$, $\underline v_{\omega}=0$ on $\partial\Omega\times[0,\omega]$,
and is a subsolution for system \eqref{nonautonomous asymptotically periodic equations} in $\Omega\times[t_0,\infty)$.
\end{proposition}

\begin{proof}
By Hopf boundary lemma, there exists $C_*>0$ such that
$$
\max_{1\le j\le m}(\phi_\omega^*)_j\le C_*\,(\phi_\omega^*)_i \quad {\rm in} \ \  \overline\Omega\times[0,\omega], \ \ 1\le i\le m.
$$
Take $\varepsilon$ such that
$0<\varepsilon<\min\{\frac{c_*}{2}, \frac{-\lambda_\omega^{\ast}c_*}{8\mu}, \frac{-\lambda_\omega^*c_*}{8C_*\|\tilde f_u(\cdot,\cdot,0)\|_{C(\overline\Omega\times[0,\omega])}}\}$.
Then there exists a sufficiently large $t_0>0$ such that
\begin{align}\label{lowersolution-for-general-threshold-result-2}
\big(A(t)-\tilde A(t)\big)\phi^*_\omega(x,t)\le \varepsilon(1,\dots,1)^\top \quad {\rm in} \ \ \overline\Omega\times[t_0,\infty).
\end{align}
By \ref{ass:a2} and time periodicity, $\tilde f(x,t,\cdot)$ is differentiable at zero
uniformly in $(x,t)\in\overline\Omega\times\mathbb R$.
Since $\tilde f(x,t,0)=0$ by \ref{ass:a3}, \ref{ass:a7} and periodicity,
there exists $r_*\in(0,1)$ such that, for all $x\in\overline\Omega$ and $t\in\mathbb R$,
\[
\tilde f(x,t,z)\ge \tilde f_u(x,t,0)z
-\frac{\varepsilon}{2C_*}\|z\|_{\mathbb R^m}(1,\dots,1)^\top,
\qquad z\ge0,\quad \|z\|_{\mathbb R^m}\le r_*.
\]
Decreasing $r_*$ if necessary, we may also assume that $z\in\mathbb M_{\bar u}$
whenever $z\ge0$ and $\|z\|_{\mathbb R^m}\le r_*$.
With $r_*$ fixed, applying \ref{ass:a3} on the sphere $\|z\|_{\mathbb R^m}=r_*$ ($z\ge 0$) and then scaling down to $0<\|z\|_{\mathbb R^m}\le r_*$ ($z\gg 0$) via the subhomogeneity in \ref{ass:a7}, we obtain, after increasing $t_0$ if necessary, that for all $x\in\overline\Omega$ and $t\ge t_0$,
\begin{equation}\label{lowersolution-for-general-threshold-result-2-1}
	f(x,t,z)\ge \tilde f_u(x,t,0)z
	-\frac{\varepsilon}{C_*}\|z\|_{\mathbb R^m}(1,\dots,1)^\top,
	\qquad z\ge0,\quad \|z\|_{\mathbb R^m}\le r_*.
\end{equation}
Here the inequality extends from $z\gg0$ to the boundary of the positive cone
by continuity, with equality at $z=0$.

Let $\eta=\frac{r_*}{\|\phi^*_\omega\|_{C(\overline\Omega\times[0,\omega])}}$ and $\delta=\eta\varepsilon$.
Combining with \eqref{phi-Phi}, we have
\begin{equation}\label{lowersolution-for-general-threshold-result-3}
\frac{c_*}{2}\eta \Phi(x)\ll \underline v_{\omega}(x,t)\ll\eta\phi^*_\omega(x,t)\le r_*(1,\dots,1)^\top \quad {\rm in} \ \ \Omega\times[t_0,\infty).
\end{equation}
A straightforward calculation yields
\begin{align}\label{lowersolution-for-general-threshold-result-4}
    &\frac{\partial \underline v_{\omega}}{\partial t}+A(t)\underline v_{\omega}-f(x,t,\underline v_{\omega})      \nonumber
    \\      \nonumber
    &  =  \eta\lambda_\omega^*\phi^*_\omega+\eta\big(A(t)-\tilde A(t)\big)\phi^*_\omega-\delta(A(t)+\mu)\Phi+\delta\mu\Phi
    \\
    & {\hspace{17pt}} +\delta\tilde f_u(x,t,0)\Phi-\big(f(x,t,\underline v_{\omega})-\tilde f_u(x,t,0)\underline v_{\omega}\big).
\end{align}
It then follows from \eqref{Sec2 auxiliary functions-1} and \eqref{lowersolution-for-general-threshold-result-2} that
\begin{align}\label{lowersolution-for-general-threshold-result-5}
\eta\big(A(t)-\tilde A(t)\big)\phi^*_\omega-\delta(A(t)+\mu)\Phi
\le \eta\varepsilon(1,\dots,1)^\top-\delta (1,\dots,1)^\top = 0.
\end{align}
In view of \eqref{phi-Phi} and the choices of $C_*,\varepsilon$, we obtain
\begin{align}\label{lowersolution-for-general-threshold-result-6}
    \delta\mu\Phi+\delta\tilde f_u(x,t,0)\Phi
    &\le \frac{\delta\mu}{c_*}\phi^*_\omega+\frac{\delta}{c_*}C_*\big\|\tilde f_u(\cdot,\cdot,0)\big\|_{C(\overline\Omega\times[0,\omega])}\phi_\omega^* 
    \le -\frac{\lambda_\omega^*}{4}\eta\phi_\omega^*.
\end{align}
By \eqref{lowersolution-for-general-threshold-result-2-1}, \eqref{lowersolution-for-general-threshold-result-3} and the choices of $C_*$, for any fixed $(x,t)\in\overline\Omega\times[t_0,\infty)$, we have
\begin{align}\label{lowersolution-for-general-threshold-result-7}
	-\big(f(x,t,\underline v_\omega)-\tilde f_u(x,t,0)\underline v_\omega\big)
	&\le \frac{\varepsilon}{C_*}\|\underline v_\omega\|_{\mathbb R^m}(1,\dots,1)^\top
	\le -\frac{\lambda_\omega^*}{4}\eta\phi_\omega^*.
\end{align}
Substituting estimates \eqref{lowersolution-for-general-threshold-result-5}--\eqref{lowersolution-for-general-threshold-result-7}
into \eqref{lowersolution-for-general-threshold-result-4} yields that
\begin{align*}
\frac{\partial \underline v_{\omega}}{\partial t}+A(t)\underline v_{\omega}-f(x,t,\underline v_{\omega})\ll0,\quad (x,t)\in\Omega\times[t_0,\infty).
\end{align*}
This completes the proof.
\end{proof}

Now we are ready to prove the following threshold type result on the global dynamics of system \eqref{nonautonomous asymptotically periodic equations}.

\begin{theorem}\label{general threshold dynamics}
Let	$\mathcal{M}_{\bar u}^{\beta}:=\{\varphi\in X^{\beta}:\, \varphi(x)\in\mathbb{M}_{\bar{u}},\ \forall x\in\overline\Omega\}$.
Then the  following statements are valid:
	\begin{itemize}
		\item [(a)] If $\lambda_\omega^{\ast}\ge0$, then the zero solution is globally asymptotically stable for \eqref{nonautonomous asymptotically periodic equations} in $\mathcal{M}_{\bar u}^{\beta}$.
		\item [(b)] If $\lambda_\omega^{\ast}<0$, then every nontrivial solution $\phi$ of  \eqref{nonautonomous asymptotically periodic equations}  with $\phi(0)\in\mathcal{M}_{\bar u}^{\beta}$
                    satisfies $\lim_{t\to+ \infty} \|\phi(t)-v_\omega^*(t)\|_{\beta}=0$, where $v_\omega^*$ is the unique positive $\omega$-periodic solution in $\mathcal{M}_{\bar u}^{\beta}$ of \eqref{nonautonomous periodic equations}.
	\end{itemize}
\end{theorem}

\begin{proof}
For $t\ge s\ge0$ and any $\varphi\in\mathcal{M}_{\bar u}^{\beta}$, let $\phi(t, s, \varphi)$ with $\phi(s, s, \varphi) = \varphi$, and $\tilde\phi(t, s, \varphi)$ with $\tilde\phi(s, s, \varphi) = \varphi$
be the unique solutions of \eqref{nonautonomous asymptotically periodic equations} and \eqref{nonautonomous periodic equations}, respectively.
It is easy to see that $\mathcal{M}_{\bar u}^{\beta}$ is positively invariant for both systems by comparison. Then for any $s \ge 0$, $\phi(t,s,\varphi)$ and $\tilde\phi(t,s,\varphi)$ exist globally on $[s,\infty)$ and are
uniformly bounded in $X^{\beta}$. Let $T_k(\varphi):= \phi(k\omega,0,\varphi)$, $k\ge 0$.
By Proposition \ref{asymptotically periodic semiflows-lemma-3}, it follows that
$\gamma^+(\varphi) = \{T_k(\varphi): k\ge 0\}$ is precompact in $X^{\beta}$, and its omega limit set $\omega(\varphi)$ exists.
Note that
$T_k: \mathcal{M}_{\bar u}^{\beta}\to\mathcal{M}_{\bar u}^{\beta}, k\ge 0$, is an asymptotically autonomous discrete process with limit discrete semiflow
$S^k: \mathcal{M}_{\bar u}^{\beta}\to\mathcal{M}_{\bar u}^{\beta}, k\ge 0$, where $S = T(\omega)$ is the Poincar\'e map associated with the periodic system \eqref{nonautonomous periodic equations}.
In view of \cite[Lemma 1.2.2]{MR3643081},  it then follows that  for any $\varphi\in\mathcal{M}_{\bar u}^{\beta}$, $\omega(\varphi)$ is an internally chain transitive set for $S$.

Let $\lambda_\omega^{\ast}\ge0$. By \cite[Theorem 2.2.2]{MR3643081}, we find that $0$ is a globally asymptotically stable fixed point of $S$,
and hence $W^s(0)=\mathcal{M}_{\bar u}^{\beta}$, where $W^s(0)$ is the stable set of $0$ for $S$. Then we have $\omega(\varphi)\cap W^s(0)\ne\emptyset$.
By \cite[Theorem 1.2.1]{MR3643081}, $\omega(\varphi)=\{0\}$, which implies that $\lim_{k\to \infty} T_k(\varphi) = 0$ for any $\varphi\in\mathcal{M}_{\bar u}^{\beta}$.
In view of \cite[Theorem 3.2.1]{MR3643081}, the zero solution is globally attractive for system \eqref{nonautonomous asymptotically periodic equations} in $\mathcal{M}_{\bar u}^{\beta}$.
Furthermore, applying the comparison argument, one can readily verify that the zero solution is also Lyapunov stable.
Hence, we obtain the global asymptotic stability of the zero solution in $\mathcal{M}_{\bar u}^{\beta}$.

Let $\lambda_\omega^{\ast}<0$.
By \cite[Theorem 3.2.1]{MR3643081}, it suffices to prove that $\lim_{k\to \infty} T_k(\varphi)=v_\omega^*(\cdot,0)$
for all $\varphi\in\mathcal{M}_{\bar u}^{\beta}\setminus \{0\}$.
Using the arguments similar to those for \cite[Theorem 3.1.5]{MR3643081}, we deduce that
the $\omega$-periodic system \eqref{nonautonomous periodic equations} has a unique positive $\omega$-periodic solution in $\mathcal{M}_{\bar u}^{\beta}$, denoted by $v_\omega^*$,
and that $v_\omega^*(\cdot,0)$ is a globally asymptotically stable fixed point of $S$ in $\mathcal{M}_{\bar u}^{\beta}\setminus\{0\}$.
Hence, $W^s(v_\omega^*(\cdot,0))=\mathcal{M}_{\bar u}^{\beta}\setminus\{0\}$, where $W^s(v_\omega^*(\cdot,0))$ is the stable set of $v_\omega^*(\cdot,0)$ for $S$.
Motivated by the proof of \cite[Proposition 3.2.3]{MR3643081}, we claim that
$\widetilde W^s(0)\cap(\mathcal{M}_{\bar u}^{\beta}\setminus\{0\})= \emptyset$,
where $\widetilde W^s(0)$ is the stable set of $0$ with respect to $T_k= \phi(k\omega,0,\cdot), k\ge 0$.
Suppose, by contradiction, that there exists a $\varphi_0\in\widetilde W^s(0)\cap(\mathcal{M}_{\bar u}^{\beta}\setminus\{0\})$,
that is, $\varphi_0\in\mathcal{M}_{\bar u}^{\beta}\setminus\{0\}$ and $T_k(\varphi_0)\to 0$ as $k\to\infty$.
It then follows from \cite[Theorem 3.2.1]{MR3643081} that $\phi(t)=\phi(t,0,\varphi_0)$ satisfies $\lim_{t\to\infty} \|\phi(t)\|_{\beta} = 0$,
and hence $\lim_{t\to\infty} \|\phi(t)\|_{C(\overline{\Omega})} = 0$.
By the strong maximum principle and the Hopf boundary lemma, there exists a sufficiently small $\alpha>0$ such that $\alpha\underline v_{\omega}(\cdot,t_0)\le \phi(\cdot,t_0)$
in $\Omega$, where $\underline v_{\omega}$ is defined in Proposition \ref{lowersolution-for-general-threshold-result}.
By the comparison principle, it follows that $\phi\ge\alpha\underline v_{\omega}\gg0$ in $\Omega\times(t_0,\infty)$,
and hence, $\liminf_{t\to\infty}\phi(x,t)>0$ for every $x\in\Omega$, which yields a contradiction.
Thus, for any $\varphi\in\mathcal{M}_{\bar u}^{\beta}\setminus\{0\}$, we have $\omega(\varphi)\cap(\mathcal{M}_{\bar u}^{\beta}\setminus\{0\})\ne\emptyset$, which means that
$\omega(\varphi)\cap W^s(v_\omega^*(\cdot,0))\ne\emptyset$. By \cite[Theorem 1.2.1]{MR3643081}, it follows that for any $\varphi\in\mathcal{M}_{\bar u}^{\beta}\setminus\{0\}$,
$\omega(\varphi)=\{v_\omega^*(\cdot,0)\}$, and hence, $\lim_{k\to\infty}\|T_k(\varphi)-v_\omega^*(0)\|_{\beta}=0$.
\end{proof}

In free boundary problems, the critical habitat size is often characterized via the principal eigenvalue
of the linearized limiting problem at the zero solution, see \cite{Du-Lin2010,Salako-Xue2026} for the one-dimensional case
and \cite{Du-Guo2011,Du-Guo-Peng2013} for the radially symmetric case in higher dimensions.
In the rest of this section, we use two examples to show that Theorem \ref{general threshold dynamics} provides
a more general framework for understanding such a critical habitat size.

\begin{example}\label{asymptotically-periodic-example-1}
\rm
Consider the following one-dimensional free boundary problem:
\begin{equation}\label{free boundary problem-vanishing-1}
\begin{cases}
u_t-d u_{xx}=u(a-bu), & t>0,\quad 0<x<h(t),\\
u_x(0,t)=u(h(t),t)=0, & t>0,\\
h(0)=h_0,\quad u(x,0)=u_0(x), & 0\le x\le h_0,
\end{cases}
\end{equation}
subject to the Stefan condition
\begin{equation}\label{Stefan condition}
h'(t)=-\mu u_x(h(t),t), \quad t>0,
\end{equation}
where $d, a, b, \mu, h_0>0$, and $u_0$ satisfies
$$
u_0 \in C^2([0, h_0]), \quad u'_0(0) = u_0(h_0) = 0, \quad \text{and }\ u_0>0\ \text{in}\ [0, h_0).
$$
Problem \eqref{free boundary problem-vanishing-1} can be used to describe the spreading of a new or invasive species with
population density $u(x,t)$ over a one-dimensional habitat. The free boundary $x=h(t)$ represents the spreading
front, while the homogeneous Neumann boundary condition at $x=0$ imposes zero population flux at the fixed endpoint.
It follows from  \cite{Du-Lin2010} that $h\in C^{1+\alpha/2}_{\rm loc}([0,\infty))$ with $\alpha\in(0,1)$,
and that if
\begin{equation}\label{free boundary problem-vanishing-2}
\lim_{t\to\infty}h(t)=h_{\infty}\le\frac{\pi}{2}\sqrt{\frac{d}{a}},
\end{equation}
then vanishing occurs, and $\lim_{t\to\infty}\|u(\cdot,t)\|_{C([0,h(t)])}=0$.

It is interesting to note that the criterion \eqref{free boundary problem-vanishing-2} coincides exactly with the threshold condition in Theorem \ref{general threshold dynamics}(a).
To see this, we make the change of variables $x=h(t)y$ and set $v(y,t)=u(h(t)y,t)$. Then \eqref{free boundary problem-vanishing-1}
can be reformulated on the fixed interval $(0,1)$ as
\begin{equation}\label{free boundary problem-vanishing-3}
\begin{cases}
v_t-\frac{d}{h^2(t)}v_{yy}-\frac{h'(t)}{h(t)}y\,v_y=v(a-bv), & t>0,\quad 0<y<1,\\
v_y(0,t)=v(1,t)=0, & t>0,\\
v(y,0)=u_0(h_0y), & 0\le y\le 1.
\end{cases}
\end{equation}
By even reflection across $y=0$, \eqref{free boundary problem-vanishing-3} extends to a Dirichlet problem on $(-1,1)$ with the same drift coefficient $-h'(t)y/h(t)$, so Theorem~\ref{general threshold dynamics} applies directly in the even subspace. Moreover, if $h_\infty<\infty$, then standard parabolic estimates yield $h'(t)\to 0$ as $t\to\infty$ (see \cite{Du-Lin2010}).
Then the limiting equation associated with \eqref{free boundary problem-vanishing-3} is
\begin{equation}\label{free boundary problem-vanishing-4}
\begin{cases}
v_t-\frac{d}{h^2_{\infty}}v_{yy}=v(a-bv), & t>0,\quad 0<y<1,\\
v_y(0,t)=v(1,t)=0, & t>0.
\end{cases}
\end{equation}
Let $\lambda^*$ denote the principal eigenvalue of the linearized operator
$L(\varphi)=\frac{d}{h^2_{\infty}}\varphi_{yy}+a\varphi$ for $\varphi\in C^2([0,1])$ satisfying $\varphi_y(0)=\varphi(1)=0$. A direct computation gives
$\lambda^*=a-\frac{d\pi^2}{4h_\infty^2}$. Note that
\[
\lambda^*\le 0\quad\Longleftrightarrow\quad h_\infty\le \frac{\pi}{2}\sqrt{\frac{d}{a}}.
\]
When $h_\infty<\infty$, Theorem~\ref{general threshold dynamics} yields the following two alternatives for every nontrivial solution of~\eqref{free boundary problem-vanishing-3}:
\begin{itemize}
	\item[(i)] if $h_{\infty}\le\frac{\pi}{2}\sqrt{\frac{d}{a}}$, then $\lim_{t\to\infty}\|v(\cdot,t)\|_{C^1([0,1])}=0$;
	\item[(ii)] if $h_{\infty}>\frac{\pi}{2}\sqrt{\frac{d}{a}}$, then $\lim_{t\to\infty}\|v(\cdot,t)-v^*(\cdot)\|_{C^1([0,1])}=0$, where $v^*$ is the unique positive steady state of~\eqref{free boundary problem-vanishing-4}.
\end{itemize}
Case~(i) corresponds precisely to the vanishing alternative for~\eqref{free boundary problem-vanishing-1}--\eqref{Stefan condition}. Under the Stefan condition~\eqref{Stefan condition}, however, Case~(ii) is impossible: the Hopf boundary lemma would then yield $\lim_{t\to\infty}h'(t)>0$, contradicting $h_\infty<\infty$. Consequently,
\[
\text{either}\quad h_{\infty}\le\frac{\pi}{2}\sqrt{\frac{d}{a}}
\qquad\text{or}\qquad
h_{\infty}=\infty,
\]
which is exactly the vanishing-spreading dichotomy established in~\cite{Du-Lin2010}. We remark that the above argument also characterizes the critical habitat size for free boundary models with spatially heterogeneous logistic growth in the autonomous~\cite{Zhou-Xiao2014} and time-periodic~\cite{Wang2016} cases, with monostable space-time periodic reaction~\cite{Ding-Du-Liang2017}, and with small time-periodic advection~\cite{Sun-Lou-Zhou2017}. The spreading case of~\eqref{free boundary problem-vanishing-1}--\eqref{Stefan condition} will be discussed in Example~\ref{degenerate-example-1}.
\end{example}

\begin{example}\label{asymptotically-periodic-example-2}
\rm
Consider~\eqref{free boundary problem-vanishing-1} equipped with the nonlocal weighted free boundary condition introduced in~\cite{Feng-Lewis-Wang-Wang2022}:
\begin{equation}\label{nonlocal weighted free boundary condition}
h'(t)=\mu\displaystyle\int_{0}^{h(t)}u(x,t)w(h(t)-x)\,\mathrm{d}x, \quad t>0,
\end{equation}
where $w(x)=c_1e^{-\alpha_1x}-c_2e^{-\alpha_2x}$ for $x\in\mathbb{R}$, with $c_1>c_2>0$ and $\alpha_1>\alpha_2>0$.
Condition \eqref{nonlocal weighted free boundary condition} generalizes the Stefan condition \eqref{Stefan condition} by incorporating nonlocal effects into boundary movement. Because $w$ changes sign on $(0,\infty)$, both range expansion and shrinkage are possible. As Theorem~\ref{general threshold dynamics} is independent of the specific boundary law, the limit analysis in Example~\ref{asymptotically-periodic-example-1} still applies. 

Since $w(0)>0$, \eqref{nonlocal weighted free boundary condition} implies that $h'(t)>0$ whenever $h(t)$ is sufficiently small, so $\inf_{t\ge 0}h(t)>0$. Moreover, since $u$ is bounded and $w\in L^1(0,\infty)$, $h'$ is bounded on $[0,\infty)$.
Suppose now that $h_\infty:=\lim_{t\to\infty}h(t)$ exists. If $h_\infty<\infty$, parabolic H\"older estimates and~\eqref{nonlocal weighted free boundary condition} imply that $h'$ is uniformly H\"older continuous on $[1,\infty)$, and hence $h'(t)\to 0$ as $t\to\infty$. Because $w$ changes sign on $(0,\infty)$, this does not rule out Case~(ii) above, which may indeed occur.
Consequently, whenever $h_\infty$ exists, all three regimes
\[
h_{\infty}\le\frac{\pi}{2}\sqrt{\frac{d}{a}},
\qquad
\frac{\pi}{2}\sqrt{\frac{d}{a}}<h_{\infty}<\infty,
\qquad\text{and}\qquad
h_{\infty}=\infty
\]
are possible, corresponding to the vanishing-balancing-spreading trichotomy established in~\cite{Feng-Lewis-Wang-Wang2022}. For results on the corresponding time-periodic problem, see~\cite{Salako-Xue2026}. We address the spreading case in Example~\ref{degenerate-example-2}.
\end{example}

\section{Asymptotically partially degenerate regime}
\label{sec:asymptotically partially degenerate}
In this section, we establish the threshold dynamics for a class of nonautonomous parabolic  systems whose differential operators exhibit asymptotic partial degeneracy in time.  For simplicity, we focus on the product-domain setting $\Omega=\Omega_1\times\Omega_2\subset\mathbb{R}^{n}$,
where $\Omega_k\subset\mathbb{R}^{n_k}$($k=1,2$) are bounded domains with boundary of class $C^{2+\theta}$ for some $0<\theta<1$,
and $n_1,n_2$ are positive integers satisfying $n_1+n_2=n$. We write $x=(x^{(1)},x^{(2)})\in\mathbb{R}^{n}$, where $x^{(1)}=(x_1,\cdots,x_{n_1})$ and $x^{(2)}=(x_{n_1+1},\cdots,x_{n})$.

Consider the following nonautonomous parabolic systems:
\begin{equation}\label{asymptotically partial degeneracy equations}
\left\{
\begin{aligned}
& \frac{\partial v_i}{\partial t} +A^{(1)}_i(t)v_i+A^{(2)}_i(t)v_i=f_i(x^{(2)},t,v) \ \ \ &{\rm in}&\ \ \ \Omega\times(0,\infty),\\
& v_i =0 \ \ \ \ &{\rm on}&\ \ \ \partial\Omega\times(0,\infty),
\end{aligned}
\right.
\end{equation}
and
\begin{equation}\label{partial degeneracy periodic equations}
\left\{
\begin{aligned}
& \frac{\partial w_i}{\partial t}+\tilde A^{(2)}_i(t)w_i=\tilde f_i(x^{(2)},t,w) \ \ \ &{\rm in}&\ \ \ \Omega_2\times(0,\infty),\\
& w_i = 0 \ \ \ \ &{\rm on}&\ \ \ \partial\Omega_2\times(0,\infty),
\end{aligned}
\right.
\end{equation}
for $1\le i\le m$, where
\begin{align*}
&A^{(1)}_i(t)v_i=-\sum_{j,k=1}^{n_1}a^{(i)}_{jk}(x,t)\partial_{jk}v_i +\sum_{j=1}^{n_1}b^{(i)}_{j}(x,t)\partial_jv_i,
\\
&A^{(2)}_i(t)v_i=-\sum_{j,k=n_1+1}^{n}a^{(i)}_{jk}(x^{(2)},t)\partial_{jk}v_i +\sum_{j=n_1+1}^{n}b^{(i)}_{j}(x^{(2)},t)\partial_jv_i,
\\
&\tilde A^{(2)}_i(t)w_i=-\sum_{j,k=n_1+1}^{n}\tilde a^{(i)}_{jk}(x^{(2)},t)\partial_{jk}w_i +\sum_{j=n_1+1}^{n}\tilde b^{(i)}_{j}(x^{(2)},t)\partial_jw_i.
\end{align*}
Throughout this section, we adopt the following hypotheses:
\begin{enumerate}[label=(B\arabic*),ref=(B\arabic*),leftmargin=*]
\item\label{ass:b1}
For each $i=1,\dots,m$, we assume that $a^{(i)}_{jk}$, $b^{(i)}_{j}\in C^{\theta, \theta/2}(\overline\Omega\times\mathbb{R}_+, \mathbb{R})$ for $1\le j, k\le n_1$,
and $a^{(i)}_{jk}$, $b^{(i)}_{j}$, $\tilde a^{(i)}_{jk}$, $\tilde b^{(i)}_{j}\in C^{\theta, \theta/2}(\overline\Omega_2\times\mathbb{R}_+, \mathbb{R})$
for $n_1+1\le j, k\le n$; $a^{(i)}_{jk}=a^{(i)}_{kj}$, $\tilde a^{(i)}_{jk}=\tilde a^{(i)}_{kj}$ for $1\le j, k\le n$; 
and there exist positive constants $\lambda_0$, $\Lambda_0$, $\lambda_i$, $\Lambda_i$, and a positive continuous function $\alpha(t)$ such that
\begin{align*}
\alpha(t)\lambda_0|\xi|^2\le\sum_{j,k=1}^{n_1}a^{(i)}_{jk}(x,t)\xi_j\xi_k\le \alpha(t)\Lambda_0|\xi|^2,& \quad (x,t)\in \Omega\times\mathbb{R}_+,
\\
\lambda_i|\zeta|^2\le\sum_{j,k=n_1+1}^{n}a^{(i)}_{jk}(x^{(2)},t)\zeta_j\zeta_k\le\Lambda_i|\zeta|^2,&\quad (x^{(2)},t)\in \Omega_2\times\mathbb{R}_+,
\end{align*}
for any $\xi=(\xi_1,\dots,\xi_{n_1})^\top$ and $\zeta=(\zeta_{n_1+1},\dots,\zeta_n)^\top$.

\item\label{ass:b2}
For $1\le i,j\le m$,
$f_i,\tilde f_i,\frac{\partial f_i}{\partial u_j},
\frac{\partial\tilde f_i}{\partial u_j}
\in C(\overline{\Omega}_2\times\mathbb{R}_+\times\mathbb{R}^m,\mathbb{R})$,
and for any fixed $T>0$,
$f_i(\cdot,\cdot,u),\tilde f_i(\cdot,\cdot,u)
\in C^{\theta,\theta/2}(\overline{\Omega}_2\times[0,T],\mathbb{R})$
uniformly for $u$ in bounded subsets of $\mathbb{R}^m$.
Moreover, $\frac{\partial\tilde f_i}{\partial u_j}(\cdot,\cdot,0)
\in C^{\theta,\theta/2}(\overline{\Omega}_2\times[0,T],\mathbb{R})$
for every $T>0$.

\item\label{ass:b3}
For each $1\le i\le m$, $\tilde A^{(2)}_i(t)$ and $\tilde f_i(\cdot,t,\cdot)$ are $\omega$-periodic in $\mathbb{R}$,
$\lim\limits_{t \to \infty}\big(\alpha(t)+\beta(t)\big) = 0$,
$$ \lim_{t \to \infty}\sum_{j,k=n_1+1}^{n}\left(\left\|a^{(i)}_{jk}(\cdot,t) - \tilde a^{(i)}_{jk}(\cdot,t)\right\|_{C(\overline\Omega_2)}+\left\| b^{(i)}_{j}(\cdot,t) - \tilde b^{(i)}_{j}(\cdot,t)\right\|_{C(\overline\Omega_2)}\right) = 0,$$
and
$$ 
\lim_{t \to \infty}\left(\left\|f_i(\cdot,t,u) -\tilde f_i(\cdot,t,u) \right\|_{C(\overline\Omega_2)}+\left\|f_u(\cdot,t,0) -\tilde f_u(\cdot,t,0) \right\|_{C(\overline \Omega_2)}\right) = 0
$$
uniformly for $u$ in any bounded sets of $\mathbb{R}^m$, where
$$
\beta(t):=\sqrt{n_1}\sup\limits_{1\le i\le m, 1\le j\le n_1}\|b^{(i)}_{j}(\cdot,t)\|_{C(\overline\Omega)}.
$$
\item\label{ass:b4}
$f_u(x^{(2)},t,u)$ is cooperative and irreducible for all $x^{(2)}\in\overline{\Omega}_2$, $t\ge0$ and $u\in\mathbb{M}_{\bar{u}}$.
For each $1\le i\le m$, $f_i(\cdot,\cdot,u)\le0$ whenever $u\in\mathbb M_{\bar u}$ and $u_i=\bar u_i$.
Moreover, $f(\cdot,\cdot,0)=0$ and $\alpha f(\cdot,\cdot,u)\le f(\cdot,\cdot,\alpha u)$ for all $u\in\mathbb{M}_{\bar{u}}\cap\mathrm{Int}(\mathbb{R}^{m}_{+})$ and $\alpha\in (0,1)$.

\item\label{ass:b5}
$\tilde f_u(x^{(2)},t,u)$ is cooperative and irreducible for all $x^{(2)}\in\overline{\Omega}_2$, $t\ge0$ and $u\in\mathbb{M}_{\bar{u}}$;
and $\alpha\tilde f(\cdot,\cdot,u)<\tilde f(\cdot,\cdot,\alpha u)$ for all  $u\in\mathbb{M}_{\bar{u}}\cap\mathrm{Int}(\mathbb{R}^{m}_{+})$ and $\alpha\in (0,1)$.
\end{enumerate}

Note that the product-type domain $\Omega$ has the exterior ball property, namely, for any
$x^0\in \partial\Omega$, there exists $R>0$ and $y\in \mathbb{R}^n\setminus\overline{\Omega}$
such that $\overline{B}_R(y)\cap \overline{\Omega}=\{x^0\}$. Moreover, there exists an increasing
sequence of subdomains $\{\Omega^k\}$ with $C^\infty$ boundaries such that
$\overline{\Omega^k}\subset \Omega^{k+1}$, $\bigcup_k \Omega^k = \Omega$, and $\partial\Omega^k\to\partial\Omega$ uniformly in the sense of Hausdorff distance as $k\to\infty$.
By adapting the parabolic regularity theory (\cite[Chapter~III.20]{MR1100011}) and standard approximation (\cite[Chapter 8.3]{WuYinW2006}), we see that for any $T>0$ and any initial datum
$\varphi\in C_0(\overline\Omega,\mathbb{M}_{\bar{u}})$, system \eqref{asymptotically partial degeneracy equations} admits a unique solution
$u\in C_{\rm loc}^{2+\theta,1+\theta/2}(\Omega\times(0,T),\mathbb{R}^m)\cap C(\overline\Omega\times[0,T],\mathbb{R}^m)$.

Let $(\lambda^*_\omega, \phi^*_\omega)$ be the principal eigen-pair of the following $\omega$-periodic parabolic eigenvalue problem:
\begin{equation}\label{parabolic eigenvalue problem-PD}
\left\{
\begin{aligned}
& \frac{\partial w}{\partial t}+\tilde A^{(2)}(t)w=\tilde f_u(x^{(2)},t,0)w+\lambda w \ \ \ &{\rm in}&\ \ \ \Omega_2\times\mathbb{R},\\
& w=0 \ \ \ \ &{\rm on}&\ \ \ \partial\Omega_2\times\mathbb{R},\\
& w(x^{(2)},\omega+t)=w(x^{(2)},t) \ \ \ &{\rm in}&\ \ \ \Omega_2\times\mathbb{R}.
\end{aligned}
\right.
\end{equation}
Then $\phi^*_\omega\in C^{2,1}(\overline\Omega_2\times\mathbb{R},\mathbb{R}^m)$ and $\phi^*_\omega\gg0$ in $\Omega_2\times\mathbb{R}$.
Since $\phi^*_\omega$ is unique up to multiplication by a positive constant, we may assume, without loss of generality,  that
$(\phi^*_\omega)_i\le1$ in $\Omega_2\times\mathbb{R}$ for each $i=1,\dots,m$.
Motivated by \cite[Lemma 5.3]{Lin-Ye-Zhao2026}, we introduce a similar auxiliary function on $\overline\Omega_1$.
\begin{lemma}\label{Sec3 auxiliary function-1}
There exist a constant $t_1>0$, a function $\phi\in C^2(\overline{\Omega}_1)$ and a positive function $H_\phi\in C([t_1,\infty))$ such that for each $i=1,\dots,m$, there holds
\begin{equation}\label{Sec3 auxiliary function-1-1}
\left\{
\begin{aligned}
& A_i^{(1)}(t)\phi\le H_\phi(t)\phi \ \ \ &{\rm in}&\ \ \ \overline\Omega\times[t_1,\infty),\\
& \phi>0 \ \ \ \ &{\rm in}&\ \ \ \Omega_1,\\
& \phi=0 \ \ \ \ &{\rm on}&\ \ \ \partial\Omega_1.
\end{aligned}
\right.
\end{equation}
Moreover, if $\lambda^*_\omega<0$ and there exists $T>0$ such that
\begin{equation}\label{Sec3 auxiliary function-1-2}
\limsup_{t\to\infty}\frac{1}{T}\int_t^{t+T}\frac{\beta(s)^2}{\alpha(s)}\,\mathrm{d}s<4\lambda_0(-\lambda^*_\omega),
\end{equation}
then $\phi$ and $H_\phi$ in \eqref{Sec3 auxiliary function-1-1} can be chosen so that
\begin{equation}\label{Sec3 auxiliary function-1-3}
\limsup_{t\to\infty}\frac{1}{T}\int_t^{t+T}H_\phi(s)\,\mathrm{d}s<-\lambda^*_\omega.
\end{equation}
\end{lemma}

\begin{proof}
For any symmetric $n_1\times n_1$ matrix $M$,  we denote by $\mathcal{M}^-_{\lambda_0,\Lambda_0}$
Pucci's lower extremal operator with parameters $0<\lambda_0\le\Lambda_0$, that is,
$\mathcal{M}^-_{\lambda_0,\Lambda_0}(M)=\inf_{A\in\mathfrak{A}}{\rm tr}(AM)$,
where $\mathfrak{A}$ denotes the set of all symmetric matrices whose eigenvalues lie in the interval $[\lambda_0,\Lambda_0]$.
By the uniform ellipticity of $\frac{1}{\alpha(t)}A^{(1)}_i(t)$ for each $1\le i\le m$, we have
$$
\sum_{j,k=1}^{n_1}a^{(i)}_{jk}(x,t)\partial_{jk}v={\rm tr}(\mathcal{A}^{(1)}_i\mathcal{D}_1^2v)\ge \alpha(t)\mathcal{M}^-_{\lambda_0,\Lambda_0}(\mathcal{D}_1^2v),
$$
where $v\in C^{2}(\overline\Omega_1, \mathbb{R})$, $\mathcal{A}^{(1)}_i=(a^{(i)}_{jk})_{n_1\times n_1}$ and $\mathcal{D}_1^2$ denotes Hessian matrix with respect to the variables $x^{(1)}$.
It then follows that
\begin{equation*}
A^{(1)}_i(t)v\le-\alpha(t)\mathcal{M}^-_{\lambda_0,\Lambda_0}(\mathcal{D}_1^2v)+\beta(t)|\nabla_1 v|,
\end{equation*}
where $\nabla_1$ denotes the gradient with respect to the variables $x^{(1)}$.

Fix $p\ge 2$ and define $\phi(x^{(1)})=\left(\rho(x^{(1)})\right)^p$,
where $\rho\in C^2(\overline{\Omega}_1)$ satisfies $0<\rho\le1$ in $\Omega_1$, $\rho=0$ on $\partial\Omega_1$,
and for some sufficiently small $\eta\in(0,1)$, $\rho(x^{(1)})=\operatorname{dist}(x^{(1)},\partial\Omega_1)$,
$x^{(1)}\in\Omega_1\setminus\Omega_1^\eta$,
with $\Omega_1^\eta:=\left\{x^{(1)}\in\Omega_1: {\rm dist}(x^{(1)},\partial\Omega_1)> \eta\right\}$.
The auxiliary function introduced in \cite[Lemma 5.3]{Lin-Ye-Zhao2026} corresponds to the special case $p=2$ in the above construction.

Define
\begin{align*}
h_\phi(t):=\sup_{x^{(1)}\in\Omega_1}\frac{-\alpha(t)\mathcal{M}^-_{\lambda_0,\Lambda_0}(\mathcal{D}_1^2\phi)+\beta(t)|\nabla_1\phi|}{\phi}.
\end{align*}
We will show that for each fixed $p\ge 2$, $h_\phi(t)<\infty$ for all $t\in[0,\infty)$.
To this end, we treat separately the interior region $x^{(1)}\in\Omega_1^\eta$ and the boundary layer $x^{(1)}\in\Omega_1\setminus\Omega_1^\eta$.
In the case where $x^{(1)}\in\Omega_1^\eta$, since $\phi$ has a positive lower bound on $\Omega_1^\eta$,  it follows from assumption \ref{ass:b3} that
\begin{align*}
\limsup_{t\to\infty}\sup_{x^{(1)}\in\Omega_1^\eta}\frac{-\alpha(t)\mathcal{M}^-_{\lambda_0,\Lambda_0}(\mathcal{D}_1^2\phi)+\beta(t)|\nabla_1\phi|}{\phi}=0.
\end{align*}

In the case where $x^{(1)}\in\Omega_1\setminus\Omega_1^\eta$, we have
$|\nabla_1\rho|=1$ in $\Omega_1\setminus\Omega_1^\eta$, and $|\nabla_1\phi|=p\rho^{p-1}|\nabla_1\rho|=0$ on $\partial\Omega_1$.
A direct calculation gives
$$
\mathcal{D}_1^2\phi = p(p-1)\rho^{p-2}\nabla_1\rho\otimes\nabla_1\rho + p\rho^{p-1}\mathcal{D}_1^2\rho,
$$
where $\otimes$ means tensor product.
Note that $\nabla_1\rho\otimes\nabla_1\rho$ in $\Omega_1\setminus\Omega_1^\eta$ is a rank-one and positive semidefinite matrix,
with eigenvalue 1 in the normal direction and eigenvalue 0 in the tangential directions.
Therefore, by the properties of Pucci's extremal operators (see \cite[Section 2.2]{Caffarelli1995}), we obtain
\begin{align*}
    \mathcal{M}_{\lambda_0, \Lambda_0}^{-}(\mathcal{D}_1^2\phi) &= \mathcal{M}_{\lambda_0, \Lambda_0}^{-}\Big(p(p-1)\rho^{p-2}\nabla_1\rho\otimes\nabla_1\rho + p\rho^{p-1}\mathcal{D}_1^2\rho\Big)      \nonumber
    \\      \nonumber
    & \geq \mathcal{M}_{\lambda_0, \Lambda_0}^{-}\Big(p(p-1)\rho^{p-2}\nabla_1\rho\otimes\nabla_1\rho\Big) + \mathcal{M}_{\lambda_0, \Lambda_0}^{-}\Big(p\rho^{p-1}\mathcal{D}_1^2\rho\Big)      \nonumber
    \\      \nonumber
    &  =  p(p-1)\rho^{p-2}\mathcal{M}_{\lambda_0, \Lambda_0}^{-}\Big(\nabla_1\rho\otimes\nabla_1\rho\Big) + p\rho^{p-1}\mathcal{M}_{\lambda_0, \Lambda_0}^{-}\Big(\mathcal{D}_1^2\rho\Big)      \nonumber
    \\
    & \geq \lambda_0p(p-1)\rho^{p-2}-c\rho^{p-1}\ \ \text{ in }\ \ \Omega_1\setminus\Omega_1^\eta,
\end{align*}
where $c>0$ is a constant depending on $p$ and $\mathcal{D}_1^2\rho$.
Then for $x^{(1)}\in\Omega_1\setminus\Omega_1^\eta$, there holds
\begin{align}\label{Sec3 auxiliary function-1-7}
\frac{-\alpha(t)\mathcal{M}^-_{\lambda_0,\Lambda_0}(\mathcal{D}_1^2\phi)+\beta(t)|\nabla_1\phi|}{\phi}
\le \frac{-\lambda_0p(p-1)\alpha(t)}{\rho^2} +\frac{c\alpha(t)+p\beta(t)}{\rho}
\le h_p(t),
\end{align}
where the last inequality follows by maximizing the quadratic in $1/\rho$, and
\begin{align*}
h_p(t):&=\frac{\bigl(c\alpha(t)+p\beta(t)\bigr)^2}{4\lambda_0 p(p-1)\alpha(t)}
=  \frac{c^2}{4\lambda_0 p(p-1)}\alpha(t)+\frac{c}{2\lambda_0 (p-1)}\beta(t)+\frac{p}{4\lambda_0 (p-1)}\frac{\beta^2(t)}{\alpha(t)}.
\end{align*}
Clearly, after increasing $c$ in $h_p$ if necessary, there exists $t_1>0$ independent of $p$, such that $h_\phi(t)\le h_p(t)$ for all $t\ge t_1$. Consequently, we have
$$
A^{(1)}_i(t)\phi\le-\alpha(t)\mathcal{M}_{\lambda_0,\Lambda_0}^{-}\bigl(\mathcal{D}_1^2\phi\bigr)+\beta(t)|\nabla_1\phi|\le h_p(t)\phi, \quad x^{(1)}\in\overline\Omega_1,\ \ t\ge t_1.
$$
Thus, we may take $H_\phi=h_p$ on $[t_1,\infty)$.
Furthermore, combining \eqref{Sec3 auxiliary function-1-2} with \eqref{Sec3 auxiliary function-1-7},
we can choose $p$ sufficiently large so that
\begin{align*}
\limsup_{t\to\infty}\frac{1}{T}\int_t^{t+T}H_\phi(s)\,\mathrm{d}s=\limsup_{t\to\infty}\frac{1}{T}\int_t^{t+T}h_p(s)\,\mathrm{d}s<-\lambda^*_\omega.
\end{align*}
This completes the proof.
\end{proof}

Applying Lemma \ref{Sec2 auxiliary functions} to the operator $A^{(2)}(t)$, we obtain a positive constant and an auxiliary vector-valued barrier function,
which are still denoted by $\mu$ and $\Phi=(\Phi_1,\dots,\Phi_m)^\top$, respectively,
such that $\Phi\in C^{2}(\overline\Omega_2,\mathbb{R}^m)$, $\Phi\gg0$ in $\Omega_2$, $\Phi=0$ on $\partial\Omega_2$, and
\begin{equation}\label{eigenvalue-PD-2}
-(A^{(2)}(t)+\mu)\Phi\le (-1,\dots,-1)^\top \quad {\rm in}\ \ \Omega_2.
\end{equation}
Based on Lemma \ref{Sec3 auxiliary function-1}, we next construct a subsolution $\underline v$ of system \eqref{asymptotically partial degeneracy equations}.

\begin{proposition}\label{Sec3 subsolution-1}
Assume that $\lambda_\omega^*<0$ and \eqref{Sec3 auxiliary function-1-2} holds.
Then there exist $t_2>t_1$ and a function $\underline v\in C^{2,1}(\overline\Omega\times[t_2,\infty),\mathbb{R}^m)$ such that
$\underline v$ is a subsolution of system \eqref{asymptotically partial degeneracy equations} in $\Omega\times(t_2,\infty)$
and satisfies $\underline v\gg0$ in $\Omega\times[t_2,\infty)$, $\underline v=0$ on $\partial\Omega\times[t_2,\infty)$.
Moreover, for every compact set $K\subset\Omega$, there holds
\begin{equation}\label{Sec3 subsolution-1-1}
\min_{1\le i\le m}\inf_{x\in K,\, t\ge t_2}\underline v_i(x,t)>0.
\end{equation}
\end{proposition}

\begin{proof}
We divide the proof into two steps.

\medskip
{\noindent\it Step 1.} 
By \eqref{Sec3 auxiliary function-1-3}, we choose $\varepsilon>0$ sufficiently small such that
\[
\limsup_{t\to\infty}\frac{1}{T}\int_t^{t+T}H_\phi(s)\,\mathrm{d}s<-\lambda^*_\omega-5\varepsilon,
\]
and fix $\tilde t_1>t_1$ such that
\begin{equation}\label{Sec3 subsolution-1-2}
	\frac{1}{T}\int_t^{t+T}\bigl(\lambda^*_\omega+4\varepsilon+H_\phi(s)\bigr)\,\mathrm{d}s<-\varepsilon,\qquad \forall t\ge \tilde t_1.
\end{equation}
By Hopf boundary lemma, there exists $C_*>0$ such that $\max_{1\le j\le m}(\phi_\omega^*)_j\le C_*\,(\phi_\omega^*)_i$ in $\overline\Omega_2\times[0,\omega]$ for all $1\le i\le m$. Arguing as in the derivation of \eqref{lowersolution-for-general-threshold-result-2-1} in Proposition~\ref{lowersolution-for-general-threshold-result}, we deduce from \ref{ass:b2}--\ref{ass:b4} that there exists $r_*\in(0,1)$ with $\{z\in\mathbb R^m: z\ge 0,\ \|z\|_{\mathbb R^m}\le r_*\}\subset\mathbb M_{\bar u}$ such that, for all sufficiently large $t$ and all $x^{(2)}\in\overline\Omega_2$,
\begin{equation}\label{Sec3 subsolution-1-3}
	f(x^{(2)},t,z)
	\ge f_u(x^{(2)},t,0)z
	-\frac{\varepsilon}{C_*}\|z\|_{\mathbb{R}^m}(1,\dots,1)^\top,
	\qquad z\ge 0,\quad \|z\|_{\mathbb{R}^m}\le r_*.
\end{equation}

Let $q(t)>0$ solve $q'=q\bigl(-\lambda^*_\omega-4\varepsilon-H_\phi(t)-Cq\bigr)$ on $[\tilde t_1,\infty)$ with $q(\tilde t_1)\in(0,r_*)$, where $C\ge(-\lambda^*_\omega-4\varepsilon)/r_*$ ensures $0<q(t)\le r_*$ for all $t\ge \tilde t_1$. Setting $g:=1/q$ gives $g'=\bigl(\lambda_\omega^*+4\varepsilon+H_\phi(t)\bigr)g+C$, and hence
$$
g(t)=\mathrm{e}^{\int_{\tilde t_1}^{t}\bigl(\lambda_\omega^*+4\varepsilon+H_\phi(r)\bigr)\,\mathrm{d}r}g(\tilde t_1)+C\int_{\tilde t_1}^{t}\mathrm{e}^{\int_s^t \bigl(\lambda_\omega^*+4\varepsilon+H_\phi(r)\bigr)\,\mathrm{d}r}\,\mathrm{d}s,
\qquad t\ge \tilde t_1.
$$
Since $H_\phi>0$, \eqref{Sec3 subsolution-1-2} yields $\sup_{t\ge \tilde t_1}g(t)<\infty$, which implies that $\underline q:=\inf_{t\ge \tilde t_1}q(t)>0$. By Hopf boundary lemma, there exists $c_*>0$ such that
\begin{equation}\label{phi-Phi-PD}
	c_*\Phi\le\underline q\phi^*_\omega \quad {\rm on}\ \ \overline\Omega_2\times\mathbb{R}.
\end{equation}
By \ref{ass:b3}, we can choose $t_2\ge \tilde t_1$ sufficiently large so that \eqref{Sec3 subsolution-1-3} holds for all $t\ge t_2$, and in $\overline\Omega_2\times[t_2,\infty)$,
\begin{equation}\label{Sec3 subsolution-1-5}
	\sum_{i,j=1}^m\left\|\frac{\partial f_i}{\partial u_j}(\cdot,\cdot,0)-\frac{\partial \tilde f_i}{\partial u_j}(\cdot,\cdot,0)\right\|_{C(\overline\Omega_2\times[t_2,\infty))}\le \frac{\varepsilon}{C_*},
\end{equation}
\begin{equation}\label{Sec3 subsolution-1-6}
	\left(A^{(2)}(t)-\tilde A^{(2)}(t)\right)\phi^*_\omega\le\delta_2(1,\dots,1)^\top,
\end{equation}
where
\begin{equation}\label{varepsilon-1-PD}
	\delta_2=\min\left\{\frac{c_*}{2},\ \frac{\varepsilon c_*}{\mu},\ \frac{\varepsilon c_*}{C_*\|f_u(\cdot,\cdot,0)\|_{C(\overline\Omega_2\times\mathbb{R}_+)}}\right\}.
\end{equation}
Define
\begin{equation}\label{Sec3 subsolution-1-7}
	\underline v(x,t):=\delta_1\phi(x^{(1)})\Big(q(t)\phi^*_\omega(x^{(2)},t)-\delta_2\Phi(x^{(2)})\Big),
\end{equation}
where $\delta_1\in(0,\frac{1}{\|\phi\|_{C(\overline\Omega_1)}\|\phi^*_\omega\|_{C(\overline\Omega_2\times\mathbb{R})}})$. 
By \eqref{phi-Phi-PD} and \eqref{varepsilon-1-PD}, we have
\begin{equation}\label{Sec3 subsolution-1-9}
	\frac{c_*}{2}\delta_1\phi(x^{(1)})\Phi(x^{(2)})\le \underline v(x,t)\ll\delta_1\phi(x^{(1)})q(t)\phi^*_\omega(x^{(2)},t)\le r_*(1,\dots,1)^\top
\end{equation}
in $\Omega\times[t_2,\infty)$, from which \eqref{Sec3 subsolution-1-1} follows for every compact set $K\subset\Omega$.

\medskip
{\noindent\it Step 2.} 
Since $q\phi^*_\omega-\delta_2\Phi\ge 0$ and $H_\phi>0$, Lemma~\ref{Sec3 auxiliary function-1} gives $A^{(1)}(t)\underline v\le\delta_1 qH_\phi\phi\phi^*_\omega$ in $\Omega\times[t_2,\infty)$. Using \eqref{parabolic eigenvalue problem-PD} and $q'+\bigl(\lambda_\omega^*+H_\phi\bigr)q=-\bigl(4\varepsilon+Cq\bigr)q$, a direct calculation in $\Omega\times(t_2,\infty)$ yields
\begin{align}\label{Sec3 subsolution-1-8}
    &\frac{\partial \underline v}{\partial t}+A^{(1)}(t)\underline v+A^{(2)}(t)\underline v-f(x^{(2)},t,\underline v)      \nonumber
    \\      \nonumber
    &  \le  \delta_1q\phi\frac{\partial \phi^*_\omega}{\partial t}+\delta_1q'\phi\phi^*_\omega+\delta_1qH_\phi\phi\phi^*_\omega+\delta_1q\phi A^{(2)}(t)\phi^*_\omega-\delta_1\delta_2\phi A^{(2)}(t)\Phi-f(x^{(2)},t,\underline v)
    \\      \nonumber
    &  =  \Big(q'+\lambda_\omega^*q+qH_\phi\Big)\delta_1\phi\phi^*_\omega+\mu\delta_1\delta_2\phi\Phi-\delta_1\delta_2\phi \Big(A^{(2)}(t)+\mu\Big)\Phi
    \\      \nonumber
    & {\hspace{17pt}} +\delta_1q\phi\left(A^{(2)}(t)-\tilde A^{(2)}(t)\right)\phi^*_\omega+\delta_1\delta_2\phi f_u(x^{(2)},t,0)\Phi
    \\      \nonumber
    & {\hspace{17pt}} -\Big(f(x^{(2)},t,\underline v)-f_u(x^{(2)},t,0)\underline v\Big)+\delta_1q\phi\Big(\tilde f_u(x^{(2)},t,0)-f_u(x^{(2)},t,0)\Big)\phi^*_\omega
    \\       \nonumber
    &  \le   \Big(-4\varepsilon q-Cq^2\Big)\delta_1\phi\phi^*_\omega+\mu\delta_1\delta_2\phi\Phi+\Big(\tilde f_u(x^{(2)},t,0)-f_u(x^{(2)},t,0)\Big)\delta_1q\phi\phi^*_\omega
    \\       \nonumber
    & {\hspace{17pt}} +\delta_1\phi\left[q\left(A^{(2)}(t)-\tilde A^{(2)}(t)\right)\phi^*_\omega-\delta_2\Big(A^{(2)}(t)+\mu\Big)\Phi\right]
    \\
    & {\hspace{17pt}} +\delta_1\delta_2\phi f_u(x^{(2)},t,0)\Phi-\Big(f(x^{(2)},t,\underline v)-f_u(x^{(2)},t,0)\underline v\Big).
\end{align}
In view of  \eqref{phi-Phi-PD} and \eqref{varepsilon-1-PD}, we have
$\mu\delta_1\delta_2\phi\Phi\le\frac{\mu\delta_1\delta_2}{c_*}q\phi\phi^*_\omega\le\varepsilon\delta_1q\phi\phi^*_\omega$.
Then \eqref{Sec3 subsolution-1-5} implies
$
\delta_1q\phi\big(\tilde f_u(x^{(2)},t,0)-f_u(x^{(2)},t,0)\big)\phi^*_\omega\le\varepsilon\delta_1q\phi\phi^*_\omega.
$
Combining  \eqref{eigenvalue-PD-2} and \eqref{Sec3 subsolution-1-6}, we obtain
\begin{align*}
\delta_1\phi\left[q\left(A^{(2)}(t)-\tilde A^{(2)}(t)\right)\phi^*_\omega-\delta_2\Big(A^{(2)}(t)+\mu\Big)\Phi\right]\le0.
\end{align*}
Further, \eqref{phi-Phi-PD}, \eqref{varepsilon-1-PD} and the choice of $C_*$ give rise to
\begin{align*}
    \delta_1\delta_2\phi f_u(x^{(2)},t,0)\Phi
    & \le \frac{\delta_1\delta_2}{c_*}q\phi C_*\big\|f_u(\cdot,\cdot,0)\big\|_{C(\overline\Omega_2\times\mathbb{R}_+)}\phi^*_\omega
        \le \varepsilon\delta_1q\phi\phi^*_\omega.
\end{align*}
By \eqref{Sec3 subsolution-1-3} and \eqref{Sec3 subsolution-1-9}, we have
$$
-\big(f(x^{(2)},t,\underline v)-f_u(x^{(2)},t,0)\underline v\big)
\le \frac{\varepsilon}{C_*}\|\underline v\|_{\mathbb R^m}(1,\dots,1)^\top
\le \varepsilon\delta_1q\phi\,\phi_\omega^*.
$$
Combining the above estimates and substituting them into \eqref{Sec3 subsolution-1-8}, we derive
\begin{align*}
\frac{\partial \underline v}{\partial t}+A^{(1)}(t)\underline v+A^{(2)}(t)\underline v-f(x^{(2)},t,\underline v)\ll0,\quad (x,t)\in\Omega\times(t_2,\infty).
\end{align*}
This completes the proof.
\end{proof}

The subsolution $\underline v$ above keeps nontrivial solutions uniformly bounded away from zero on each compact subset of $\Omega$. To recover the full limiting profile, we next construct a second subsolution.

Define $\Omega_2^\eta:=\left\{x^{(2)}\in\Omega_2: d(x^{(2)},\partial\Omega_2)> \eta\right\}$ with $\eta>0$.
It is clear that $\Omega_2^\eta\uparrow\Omega_2$ as $\eta\downarrow0^{+}$, and $\partial\Omega_2^\eta=\left\{x^{(2)}-\eta\nu(x^{(2)}); x^{(2)}\in\partial\Omega_2\right\}$
for sufficiently small $\eta$, where $\nu(x^{(2)})$ is the unit outward normal vector on $\partial\Omega_2$. Thus, there exists small $\eta_0>0$ such that the $C^{2+\theta}$-norm of $\partial\Omega_2^{\eta}$ is uniformly bounded with respect to $\eta\in(0,\eta_0)$.
Consider the following auxiliary problem:
\begin{equation}\label{Omega-2-eta-problem}
\left\{
\begin{aligned}
& \frac{\partial w}{\partial t}+A^{(2)}(t)w=f(x^{(2)},t,w) \ \ \ &{\rm in}&\ \ \ \Omega_2^\eta\times(t_2,\infty),\\
& w = 0 \ \ \ \ &{\rm on}&\ \ \ \partial\Omega_2^\eta\times(t_2,\infty),
\end{aligned}
\right.
\end{equation}
and its limiting problem
\begin{equation}\label{limiting-Omega-2-eta-problem}
\left\{
\begin{aligned}
& \frac{\partial w}{\partial t}+\tilde A^{(2)}(t)w=\tilde f(x^{(2)},t,w) \ \ \ &{\rm in}&\ \ \ \Omega_2^\eta\times(t_2,\infty),\\
& w = 0 \ \ \ \ &{\rm on}&\ \ \ \partial\Omega_2^\eta\times(t_2,\infty),
\end{aligned}
\right.
\end{equation}
where $A^{(2)}(t)w=(A^{(2)}_1(t)w_1,\dots,A^{(2)}_m(t)w_m)^\top$ and $\tilde A^{(2)}(t)w=(\tilde A^{(2)}_1(t)w_1,\dots,\tilde A^{(2)}_m(t)w_m)^\top$.

Let $\lambda_\omega^{\ast}$ and $\lambda^\eta_\omega$ be the principal eigenvalues of the $\omega$-periodic parabolic
eigenvalue problems associated with the linearizations of systems
\eqref{partial degeneracy periodic equations} and \eqref{limiting-Omega-2-eta-problem} at the zero solution,
respectively. Assume that $\lambda_\omega^{\ast}<0$. By Theorem~\ref{general threshold dynamics},
system \eqref{partial degeneracy periodic equations} admits a unique positive $\omega$-periodic solution $w_\omega^*$ in $\mathcal{M}_{\bar u}^{\beta}$,
which is globally asymptotically stable in $\mathcal{M}_{\bar u}^{\beta}\setminus \{0\}$.
Moreover, by the spectral continuity for periodic-parabolic operators (see \cite[Theorem~5.3]{Daners1996JDE}),
the sign of the principal eigenvalue is preserved under small domain perturbations, so we may choose $\eta_0>0$ sufficiently small such that $\lambda^\eta_\omega<0$ for all $\eta\in(0,\eta_0)$.
Denote the unique positive $\omega$-periodic solution of system \eqref{limiting-Omega-2-eta-problem} by $w_\omega^\eta$, and
extend it by zero to $\overline\Omega_2\setminus\overline\Omega_2^{\eta}$ for all $t\in\mathbb R$. Thus, we have the following convergence result.

\begin{lemma}\label{convergence-domain perturbation}
	Assume that $\lambda_\omega^*<0$. Then $w_\omega^\eta \uparrow w_\omega^*$ uniformly on $\overline{\Omega}_2 \times [0, \omega]$ as $\eta \downarrow 0$.
\end{lemma}

\begin{proof}
For $0<\underline\eta<\overline\eta<\eta_0$, we have $\Omega_2^{\overline\eta}\subset\Omega_2^{\underline\eta}\subset\Omega_2$. Consider \eqref{Omega-2-eta-problem} on these three domains, starting at $t=t_2$ from the same nonzero smooth initial datum taking values in $\mathbb M_{\bar u}$ and compactly supported in $\Omega_2^{\overline\eta}$. By the comparison principle, the corresponding solutions remain nonnegative and ordered by domain inclusion. Applying Theorem~\ref{general threshold dynamics} at $t+k\omega$ ($t\in[0,\omega]$) and letting $k\to\infty$, we obtain, after zero extension to $\overline\Omega_2$,
\begin{equation}\label{convergence-domain perturbation-2}
	w_\omega^*\ge w_\omega^{\underline\eta}\ge w_\omega^{\overline\eta}\ge 0 \quad \text{on } \overline\Omega_2\times\mathbb R.
\end{equation}

Fix a sequence $\eta_l\downarrow 0$ in $(0,\eta_0)$. By \eqref{convergence-domain perturbation-2}, the pointwise limit $W:=\lim_{l\to\infty} w_\omega^{\eta_l}$ exists on $\overline\Omega_2\times\mathbb R$ and is $\omega$-periodic in $t$. For each compact subset $K\subset\Omega_2$, periodicity and interior parabolic $L^p$ and Schauder estimates give an index $l_K$ and $\sigma\in(0,\theta)$ such that
\[
\sup_{l\ge l_K} \|w_\omega^{\eta_l}\|_{C^{2+\sigma,1+\sigma/2}(K\times[0,\omega])} < \infty.
\]
By a diagonal argument, passing to a subsequence if necessary, $w_\omega^{\eta_l}\to W$ in $C_{\mathrm{loc}}^{2,1}(\Omega_2\times\mathbb R)$. Passing to the limit in \eqref{limiting-Omega-2-eta-problem} along this subsequence shows that $W$ is an $\omega$-periodic classical solution of the differential equation in \eqref{partial degeneracy periodic equations} on $\Omega_2\times\mathbb R$.
By \eqref{convergence-domain perturbation-2}, for each $l\ge 1$,
\[
0 \ll w_\omega^{\eta_l} \le W \le w_\omega^* \quad \text{in } \Omega_2^{\eta_l}\times\mathbb R.
\]
Since $\Omega_2^{\eta_l}\uparrow\Omega_2$, we have $W\gg 0$ in $\Omega_2\times\mathbb R$. Moreover, since $0\le W\le w_\omega^*$ on $\overline\Omega_2\times\mathbb R$ and $w_\omega^*=0$ on $\partial\Omega_2\times\mathbb R$, $W$ is continuous on $\overline\Omega_2\times\mathbb R$ and vanishes on $\partial\Omega_2\times\mathbb R$. Standard global parabolic regularity then ensures $W(\cdot,t)\in\mathcal{M}_{\bar u}^\beta$ for all $t\in\mathbb R$, so the uniqueness of the positive $\omega$-periodic solution in $\mathcal{M}_{\bar u}^\beta$ yields $W=w_\omega^*$. Finally, \eqref{convergence-domain perturbation-2} and Dini's theorem give
\[
\sup_{0<\eta\le\eta_l}\|w_\omega^*-w_\omega^\eta\|_{C(\overline\Omega_2\times[0,\omega])}
\le \|w_\omega^*-w_\omega^{\eta_l}\|_{C(\overline\Omega_2\times[0,\omega])} \to 0 \quad \text{as } l\to\infty. 
\]
\end{proof}

Choose $\eta_0$ small enough that $\Omega_1^{2\eta_0}$ is nonempty, and fix $\eta_1, \eta_2 \in (0, \eta_0)$. Similarly as $\Omega_2^{\eta_2}$, we may define $\Omega_1^{\eta_1}$ for sufficiently small $\eta_0$.
Let $\varrho_{\eta_1}$ be a smooth cut-off function on $\mathbb{R}^{n_1}$ satisfying $\varrho_{\eta_1}=1$ in $\Omega_1^{2\eta_1}$,
$\varrho_{\eta_1}=0$ in $\mathbb{R}^{n_1}\setminus\Omega_1^{\eta_1}$, and $0\le\varrho_{\eta_1}\le1$.
Define
$$
\mathcal{W}:= \left\{w\in C(\overline\Omega_2^{\eta_2}\times[t_2,\infty),\mathbb{M}_{\bar{u}}): \ w=0 \ \ {\rm on}\ \ \partial\Omega_2^{\eta_2}\times[t_2,\infty)\right\},
$$
where $t_2$ is given in Proposition \ref{Sec3 subsolution-1}. By definition, $w_\omega^{\eta_2}\in\mathcal{W}$.
Combining Proposition \ref{Sec3 subsolution-1} with Lemma \ref{convergence-domain perturbation}, we can choose $\delta_1$ sufficiently small in \eqref{Sec3 subsolution-1-7}
so that there exist $x_0^{(1)}\in\Omega_1^{2\eta_1}$ and $x_0^{(2)}\in\Omega_2^{\eta_2}$ satisfying
\begin{equation}\label{nonempty}
w_\omega^{\eta_2}(x_0^{(2)},t)\gg\underline v(x_0^{(1)},x_0^{(2)},t),\qquad \forall t\ge t_2,\ \ \forall\eta_1, \eta_2\in(0,\eta_0).
\end{equation}
For any fixed $w\in\mathcal{W}$,
extend $w$ to be zero on $\overline\Omega_2\setminus\overline\Omega_2^{\eta_2}$ for any $t\ge t_2$.
Define
$$
\mathcal{D}_w:= \left\{(x,t)\in\Omega\times[t_2,\infty): \varrho_{\eta_1}(x^{(1)})w(x^{(2)},t)>\underline v(x,t)\right\}.
$$
It follows from \eqref{nonempty} that $\mathcal{D}_{w}\ne\varnothing$ for some $w\in\mathcal{W}$.

\begin{lemma}\label{lowersolution-2-PD}
Assume that $\lambda_\omega^*<0$, and $\eta_1, \eta_2\in(0,\eta_0)$.
Then we have
$$
\inf_{w\in\mathcal{W},\mathcal{D}_w\ne\emptyset}\,\inf_{(x,t)\in\mathcal{D}_w}\varrho_{\eta_1}(x^{(1)})w(x^{(2)},t)\gg0.
$$
\end{lemma}

\begin{proof}
Suppose, for contradiction, that this is not the case. There exist $i_0\in\{1,\dots,m\}$, $w_k\in\mathcal{W}$ and $(x_k,t_k)\in\mathcal{D}_{w_k}$ such that
$$
\underline v_{i_0}(x_k,t_k)\le\varrho_{\eta_1}(x_k^{(1)})\left(w_k\right)_{i_0}(x_k^{(2)},t_k)\to0 \quad {\rm as}\ \ \ k\to\infty.
$$
Note that $x_k\in\Omega_1^{\eta_1}\times\Omega_2^{\eta_2}$ for all $k$.
It then follows from the compactness that there exist a limiting point $\tilde x\in\overline\Omega_1^{\eta_1}\times\overline\Omega_2^{\eta_2}$
and a subsequence $\{(x_{k_l},t_{k_l})\}_{l\ge1}$ such that
$\lim_{l\to\infty}x_{k_l}=\tilde x$ and $\lim_{l\to\infty}\underline v_{i_0}(x_{k_l},t_{k_l})=0$,
and hence, $\tilde x\in\partial\Omega$. This  is a contradiction.
\end{proof}

\begin{remark}\label{lowersolution-2-PD-1}
In view of \eqref{nonempty} and Lemma \ref{lowersolution-2-PD}, one can choose a positive function $\epsilon\in C^1([t_2,\infty))$ with $\epsilon(t)\downarrow0$ as $t\to\infty$, such that
for every $w\in\mathcal{W}$ with $\mathcal{D}_w\ne\emptyset$, the following estimate holds with the same $\epsilon$:
\begin{equation}\label{lowersolution-2-PD-1-1}
A^{(1)}(t)\left(\varrho_{\eta_1}(x^{(1)})w(x^{(2)},t)\right)\le \epsilon(t)\varrho_{\eta_1}(x^{(1)})w(x^{(2)},t), \quad (x,t)\in\mathcal{D}_w.
\end{equation}
\end{remark}

Now, we let $w^{\eta_2}_\varphi$ denote the nontrivial solution of the following problem
\begin{equation}\label{Omega-2-eta-problem-2}
\left\{
\begin{aligned}
& \frac{\partial w}{\partial t}+A^{(2)}(t)w=f(x^{(2)},t,w)-\epsilon(t)w \ \ \ &{\rm in}&\ \ \ \Omega_2^{\eta_2}\times(t_2,\infty),\\
& w = 0 \ \ \ \ &{\rm on}&\ \ \ \partial\Omega_2^{\eta_2}\times(t_2,\infty),\\
& w(\cdot,t_2) = \varphi \ \ \ \ &{\rm in}&\ \ \ \Omega_2^{\eta_2},
\end{aligned}
\right.
\end{equation}
where $\varphi\in C_0(\overline\Omega_2^{\eta_2},\mathbb{M}_{\bar{u}})\setminus\{0\}$. Clearly, $w^{\eta_2}_\varphi\in\mathcal{W}$.
If $\lambda_\omega^*<0$, then the definition of $\eta_0$ and Theorem \ref{general threshold dynamics} imply that
\begin{equation}\label{asymptotic-Omega-2-eta-problem-2}
\lim\limits_{t\to+ \infty} \|w^{\eta_2}_\varphi(\cdot,t)-w_\omega^{\eta_2}(\cdot,t)\|_{C^1(\overline\Omega_2^{\eta_2})}=0,
\end{equation}
where $w_\omega^{\eta_2}$ denotes the unique positive $\omega$-periodic solution of \eqref{limiting-Omega-2-eta-problem} in $\Omega_2^{\eta_2}\times[0,\omega]$.
We extend both $w^{\eta_2}_\varphi$ and $w_\omega^{\eta_2}$ by zero to $\overline\Omega_2\setminus\overline\Omega_2^{\eta_2}$ for all $t\ge t_2$.

\begin{proposition}\label{lowersolution-3-PD}
Assume that $\lambda_\omega^*<0$, $\eta_1, \eta_2\in(0,\eta_0)$ and \eqref{Sec3 auxiliary function-1-2} holds. Then the componentwise maximum $\underline u(x,t):=\max\left\{\varrho_{\eta_1}(x^{(1)})w^{\eta_2}_\varphi(x^{(2)},t),\ \underline v(x,t)\right\}$
is a generalized subsolution for system \eqref{asymptotically partial degeneracy equations} in $\Omega\times(t_2,\infty)$.
\end{proposition}

Here and in Proposition \ref{Sec4 subsolution-2}, we adopt the definitions of generalized subsolutions as given in \cite[Definition 1.1.1]{lam2022},
which enables the gluing of classical subsolutions.
The corresponding comparison principle for cooperative parabolic systems is provided in \cite[Theorem 3.1.5]{lam2022}.

\begin{proof}
Note that when $\varphi$ satisfies $\varrho_{\eta_1}(x^{(1)})\varphi(x^{(2)})\le\underline v(x,t_2)$ for $x\in\Omega$, we have $\underline u(x,t_2)=\underline v(x,t_2)$ for all $x\in\Omega$.
Whenever $\varrho_{\eta_1}(x^{(1)})(w^{\eta_2}_\varphi)_i(x^{(2)},t)\ge\underline v_i(x,t)$ for some $1\le i\le m$ in $\Omega\times(t_2,\infty)$, we have $x\in\Omega_1^{\eta_1}\times\Omega_2^{\eta_2}$. By \eqref{Sec3 subsolution-1-1} on the compact set $\overline{\Omega_1^{\eta_1}}\times\overline{\Omega_2^{\eta_2}}$ and the bound $(w^{\eta_2}_\varphi)_i\le\bar u_i$, there exists a constant $c_0>0$, independent of $\varphi$ and $\epsilon$, such that $\varrho_{\eta_1}(x^{(1)})\ge c_0$. Hence every such point lies in the open set
\[
U:=\bigl\{x^{(1)}\in\Omega_1:\varrho_{\eta_1}(x^{(1)})>c_0/2\bigr\}\times\Omega_2^{\eta_2}\times(t_2,\infty).
\]

Choose $\epsilon(t)$ in \eqref{lowersolution-2-PD-1-1} large enough so that $A_j^{(1)}(t)\varrho_{\eta_1}\le\epsilon(t)\varrho_{\eta_1}$ on $\{x\in\overline\Omega:\varrho_{\eta_1}(x^{(1)})>c_0/2\}\times[t_2,\infty)$ for $1\le j\le m$. Since $w^{\eta_2}_\varphi\gg 0$ in $\Omega_2^{\eta_2}\times(t_2,\infty)$, we obtain from \eqref{Omega-2-eta-problem-2} and \ref{ass:b4} that, in $U$,
\begin{align*}
	&\left(\frac{\partial}{\partial t}+A^{(1)}(t)+A^{(2)}(t)\right)\bigl(\varrho_{\eta_1}w^{\eta_2}_\varphi\bigr)-f(x^{(2)},t,\varrho_{\eta_1}w^{\eta_2}_\varphi)
	\\
	&\quad= \varrho_{\eta_1}\frac{\partial w^{\eta_2}_\varphi}{\partial t}+A^{(1)}(t)\bigl(\varrho_{\eta_1}w^{\eta_2}_\varphi\bigr)+\varrho_{\eta_1}A^{(2)}(t)w^{\eta_2}_\varphi-f(x^{(2)},t,\varrho_{\eta_1}w^{\eta_2}_\varphi)
	\\
	&\quad= \varrho_{\eta_1}f(x^{(2)},t,w^{\eta_2}_\varphi)-\epsilon(t)\varrho_{\eta_1}w^{\eta_2}_\varphi+A^{(1)}(t)\bigl(\varrho_{\eta_1}w^{\eta_2}_\varphi\bigr)-f(x^{(2)},t,\varrho_{\eta_1}w^{\eta_2}_\varphi)
	\\
	&\quad\le \varrho_{\eta_1}f(x^{(2)},t,w^{\eta_2}_\varphi)-f(x^{(2)},t,\varrho_{\eta_1}w^{\eta_2}_\varphi)\le 0.
\end{align*}
Thus $\varrho_{\eta_1}w^{\eta_2}_\varphi$ is a classical subsolution in $U$, whereas $\underline u=\underline v$ on the open set
\[
V:=\bigl\{(x,t)\in\Omega\times(t_2,\infty):\varrho_{\eta_1}(x^{(1)})(w^{\eta_2}_\varphi)_i(x^{(2)},t)<\underline v_i(x,t)\ \text{for all }1\le i\le m\bigr\}.
\]
Since $\Omega\times(t_2,\infty)=U\cup V$, Proposition~\ref{Sec3 subsolution-1} implies that $\underline u$ is a generalized subsolution of system \eqref{asymptotically partial degeneracy equations} in $\Omega\times(t_2,\infty)$.
\end{proof}

\begin{proposition}\label{lowersolution-4-PD}
Assume that $\lambda_\omega^*<0$ and \eqref{Sec3 auxiliary function-1-2} holds.
Let $K$ be any compact subset of $\Omega_1$.
Then for every solution $u(x,t)$ of system \eqref{asymptotically partial degeneracy equations}
with initial data $u(\cdot,0)\in C_0(\overline\Omega,\mathbb{M}_{\bar{u}})\setminus\{0\}$,
$\lim_{t\to+ \infty}|u(x,t)-w_\omega^*(x^{(2)},t)|=0 \ \ {\rm uniformly\ for} \ \ (x^{(1)},x^{(2)})\in K\times\overline\Omega_2$.
\end{proposition}

\begin{proof}
Fix a nonnegative, nontrivial solution $u(x,t)$ of system
\eqref{asymptotically partial degeneracy equations}, and choose a smooth domain
$\widetilde K$ containing $K$ whose closure is a compact subset of $\Omega_1$.
Apply the constructions of Propositions~\ref{Sec3 subsolution-1}
and~\ref{lowersolution-3-PD} on $\widetilde K\times\Omega_2$.
Choose $\eta_1$ sufficiently small that $\varrho_{\eta_1}=1$ on $K$,
and fix $\eta_2\in(0,\eta_0)$.
By the strong maximum principle, the Hopf boundary lemma on
$\overline{\widetilde K}\times\partial\Omega_2$,
and \eqref{Sec3 subsolution-1-7}, we may choose $\delta_1>0$
sufficiently small that
\[
u(x,t_2)\ge\underline v(x,t_2)
=\delta_1\phi(x^{(1)})
\bigl(q(t_2)\phi_\omega^*(x^{(2)},t_2)-\delta_2\Phi(x^{(2)})\bigr)
\quad\text{in }\widetilde K\times\Omega_2.
\]
With $\delta_1$ fixed, take the nonzero initial datum $\varphi$ in
Proposition~\ref{lowersolution-3-PD} sufficiently small that
$\varrho_{\eta_1}\varphi\le\underline v(\cdot,t_2)$.
Then $\underline u(\cdot,t_2)=\underline v(\cdot,t_2)\le u(\cdot,t_2)$.
By the comparison principle on $\widetilde K\times\Omega_2$, we obtain
\begin{align}\label{lowersolution-4-PD-1}
	u(x,t+k\omega)
	&\ge \underline u(x,t+k\omega)
	\ge \varrho_{\eta_1}(x^{(1)})
	w^{\eta_2}_\varphi(x^{(2)},t+k\omega)
	= w^{\eta_2}_\varphi(x^{(2)},t+k\omega)
\end{align}
for $x\in K\times\overline\Omega_2$, $t\in[0,\omega]$ and
$t+k\omega\ge t_2$, using the zero extension of $w^{\eta_2}_\varphi$.
Given $\varepsilon>0$, choose $\eta_2$ sufficiently small that
\[
\|w_\omega^{\eta_2}-w_\omega^*\|_{C(\overline\Omega_2\times[0,\omega])}
<\varepsilon/2
\]
by Lemma~\ref{convergence-domain perturbation}. For this $\eta_2$,
\eqref{asymptotic-Omega-2-eta-problem-2}, the periodicity of
$w_\omega^{\eta_2}$ and \eqref{lowersolution-4-PD-1} give
\begin{equation}\label{lowersolution-4-PD-2}
	u(x,t+k\omega)\ge
	w_\omega^*(x^{(2)},t)-\varepsilon(1,\dots,1)^\top
\end{equation}
for all $x\in K\times\overline\Omega_2$, all $t\in[0,\omega]$
and all sufficiently large $k$.

Now, let $v_{\overline\varphi}$ denote the solution of
\begin{equation*}
	\left\{
	\begin{aligned}
		& \frac{\partial v}{\partial t}+A^{(2)}(t)v=f(x^{(2)},t,v) \ \ \ &{\rm in}&\ \ \ \Omega_2\times(t_2,\infty),\\
		& v =0 \ \ \ \ &{\rm on}&\ \ \ \partial\Omega_2\times(t_2,\infty),\\
		& v(\cdot,t_2) = \overline\varphi \ \ \ \ &{\rm in}&\ \ \ \Omega_2,
	\end{aligned}
	\right.
\end{equation*}
where
$\overline\varphi(x^{(2)})
=\max_{x^{(1)}\in\overline\Omega_1}u(x,t_2)$
componentwise. Since $v_{\overline\varphi}$ is a supersolution of
system \eqref{asymptotically partial degeneracy equations},
the comparison principle gives
\begin{equation}\label{lowersolution-4-PD-3}
	u(x,t+k\omega)\le v_{\overline\varphi}(x^{(2)},t+k\omega)
\end{equation}
for $x\in\overline\Omega$, $t\in[0,\omega]$ and $t+k\omega\ge t_2$.
By Theorem~\ref{general threshold dynamics} on $\Omega_2$, we obtain
\[
\lim_{k\to\infty}\sup_{t\in[0,\omega]}
\bigl\|v_{\overline\varphi}(\cdot,t+k\omega)-w_\omega^*(\cdot,t)
\bigr\|_{C(\overline\Omega_2)}=0.
\]
Together with \eqref{lowersolution-4-PD-2} and
\eqref{lowersolution-4-PD-3}, this proves the asserted uniform convergence.
\end{proof}

If $\lambda_\omega^*\ge0$, the same upper comparison as in
\eqref{lowersolution-4-PD-3}, together with
Theorem~\ref{general threshold dynamics}, shows that every solution
$u(x,t)$ of system~\eqref{asymptotically partial degeneracy equations}
with $u(\cdot,0)\in C_0(\overline\Omega,\mathbb M_{\bar u})$
satisfies $u(x,t)\to0$ uniformly for $x\in\overline\Omega$
as $t\to\infty$.
As in \cite[Remark~3.7]{Lin-Ye-Zhao2026}, the comparison principle
then implies that the zero solution is also Lyapunov stable.

\begin{theorem}\label{threshold dynamics-PD1} The following statements are valid:
\begin{itemize}
\item [(a)] If $\lambda_\omega^{\ast}\ge0$, then the zero solution is globally asymptotically stable
		for system \eqref{asymptotically partial degeneracy equations} in $C_0(\overline\Omega,\mathbb{M}_{\bar{u}})$.
\item [(b)] If $\lambda_\omega^{\ast}<0$ and there exists $T>0$ such that
            \begin{equation}\label{threshold dynamics-PD1-1}
            \limsup_{t\to\infty}\frac{1}{T}\int_t^{t+T}\frac{\beta(s)^2}{\alpha(s)}\,\mathrm{d}s<4\lambda_0(-\lambda_\omega^*),
            \end{equation}
            then every nontrivial solution $u$ of system \eqref{asymptotically partial degeneracy equations} with $u(\cdot,0)\in C_0(\overline\Omega,\mathbb{M}_{\bar{u}})$ satisfies
                    $$
                    \lim\limits_{t\to+ \infty}\sup_{(x^{(1)},x^{(2)})\in K\times\overline\Omega_2}|u(x,t)-w_\omega^*(x^{(2)},t)|=0
                    $$
            for every compact subset $K$ of $\Omega_1$.
\end{itemize}
\end{theorem}

In the case where  the condition \eqref{threshold dynamics-PD1-1} does not hold, the asymptotic behavior of solutions becomes
more delicate. To address such a case,  we use the following assumptions \ref{ass:b4'} and \ref{ass:b5'}
to replace \ref{ass:b4} and \ref{ass:b5}, respectively, and impose an additional drift condition \ref{ass:b6}.

\begin{enumerate}[label=(B\arabic*)$'$,ref=(B\arabic*)$'$,leftmargin=*,start=4]
\item\label{ass:b4'}
$f_u(x^{(2)},t,u)$ is cooperative and irreducible for all $x^{(2)}\in\overline{\Omega}_2$, $t\ge0$ and $u\in\mathbb{R}^{m}_{+}$;
for each $1 \le i \le m$, $f_i(\cdot,\cdot,u)\ge0$ whenever $u\in\mathbb R^m_+$ with $u_i = 0$;
$f(\cdot,\cdot,0) = 0$; $\alpha f(\cdot,\cdot,u)\le f(\cdot,\cdot,\alpha u)$ for all  $u\in\mathrm{Int}(\mathbb{R}^{m}_{+})$ and $\alpha\in (0,1)$.

\item\label{ass:b5'}
$\tilde f_u(x^{(2)},t,u)$ is cooperative and irreducible for all $x^{(2)}\in\overline{\Omega}_2$, $t\ge0$ and $u\in\mathbb{R}^{m}_{+}$;
for each $1 \le i \le m$, $\tilde f_i(\cdot,\cdot,u)\ge0$ whenever $u\in\mathbb R^m_+$ with $u_i = 0$;
$\tilde f(\cdot,\cdot,0) = 0$; $\alpha\tilde f(\cdot,\cdot,u)<\tilde f(\cdot,\cdot,\alpha u)$ for all  $u\in\mathrm{Int}(\mathbb{R}^{m}_{+})$ and $\alpha\in (0,1)$.
\end{enumerate}
\begin{enumerate}[label=(B\arabic*),ref=(B\arabic*),leftmargin=*,start=6]
\item\label{ass:b6}
There exist a unit vector $e\in\mathbb{R}^{n_1}$, $t_0\ge0$ and a positive function $\chi$ such that
$$
\sum_{j=1}^{n_1}b^{(i)}_{j}(x,t)e_j\ge\chi(t), \quad x\in\overline\Omega,\ t\ge t_0, \ i=1,\dots,m.
$$
\end{enumerate}
Here $\chi(t)$ represents a uniform positive lower bound for the projection of each species' drift field onto the direction $e$.
For convenience, we define $p(t):=\frac{\chi(t)}{\alpha(t)}$ for $t\ge t_0$.

\begin{theorem}\label{threshold dynamics-PD2}
Assume that $\lambda_\omega^{\ast}<0$. If $p$ is locally absolutely continuous on $[t_0,\infty)$ and there exists $T>0$ such that
\begin{equation}\label{threshold dynamics-PD2-1}
\liminf_{t\to\infty}\frac{1}{T}\int_{t}^{t+T}\left(\frac{\chi(s)^2}{\alpha(s)}-2d_e\bigl(p'(s)\bigr)_+\right)\,\mathrm{d}s>4\Lambda_0(-\lambda_\omega^{\ast}),
\end{equation}
where $(x)_+=\max\{x,0\}$, and $d_e=\max_{x^{(1)}\in\overline{\Omega}_1} x^{(1)}\cdot e-\min_{x^{(1)}\in\overline{\Omega}_1} x^{(1)}\cdot e$,
then every nonnegative solution of system \eqref{asymptotically partial degeneracy equations} converges to zero uniformly on $\overline\Omega$ as $t\to\infty$.
\end{theorem}

\begin{proof}
By \ref{ass:b4'} and \ref{ass:b5'}, and the comparison principle, it suffices to construct a supersolution that converges to zero uniformly in $\overline\Omega$ as $t\to\infty$. We then proceed in two steps.

\medskip
{\noindent\it Step 1.} 
Define $\Phi(x^{(1)},t):=e^{\frac{1}{2\Lambda_0}p(t)\xi(x^{(1)})}$ for $(x^{(1)},t)\in\overline\Omega_1\times[t_0,\infty)$,
where $\xi(x^{(1)}):=x^{(1)}\cdot e-\max_{x^{(1)}\in\overline{\Omega}_1} x^{(1)}\cdot e$. Then $\Phi(x^{(1)},t)\in(0,1]$ for $x^{(1)}\in\overline\Omega_1$ and $t\ge t_0$.
A direct computation gives
$$
\frac{A_i^{(1)}(t)\Phi}{\Phi}=-\frac{p(t)^2}{4\Lambda_0^2}\sum_{j,k=1}^{n_1} a^{(i)}_{jk}(x,t)e_je_k+\frac{p(t)}{2\Lambda_0}\sum_{j=1}^{n_1} b^{(i)}_j(x,t)e_j.
$$
By \ref{ass:b1} and \ref{ass:b6}, we obtain
$$
\frac{A_i^{(1)}(t)\Phi}{\Phi}\ge-\frac{p(t)^2}{4\Lambda_0}\alpha(t)+\frac{p(t)}{2\Lambda_0}\chi(t)=\frac{\chi(t)^2}{4\Lambda_0\alpha(t)}.
$$
In addition, noting that $p'(t)\xi(x^{(1)})\ge -d_e\bigl(p'(t)\bigr)_+$ for $x\in\overline\Omega$ and almost every $t\ge t_0$, we have
\begin{equation}\label{threshold dynamics-PD2-3}
\frac{\Phi_t+A_i^{(1)}(t)\Phi}{\Phi}\ge\frac{\chi(t)^2}{4\Lambda_0\alpha(t)}-\frac{d_e}{2\Lambda_0}\bigl(p'(t)\bigr)_+.
\end{equation}

Define $\Omega_{2,\varepsilon}:=\left\{x^{(2)}\in\mathbb{R}^{n_2}: d(x^{(2)}, \Omega_2)<\varepsilon\right\}$ with $\varepsilon>0$.
It is clear that $\Omega_{2,\varepsilon}\downarrow\overline\Omega_2$ as $\varepsilon\downarrow0^{+}$, and $\partial\Omega_{2,\varepsilon}:=\left\{x^{(2)}+\varepsilon\nu(x^{(2)}); x^{(2)}\in\partial\Omega_2\right\}$
for sufficiently small $\varepsilon$, where $\nu(x^{(2)})$ is the unit outward normal vector on $\partial\Omega_2$. Thus, there exists small $\varepsilon_0>0$ such that the $C^{2+\theta}$-norm of $\partial\Omega_{2,\varepsilon}$ is uniformly bounded with respect to $\varepsilon\in(0,\varepsilon_0)$.
Extend the coefficients of the periodic operator $\tilde A_i^{(2)}(t)$ and $\tilde f_u(\cdot,\cdot,0)$
from $\overline\Omega_2$ to $\overline\Omega_{2,\varepsilon}$, preserving their $\omega$-periodicity,
H\"older regularity, symmetry of the diffusion matrices, uniform ellipticity, and the cooperativity
and irreducibility of $\tilde f_u(\cdot,\cdot,0)$.
For convenience, the extended coefficients are still denoted by the same symbols.
We then consider the corresponding $\omega$-periodic parabolic eigenvalue problem on the enlarged domain $\overline\Omega_{2,\varepsilon}$, namely,
\begin{equation}\label{parabolic eigenvalue problem-PD-varepsilon}
\left\{
\begin{aligned}
& \frac{\partial w}{\partial t}+\tilde A^{(2)}(t)w=\tilde f_u(x^{(2)},t,0)w+\lambda w \ \ \ &{\rm in}&\ \ \ \Omega_{2,\varepsilon}\times\mathbb{R},\\
& w=0 \ \ \ \ &{\rm on}&\ \ \ \partial\Omega_{2,\varepsilon}\times\mathbb{R},\\
& w(x^{(2)},\omega+t)=w(x^{(2)},t) \ \ \ &{\rm in}&\ \ \ \Omega_{2,\varepsilon}\times\mathbb{R}.
\end{aligned}
\right.
\end{equation}
Let $\lambda^*_{\omega,\varepsilon}$ denote the principal eigenvalue of \eqref{parabolic eigenvalue problem-PD-varepsilon}. Then
\begin{equation}\label{threshold dynamics-PD2-4}
\lim_{\varepsilon\to0^+}\lambda^*_{\omega,\varepsilon}=\lambda^*_{\omega}.
\end{equation}
Since $\lambda^*_{\omega}<0$, we may choose $\varepsilon_0>0$
sufficiently small such that $\lambda^*_{\omega,\varepsilon}<0$ for all $\varepsilon\in(0,\varepsilon_0]$.
Denote by $\phi^*_{\omega,\varepsilon}$ the corresponding positive principal eigenfunction.
Then $\phi^*_{\omega,\varepsilon}\in C^{2,1}(\overline\Omega_{2,\varepsilon}\times\mathbb{R},\mathbb{R}^m)$ and $\phi^*_{\omega,\varepsilon}\gg0$ in $\overline\Omega_{2}\times\mathbb{R}$.
For $t\ge t_0$, denote
$$
\Gamma_\varepsilon(t):=\max_{1\le i\le m}\sup_{x^{(2)}\in\overline\Omega_2}
\frac{\left|\left[\left(f_u(x^{(2)},t,0)-\tilde f_u(x^{(2)},t,0)\right)\phi^*_{\omega,\varepsilon}-\left(A^{(2)}(t)-\tilde A^{(2)}(t)\right)\phi^*_{\omega,\varepsilon}\right]_i\right|}{(\phi^*_{\omega,\varepsilon})_i}
$$
and
$$
\mu_\varepsilon(t):=\exp\left\{-\int_{t_0}^{t}\left(\frac{\chi(s)^2}{4\Lambda_0\alpha(s)}+\lambda^*_{\omega,\varepsilon}-\frac{d_e}{2\Lambda_0}\bigl(p'(s)\bigr)_+-\Gamma_\varepsilon(s)\right)\,\mathrm{d}s\right\}.
$$
By \ref{ass:b3}, together with the positivity and periodicity of $\phi^*_{\omega,\varepsilon}$, we have $\lim_{t\to+\infty}\Gamma_\varepsilon(t)=0$.
Hence, in view of \eqref{threshold dynamics-PD2-1} and \eqref{threshold dynamics-PD2-4}, we may, after decreasing $\varepsilon_0>0$
if necessary, assume that for every $\varepsilon\in(0,\varepsilon_0]$
$$
\int_{t_0}^{t}\left(\frac{\chi(s)^2}{4\Lambda_0\alpha(s)}+\lambda^*_{\omega,\varepsilon}-\frac{d_e}{2\Lambda_0}\bigl(p'(s)\bigr)_+-\Gamma_\varepsilon(s)\right)\,\mathrm{d}s\to+\infty \ \ {\rm as}\ \ t\to+\infty.
$$
Consequently, $\lim_{t\to+\infty}\mu_\varepsilon(t)=0$ for every $\varepsilon\in(0,\varepsilon_0]$.
Now, we define
\begin{equation}\label{threshold dynamics-PD2-5}
\overline{u}(x,t):=\mu_\varepsilon(t)\Phi(x^{(1)},t)\phi^*_{\omega,\varepsilon}(x^{(2)},t), \quad (x,t)\in\overline\Omega\times[t_0,\infty).
\end{equation}
It is clear that $\overline{u}(\cdot,t)\to0$ uniformly in $\overline\Omega$ as $t\to\infty$.

\medskip
{\noindent\it Step 2.} 
By \eqref{threshold dynamics-PD2-3} and \eqref{parabolic eigenvalue problem-PD-varepsilon}, it follows that for each \(i=1,\dots,m\) and each $x\in\overline\Omega$ and almost every $t\ge t_0$,
\begin{align*}
    &\frac{(\overline u_i)_t+A^{(1)}_i(t)\overline u_i+A^{(2)}_i(t)\overline u_i}{\mu_\varepsilon\Phi}
    \\
    &  \ge  \frac{\mu_\varepsilon'}{\mu_\varepsilon}(\phi^*_{\omega,\varepsilon})_i+\left(A^{(2)}_i(t)-\tilde A^{(2)}_i(t)\right)(\phi^*_{\omega,\varepsilon})_i+\left(\tilde f_u(x^{(2)},t,0)\phi^*_{\omega,\varepsilon}\right)_i
    \\
    & {\hspace{17pt}} +\lambda^*_{\omega,\varepsilon} (\phi^*_{\omega,\varepsilon})_i+\frac{\chi(t)^2}{4\Lambda_0\alpha(t)}(\phi^*_{\omega,\varepsilon})_i-\frac{d_e}{2\Lambda_0}\bigl(p'(t)\bigr)_+(\phi^*_{\omega,\varepsilon})_i
    \\
    &  \ge  \frac{\mu_\varepsilon'}{\mu_\varepsilon}(\phi^*_{\omega,\varepsilon})_i-\Gamma_\varepsilon(t)(\phi^*_{\omega,\varepsilon})_i+\Big(f_u(x^{(2)},t,0)\phi^*_{\omega,\varepsilon}\Big)_i+\lambda^*_{\omega,\varepsilon} (\phi^*_{\omega,\varepsilon})_i
    \\
    & {\hspace{17pt}} +\frac{\chi(t)^2}{4\Lambda_0\alpha(t)}(\phi^*_{\omega,\varepsilon})_i-\frac{d_e}{2\Lambda_0}\bigl(p'(t)\bigr)_+(\phi^*_{\omega,\varepsilon})_i
    \\
    &  =  \Big(f_u(x^{(2)},t,0)\phi^*_{\omega,\varepsilon}\Big)_i.
\end{align*}
Combining this with the subhomogeneity of $f$, we obtain
$$
\overline u_t+A^{(1)}(t)\overline u+A^{(2)}(t)\overline u\ge f(x^{(2)},t,\overline u),\quad (x,t)\in\Omega\times(t_0,\infty).
$$
This, together with the boundary condition,  shows that $\overline u$ is a supersolution of system \eqref{asymptotically partial degeneracy equations} in $\Omega\times(t_0,\infty)$.
\end{proof}

Note that in the proof of Theorem \ref{threshold dynamics-PD2}, the construction of the supersolution $\overline u$ is independent of the sign of $\lambda_\omega^{\ast}$. We retain the assumption $\lambda_\omega^{\ast}<0$ solely for a clear comparison with Theorem \ref{threshold dynamics-PD1}(b).

In the rest of this section, we present an example to illustrate the sharpness of both the upper bound $4\lambda_0(-\lambda_\omega^*)$ in \eqref{threshold dynamics-PD1-1}
and the lower bound $4\Lambda_0(-\lambda_\omega^*)$ in \eqref{threshold dynamics-PD2-1}.
\begin{example}\label{partial-degenerate-example}
\rm
Fix $r\in\mathbb{R}$.
Consider the following two-dimensional problem:
\begin{equation}\label{Example-PD}
\begin{cases}
u_t-\frac{1}{1+t}u_{xx}-u_{yy}+\frac{r}{\sqrt{1+t}}u_x=u(1-u), & -4<x,y<4, \quad t>0,\\
u(x,\pm4,t)=u(\pm4,y,t)=0, & -4<x,y<4, \quad t>0,\\
u(x,y,0)=u_0(x,y)\ge0, & -4\le x,y\le4,
\end{cases}
\end{equation}
where $u_0\in C_0([-4,4]^2)$, and $\lambda_0$, $\Lambda_0$, $d_e$, $\alpha$, $\beta$, and $\chi$ are  given by
$$
\lambda_0=\Lambda_0=1, \ \ \ d_e=8,  \ \ \  \alpha(t)=\frac{1}{1+t},  \ \ \ \chi(t)=\beta(t)=\frac{|r|}{\sqrt{1+t}}.
$$
Let $p(t):=\frac{\beta(t)}{\alpha(t)}=|r|\sqrt{1+t}$. Then $p\in C^1([0,\infty))$ and $p'(t)=\frac{|r|}{2\sqrt{1+t}}\to0$ as $t\to\infty$.

Note that the formal limiting problem associated with \eqref{Example-PD} becomes
\begin{equation}\label{Example-PD-1}
\begin{cases}
v_t-v_{yy}=v(1-v), & -4<y<4, \quad t>0,\\
v(\pm4,t)=0, & t>0.
\end{cases}
\end{equation}
It is easy to verify that the principal eigenvalue associated with the linearization of \eqref{Example-PD-1} at $v=0$ is
$$
\lambda^*=\frac{\pi^2}{64}-1\thickapprox-0.8458<0.
$$
By Theorem \ref{general threshold dynamics}, system \eqref{Example-PD-1} admits a unique 
positive steady state $v^*$,
and every nonnegative nontrivial solution $v$ of it  converges to $v^*$ uniformly on $[-4,4]$ as $t\to\infty$.
Further,  the conditions \eqref{threshold dynamics-PD1-1} and \eqref{threshold dynamics-PD2-1} in Theorems
\ref{threshold dynamics-PD1} and \ref{threshold dynamics-PD2} reduce to
$r^2<\frac{64-\pi^2}{16}$ and $r^2>\frac{64-\pi^2}{16}$, respectively.
We then obtain the following threshold-type result for system \eqref{Example-PD}.

\begin{proposition}\label{threshold dynamics-PD3}
The following statements are valid for system \eqref{Example-PD}:
	\begin{itemize}
		\item [(a)] If $|r|>\frac{\sqrt{64-\pi^2}}{4}$, then every solution $u$ satisfies
		$\lim\limits_{t\to+ \infty}\sup_{(x,y)\in [-4,4]^2}u(x,y,t)=0$.
		\item [(b)] If $|r|<\frac{\sqrt{64-\pi^2}}{4}$, then every nontrivial solution $u$
		satisfies  $$\lim\limits_{t\to+ \infty}\sup_{(x,y)\in K\times[-4,4]}|u(x,y,t)-v^*(y)|=0$$ for every compact set $K\subset(-4,4)$.
	\end{itemize}
\end{proposition}

The following numerical simulations illustrate the threshold dynamics for system \eqref{Example-PD}
with the initial value being $u_0(x,y)=\max\{(1-x^2),0\}\max\{(1-y^2),0\}$.
As shown in Fig. 1, the solution exhibits extinction or persistence depending on the value of $|r|$, which is
in agreement with Proposition \ref{threshold dynamics-PD3}.

\captionsetup{labelsep=period}
\renewcommand{\figurename}{Fig.}
\begin{figure}[htbp]
	\centering
	\subfigure[$r=1.5<\frac{\sqrt{64-\pi^2}}{4}\thickapprox1.839$.]{
		\begin{minipage}{6.5cm}
			\centering
			\includegraphics[scale=0.26]{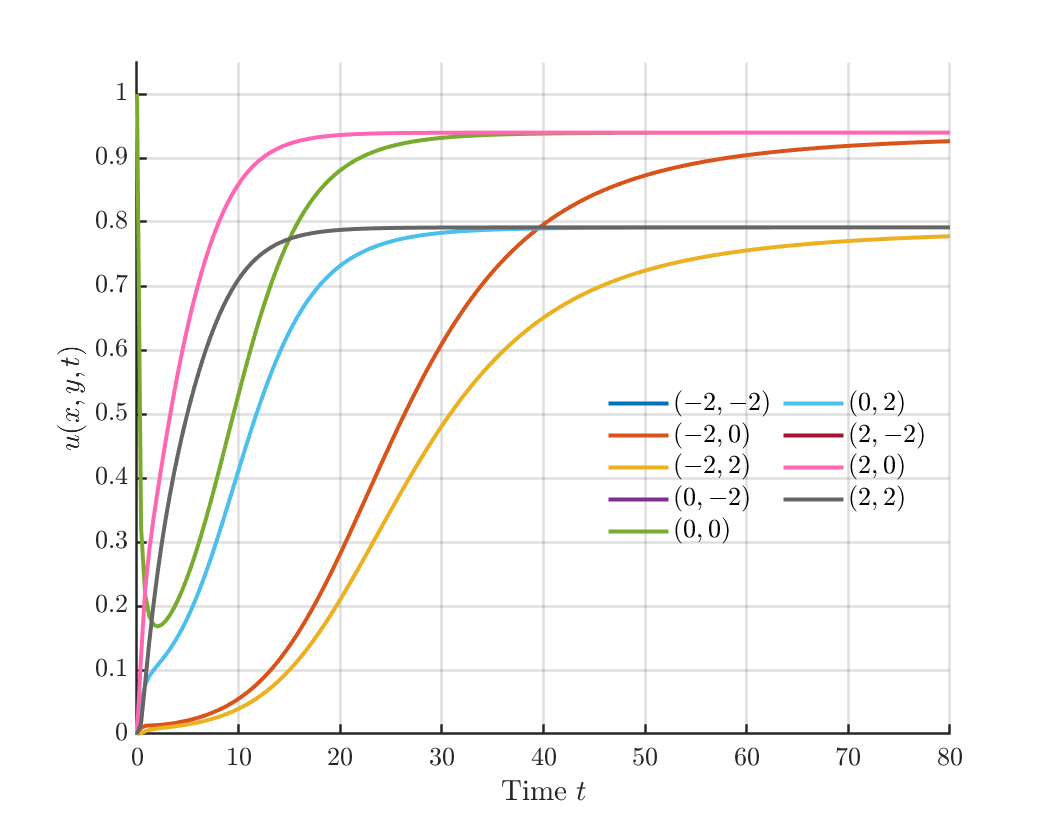}
		\end{minipage}
	}
	\hspace{0.09\textwidth}
	\subfigure[$r=2.5>\frac{\sqrt{64-\pi^2}}{4}\thickapprox1.839$.]{
		\begin{minipage}{6.5cm}
			\centering
			\includegraphics[scale=0.26]{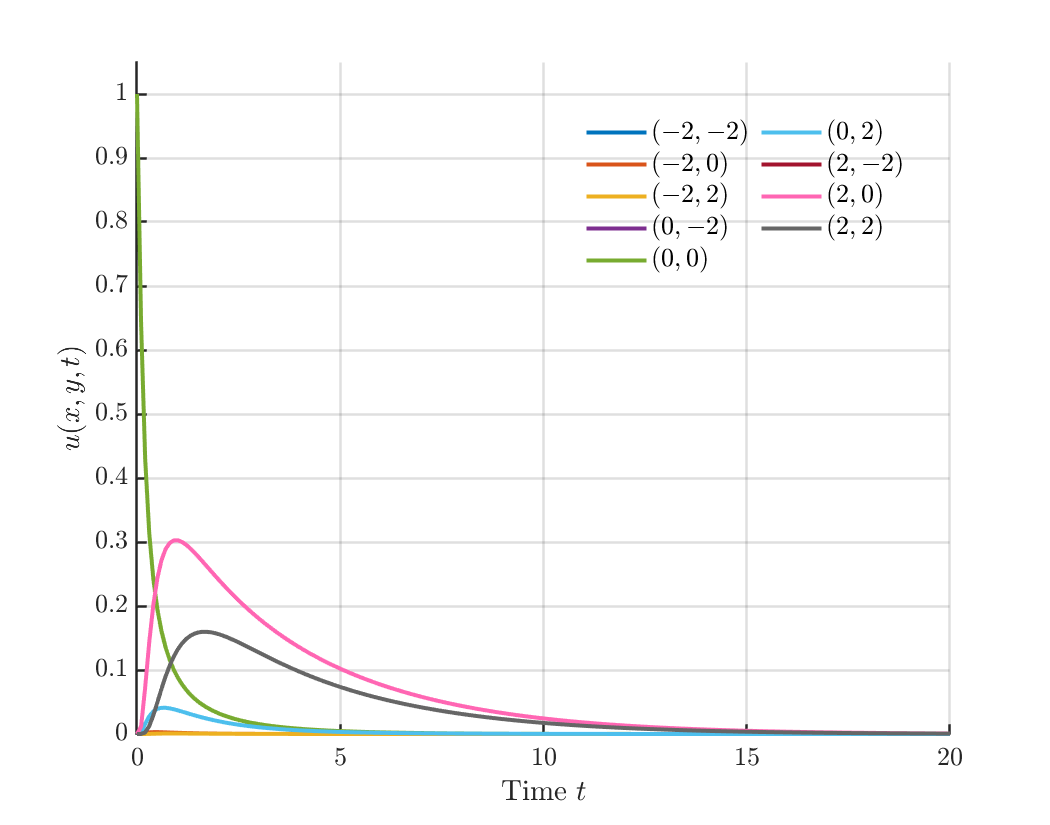}

		\end{minipage}
	}
	\caption{Temporal dynamics of the solution $u$ at different spatial points.}
	\label{fig:fixed}
\end{figure}
\end{example}

\section{Asymptotically fully degenerate regime}
\label{sec:asymptotically fully degenerate}
Let $\Omega\subset\mathbb{R}^{n}$ be a bounded domain with boundary of class $C^{2+\theta}$ for some $0<\theta<1$.
Consider the following nonautonomous parabolic system on $\Omega$:
\begin{equation}\label{asymptotically degenerate equations}
\left\{
\begin{aligned}
& \frac{\partial v}{\partial t}   +A(t)v=f(t,v) \ \ \ &{\rm in}&\ \ \ \Omega\times(0,\infty),\\
& v =0 \ \ \ \ &{\rm on}&\ \ \ \partial\Omega\times(0,\infty),
\end{aligned}
\right.
\end{equation}
and an ordinary differential system:
\begin{equation}\label{limiting ODE system}
\dfrac{\mathrm{d} w}{\mathrm{d} t}=\tilde f(t,w), \qquad t>0,
\end{equation}
where $A(t)v=(A_1(t)v_1,\cdots,A_m(t)v_m)^\top$, $f(t,v)=(f_1(t,v),\cdots,f_m(t,v))^\top$, $\tilde f(t,v)=(\tilde f_1(t,v), \cdots, \tilde f_m(t,v))^\top$,
and
\begin{align*}
A_i(t)v_i&=-\sum_{j,k=1}^{n}a^{(i)}_{jk}(x,t)\partial_{jk}v_i +\sum_{j=1}^{n}b^{(i)}_{j}(x,t)\partial_jv_i, \ \ i=1,\dots,m.
\end{align*}
The aim of this section is to establish the evolution dynamics of system \eqref{asymptotically degenerate equations}
under the assumption that for each $i$, the operator $A_i(t)$ converges coefficientwise to zero as
$t\to\infty$.
Throughout this section, we impose the following assumptions:
\begin{enumerate}[label=(C\arabic*),ref=(C\arabic*),leftmargin=*]
\item\label{ass:c1}
$a^{(i)}_{jk}$, $b^{(i)}_{j}\in C^{\theta, \theta/2}(\overline\Omega\times[0,T], \mathbb{R})$
for every $T>0$,
$a^{(i)}_{jk}=a^{(i)}_{kj}$, and there exist positive constants $\lambda_0$, $\Lambda_0$, and a positive continuous function $\alpha$ such that
\begin{equation*}
\alpha(t)\lambda_0|\xi|^2\le\sum_{j,k=1}^{n}a^{(i)}_{jk}(x,t)\xi_j\xi_k\le \alpha(t)\Lambda_0|\xi|^2, \quad (x,t)\in \Omega\times\mathbb{R}_+,
\end{equation*}
for any $\xi=(\xi_1,\dots,\xi_{n})^\top\in\mathbb{R}^n$ and $1\le i\le m$.

\item\label{ass:c2}
For $1\le i,j\le m$,
$f_i,\tilde f_i,\frac{\partial f_i}{\partial u_j},
\frac{\partial\tilde f_i}{\partial u_j}
\in C(\mathbb{R}_+\times\mathbb{R}^m,\mathbb{R})$.
For every $T>0$,
$f_i(\cdot,u)\in C^{\theta/2}([0,T],\mathbb{R})$
uniformly for $u$ in bounded subsets of $\mathbb{R}^m$.

\item\label{ass:c3}
For $\beta(t):=\sqrt{n}\sup\limits_{1\le i\le m, 1\le j\le n}\|b^{(i)}_{j}(\cdot,t)\|_{C(\overline\Omega)}$,
there holds $\lim_{t \to \infty}\big(\alpha(t)+\beta(t)\big) = 0$,  and
\begin{align*}
    & \lim_{t \to \infty}\left(\left\|f(t,u) -\tilde f(t,u) \right\|_{\mathbb{R}^m}
    + \left\| f_u(t, 0) - \tilde{f}_u(t, 0) \right\|_{\mathbb{R}^{m\times m}} \right)= 0
\end{align*}
uniformly for $u$ in bounded subsets of $\mathbb R^m$,
where $\tilde f(\cdot,u)$ is $\omega$-periodic on $\mathbb R$
for each $u\in\mathbb R^m$.

\item\label{ass:c4}
$f_u(t,u)$ is cooperative and irreducible for $t\ge0$ and $u\in\mathbb{M}_{\bar{u}}$.
For each $1 \le i \le m$ and any $t\ge0$, $f_i(t,u)\ge0$ whenever $u\in\mathbb{M}_{\bar{u}}$ with $u_i = 0$,
and $f_i(t,u)\le 0$ whenever $u \in \mathbb{M}_{\bar{u}}$ with $u_i = \bar{u}_i$;
$f(t,0) = 0$;  $\alpha f(t,u)\le f(t,\alpha u)$ for all  $u\in\mathbb{M}_{\bar{u}}\cap\mathrm{Int}(\mathbb{R}^{m}_{+})$, $t\ge0$ and $\alpha\in (0,1)$.

\item\label{ass:c5}
$\tilde f_u(t,u)$ is cooperative and irreducible for $t\ge0$ and $u\in\mathbb{M}_{\bar{u}}$.
$\alpha \tilde f(t,u)<\tilde f(t,\alpha u)$ for all $u\in\mathbb{M}_{\bar{u}}\cap\mathrm{Int}(\mathbb{R}^{m}_{+})$, $t\ge0$ and $\alpha\in (0,1)$.
\end{enumerate}

We first introduce the following auxiliary ordinary differential  system:
\begin{equation}\label{limit system-asymptotically degenerate equations}
\left\{
{\begin{array}{*{20}{l}}
	\dfrac{\mathrm{d} w}{\mathrm{d} t}=f(t,w), & t>0,
        \\[3mm]
	w(0)=w_{0}\in\mathbb{M}_{\bar{u}}.
\end{array}}
\right.
\end{equation}
Let $(\lambda^*, z)$ be  the principal eigenpair  of the following $\omega$-periodic eigenvalue problem:
\begin{equation}\label{eigenvalue problem-limit system-asymptotically degenerate equations}
\left\{
{\begin{array}{*{20}{l}}
	\dfrac{\mathrm{d} w}{\mathrm{d} t}=\tilde f_u(t,0)w+\lambda w, & t\in\mathbb{R},
        \\[3mm]
	w(t+\omega)=w(t), & t\in\mathbb{R}.
\end{array}}
\right.
\end{equation}
Then $z(t)\gg0$ for every $t\in[0,\omega]$.
By the theory of asymptotically periodic semiflows (see \cite[Section 3.2]{MR3643081}),
we have the following result.

\begin{lemma}\label{limit system-threshold dynamics}
The following statements are valid:
\begin{itemize}
	\item [(1)] If $\lambda^*\ge0$, then the zero solution is globally asymptotically stable for
                system \eqref{limit system-asymptotically degenerate equations}	in $\mathbb{M}_{\bar{u}}$.
	\item [(2)] If $\lambda^*<0$, then system \eqref{limiting ODE system} admits a unique positive $\omega$-periodic solution $w_{e}$. For every $w_0\in\mathbb{M}_{\bar u}\setminus\{0\}$, the solution
	$w(t;w_0)$ of \eqref{limit system-asymptotically degenerate equations}
	satisfies $\lim_{t\to\infty}|w(t;w_0)-w_e(t)|=0$.
\end{itemize}
\end{lemma}

For each $\phi\in C_0(\overline\Omega,\mathbb{M}_{\bar{u}})$, let $(w_{0})_i:=\max_{x\in\overline\Omega} \phi_i(x)$.
According to the comparison principle, we have $u(x,t;\phi)\le w(t;w_{0})$. This leads to the  following result.

\begin{proposition}\label{threshold dynamics-1}
If $\lambda^*\ge 0 $, then for each $\phi\in C_0(\overline\Omega,\mathbb{M}_{\bar{u}})$, the solution $u(x,t;\phi)$ of system \eqref{asymptotically degenerate equations} satisfies
$\lim_{t\to +\infty}\left\|u(\cdot,t;\phi) \right\|_{C(\overline\Omega)}=0$.
\end{proposition}

In what follows, we focus on the global dynamics of system \eqref{asymptotically degenerate equations} in the case $\lambda^*<0$.
By the arguments similar  to those for Lemma \ref{Sec3 auxiliary function-1}, we have the following observation.

\begin{lemma}\label{Sec4 auxiliary function-1}
There exist a constant $t_1>0$, a function $\phi\in C^2(\overline{\Omega})$ and a positive function $H_\phi\in C([t_1,\infty))$ such that for each $i=1,\dots,m$,
\begin{equation}\label{Sec4 auxiliary function-1-1}
\left\{
\begin{aligned}
& A_i(t)\phi\le H_\phi(t)\phi \ \ \ &{\rm in}&\ \ \ \Omega\times(t_1,\infty),\\
& \phi>0 \ \ \ \ &{\rm in}&\ \ \ \Omega,\\
& \phi=0 \ \ \ \ &{\rm on}&\ \ \ \partial\Omega.
\end{aligned}
\right.
\end{equation}
Moreover, if $\lambda^*<0$ and there exists $T>0$ such that
\begin{equation}\label{Sec4 auxiliary function-1-2}
\limsup_{t\to\infty}\frac{1}{T}\int_t^{t+T}\frac{\beta(s)^2}{\alpha(s)}\,\mathrm{d}s<4\lambda_0(-\lambda^*),
\end{equation}
then $\phi$ and $H_\phi$ in \eqref{Sec4 auxiliary function-1-1} can be chosen so that
\begin{equation*}
\limsup_{t\to\infty}\frac{1}{T}\int_t^{t+T}H_\phi(s)\,\mathrm{d}s<-\lambda^*.
\end{equation*}
\end{lemma}

The following subsolution is the fully degenerate counterpart of the one constructed in Proposition \ref{Sec3 subsolution-1}.
The proof follows from similar arguments with only minor modifications, and is therefore omitted.

\begin{proposition}\label{Sec4 subsolution-1}
Assume that $\lambda^*<0$ and \eqref{Sec4 auxiliary function-1-2} holds.
Then there exist $t_2>t_1$ and a function $\underline v\in C^{2,1}(\overline\Omega\times[t_2,\infty),\mathbb{M}_{\bar{u}})$ such that
$\underline v$ is a subsolution of system \eqref{asymptotically degenerate equations} in $\Omega\times(t_2,\infty)$, and satisfies
$\underline v\gg0 \ \ {\rm in}\ \ \Omega\times[t_2,\infty)$ and
$\underline v=0 \ \ {\rm on}\ \ \partial\Omega\times[t_2,\infty)$.
Moreover, $\min_{1\le i\le m}\inf_{x\in K,\, t\ge t_2}\underline v_i(x,t)>0$ for every compact subset $K\subset\Omega$.
\end{proposition}

Let $\varrho_{\eta}$ be a smooth cut-off function satisfying $\varrho_{\eta}=0 \text{ in } \Omega\setminus\Omega_{\eta}$,
$\varrho_{\eta}=1$ in $\Omega_{2\eta}$, and $\varrho_{\eta}$ is convex in $\mathbb{R}^{n}\setminus\Omega^{\frac{3\eta}{2}}$.
Define
$\mathcal{W}=\left\{w\in C([t_2,\infty),\mathbb{M}_{\bar{u}}): w\gg0\ {\rm for}  \ t\in(t_2,\infty)\right\}$,
where $t_2$ is given in Proposition \ref{Sec4 subsolution-1}.
By Proposition \ref{Sec4 subsolution-1}, we may choose $\delta>0$ sufficiently small
so that
\begin{equation}\label{nonempty-2}
w_e(t)\gg\delta\underline v(x,t),\qquad \forall x\in\Omega, \ t\ge t_2.
\end{equation}
For $w\in\mathcal{W}$, define
$\mathcal{D}_w= \left\{(x,t)\in\Omega\times[t_2,\infty): \varrho_{\eta}(x)w(t)>\delta\underline v(x,t)\right\}$.
It follows from \eqref{nonempty-2} that $\mathcal{D}_{w}\ne\varnothing$ for some $w\in\mathcal{W}$.
Arguing as in Lemma \ref{lowersolution-2-PD} and Remark \ref{lowersolution-2-PD-1}, and using the asymptotic degeneracy of $A(t)$,
one can choose a positive function $\epsilon\in C^1([t_2,\infty))$ with $\epsilon(t)\downarrow0$ as $t\to\infty$, such that
for every $w\in\mathcal{W}$ with $\mathcal{D}_w\ne\emptyset$, the following estimate holds with the same $\epsilon$:
\begin{equation}\label{lowersolution-2-D-1-1}
A_i(t)\varrho_{\eta}(x)\le \epsilon(t)\varrho_{\eta}(x), \quad (x,t)\in\mathcal{D}_w,\ \ 1\le i\le m.
\end{equation}
In addition, let $w_\varphi$ denote the nontrivial solution of the following system
\begin{equation*}
\left\{
{\begin{array}{*{20}{l}}
	\dfrac{\mathrm{d} w}{\mathrm{d} t}=f(t,w)-\epsilon(t)w, & t>t_2,
        \\[3mm]
	w(t_2)=\varphi\in\mathbb{M}_{\bar{u}}\setminus\{0\}.
\end{array}}
\right.
\end{equation*}
If $\lambda^*<0$, then Lemma \ref{limit system-threshold dynamics} yields
$\lim_{t\to+ \infty} |w_\varphi(t)-w_e(t)|=0$.

Following the arguments of Proposition \ref{lowersolution-3-PD} and using \eqref{lowersolution-2-D-1-1},
we can further construct a generalized subsolution of \eqref{asymptotically degenerate equations},
which plays a key role in establishing the asymptotic behavior of solutions.

\begin{proposition}\label{Sec4 subsolution-2}
Assume that $\lambda^*<0$ and \eqref{Sec4 auxiliary function-1-2} holds. Then the componentwise maximum $\underline u(x,t):=\max\left\{\varrho_{\eta}(x)w_\varphi(t), \ \delta\underline v(x,t)\right\}$
is a generalized subsolution of system \eqref{asymptotically degenerate equations} in $\Omega\times[t_2,\infty)$.
\end{proposition}

By adapting the arguments for Proposition \ref{lowersolution-4-PD} and combining it with
Proposition \ref{threshold dynamics-1}, we obtain the following result.

\begin{theorem}\label{threshold dynamics-D}
The following statements are valid:
	\begin{itemize}
		\item [(a)] If $\lambda^*\ge0$, then the zero solution is globally asymptotically stable for system \eqref{asymptotically degenerate equations} in $C_0(\overline\Omega,\mathbb{M}_{\bar{u}})$.
		\item [(b)] If $\lambda^*<0$ and there exists $T>0$ such that
                        \begin{equation}\label{threshold dynamics-D-1}
                        \limsup_{t\to\infty}\frac{1}{T}\int_t^{t+T}\frac{\beta(s)^2}{\alpha(s)}\,\mathrm{d}s<4\lambda_0(-\lambda^*),
                        \end{equation}
                    then every nontrivial solution $u(x,t)$ of system \eqref{asymptotically degenerate equations}
                    with $u(\cdot,0)\in C_0(\overline\Omega,\mathbb{M}_{\bar{u}})$
                    satisfies $\lim\limits_{t\to+ \infty}|u(\cdot,t)-w_e(t)|=0$ locally uniformly in $\Omega$.
	\end{itemize}
\end{theorem}

In order to consider the case where the condition \eqref{threshold dynamics-D-1} does not hold, we
use the following assumptions \ref{ass:c4'} and \ref{ass:c5'} to replace \ref{ass:c4} and \ref{ass:c5}, respectively,
and impose an additional drift condition \ref{ass:c6}.

\begin{enumerate}[label=(C\arabic*)$'$,ref=(C\arabic*)$'$,leftmargin=*,start=4]
\item\label{ass:c4'}
$f_u(t,u)$ is cooperative and irreducible for $t\ge0$ and $u\in\mathbb{R}^{m}_{+}$.
For each $1 \le i \le m$ and any $t\ge0$, $f_i(t,u)\ge0$ whenever $u\in\mathbb R^m_+$ with $u_i = 0$;
$f(t,0) = 0$;  $f(t,u)$ is subhomogeneous in the sense that
$\alpha f(t,u)\le f(t,\alpha u)$ for all  $u\in\mathrm{Int}(\mathbb{R}^{m}_{+})$, $t\ge0$ and $\alpha\in (0,1)$.

\item\label{ass:c5'}
$\tilde f_u(t,u)$ is cooperative and irreducible for $t\ge0$ and $u\in\mathbb{R}^{m}_{+}$.
For each $1 \le i \le m$ and any $t\ge0$, $\tilde f_i(t,u)\ge0$ whenever $u\in\mathbb R^m_+$ with $u_i = 0$;
$\tilde f(t,0) = 0$;  $\tilde f(t,u)$ is strictly subhomogeneous in the sense that
$\alpha \tilde f(t,u)<\tilde f(t,\alpha u)$ for all $u\in\mathrm{Int}(\mathbb{R}^{m}_{+})$, $t\ge0$ and $\alpha\in (0,1)$.
\end{enumerate}
\begin{enumerate}[label=(C\arabic*),ref=(C\arabic*),leftmargin=*,start=6]
\item\label{ass:c6}
There exists a unit vector $e\in\mathbb{R}^n$, $t_0\ge0$, and a positive function $\chi$ such that
$$
\sum_{j=1}^{n}b^{(i)}_{j}(x,t)e_j\ge\chi(t), \quad x\in\overline\Omega,\ t\ge t_0, \ i=1,\dots,m.
$$
\end{enumerate}
For convenience, we define $p(t):=\frac{\chi(t)}{\alpha(t)}$ for $t\ge t_0$.
\begin{theorem}\label{threshold dynamics-D2}
Assume that $\lambda^*<0$. If $p$ is locally absolutely continuous on $[t_0,\infty)$ and there exists $T>0$ such that
\begin{equation}\label{threshold dynamics-D2-1}
\liminf_{t\to\infty}\frac{1}{T}\int_{t}^{t+T}\left(\frac{\chi(s)^2}{\alpha(s)}-2d_e\bigl(p'(s)\bigr)_+\right)\,\mathrm{d}s>4\Lambda_0(-\lambda^*),
\end{equation}
where $d_e=\max_{x\in\overline{\Omega}} x\cdot e-\min_{x\in\overline{\Omega}} x\cdot e$,
then every nonnegative solution of \eqref{asymptotically degenerate equations} converges to zero uniformly on $\overline\Omega$ as $t\to\infty$.
\end{theorem}

\begin{remark}\label{DuGuoPeng}
By comparison arguments, Theorems \ref{threshold dynamics-D} and \ref{threshold dynamics-D2} remain valid
when the reaction term $f$ in system \eqref{asymptotically degenerate equations} depends on the spatial variable,
provided that $f$ and its Jacobian matrix $f_u$ satisfy suitable spatially heterogeneous assumptions. In particular, we assume that
$f(x,t,u)$ and $f_u(x,t,0)$ converge to the spatially homogeneous limiting terms $\tilde f(t,u)$ and $\tilde f_u(t,0)$
uniformly with respect to $x\in\overline{\Omega}$ and locally uniformly in $u$ as $t\to\infty$.
\end{remark}

The sharpness of both the upper bound $4\lambda_0(-\lambda^*)$ in \eqref{threshold dynamics-D-1}
and the lower bound $4\Lambda_0(-\lambda^*)$ in \eqref{threshold dynamics-D2-1}
can be illustrated by an example similar to \eqref{Example-PD}. To avoid repetition,
we finish this section with two further examples, which demonstrate not only this sharpness but also some
possible connections between spreading phenomena in free boundary problems and our threshold type results.

\begin{example}\label{degenerate-example-1}
\rm
Consider the following free boundary problem:
\begin{equation}\label{free boundary problem1}
\begin{cases}
u_t-d u_{xx}=u(a-bu), & t>0,\quad 0<x<h(t),\\
u_x(0,t)=u(h(t),t)=0, & t>0,\\
h(0)=h_0,\quad u(x,0)=u_0(x), & 0\le x\le h_0,
\end{cases}
\end{equation}
subject to the Stefan condition
\begin{equation}\label{Stefan condition1}
h'(t)=-\mu u_x(h(t),t), \quad t>0.
\end{equation}
When spreading occurs, a central question is to examine the asymptotic spreading speed $k_0$ of the moving front $h(t)$.
It follows from \cite{Du-Lin2010,Bunting-Du-Krakowski2012,Du-Matsuzawa-Zhou2014} that $h\in C^{1+\alpha/2}_{\rm loc}([0,\infty))$ with $\alpha\in(0,1)$ and
\begin{equation}\label{asymptotic spreading speed-1}
k_0:=\lim_{t\to\infty}\frac{h(t)}{t}=\lim_{t\to\infty}h'(t).
\end{equation}
Moreover, they proved that $k_0$ is strictly smaller than the spreading speed determined by the classical Fisher--KPP equation on an unbounded domain, that is,
\begin{equation}\label{asymptotic spreading speed-2}
k_0<2\sqrt{ad}.
\end{equation}
This reduction in the spreading speed reflects the population loss at the spreading front.

We claim that \eqref{asymptotic spreading speed-2} coincides exactly with the condition \eqref{threshold dynamics-D-1} in Theorem \ref{threshold dynamics-D}(b).
Indeed, this can be seen by considering the corresponding reformulation \eqref{free boundary problem-vanishing-3} on the fixed interval $(0,1)$ and reflecting evenly across $y=0$ as in Example~\ref{asymptotically-periodic-example-1}.
Here, $\lambda_0$, $\Lambda_0$, $\lambda^*$, $\alpha$ and $\beta$ are precisely given by
$$
\lambda_0=\Lambda_0=d,\ \ \ \lambda^*=-a, \ \ \ \ \alpha(t)=\frac{1}{h^2(t)},\ \ \ \ \ \ \beta(t)=\frac{|h'(t)|}{h(t)}.
$$
In view of \eqref{asymptotic spreading speed-1}, we have
$$
\limsup_{t\to\infty}\frac{1}{T}\int_t^{t+T}\frac{\beta(s)^2}{\alpha(s)}\,\mathrm{d}s
=\limsup_{t\to\infty}\frac{1}{T}\int_t^{t+T}(h'(s))^2\,\mathrm{d}s=k_0^2.
$$
Thus the sufficient condition \eqref{threshold dynamics-D-1} reduces to
the known speed bound \eqref{asymptotic spreading speed-2}. It yields
$v(y,t)\to a/b$ locally uniformly for $y\in[0,1)$.
\end{example}

\begin{example}\label{degenerate-example-2}
\rm
Consider the nonlocal weighted free boundary model \eqref{free boundary problem1} with \eqref{nonlocal weighted free boundary condition}.
When spreading occurs, the existence and sharp characterization of the asymptotic spreading speed become challenging problems for this nonlocal weighted free boundary model.
Note that the model \eqref{free boundary problem1} with \eqref{nonlocal weighted free boundary condition} was generalized in \cite{Salako-Xue2026} to
a time-periodic case, where the intrinsic growth rate $a$ and the self-limitation rate $b$ are both positive time-periodic functions.
By the same argument as in Example \ref{degenerate-example-1}, condition \eqref{threshold dynamics-D-1} in Theorem \ref{threshold dynamics-D}(b) reduces to
$$
\limsup_{t\to\infty}\frac{1}{T}\int_t^{t+T}(h'(s))^2\,\mathrm{d}s<4\hat ad,
$$
where $\hat a$ denotes the average of $a(t)$ over one period.
When spreading occurs in a time-periodic environment, this inequality imposes a strict upper bound of $2\sqrt{\hat a d}$ on the long-time average of the boundary velocity. Interestingly, this bound matches the result for the time-periodic Stefan problem in \cite[Theorem~1.4 and Remark~1.5]{Du-Guo-Peng2013}.
\end{example}

\section{Applications}
\label{sec:applications}
In this section, we adopt a general framework, independent of specific biological processes,
to study reaction-diffusion-advection systems on time-varying domains.
Accordingly, no constitutive equations are introduced for the tissue mechanics, and the tissue is assumed to be incompressible
in the sense that its growth and deformation are not accompanied by any change in density.
For simplicity, we focus on two-dimensional domains that undergo dilation along each Euclidean coordinate axis.

Let $\Omega(t)\subset\mathbb{R}^2$ denote the tissue domain at time $t\ge0$.
Assume that, for each $t\ge0$, $\Omega(t)$ is a bounded and simply connected domain with boundary of class $C^{2+\theta}$ for some $0<\theta\le1$. 
The deformation of the tissue domain is described by a family of diffeomorphisms
\[
\Gamma(\cdot,t):\Omega_0\to\Omega(t),\qquad x=\Gamma(\mathcal X,t),
\]
which maps each material point $\mathcal X\in\Omega_0=\Omega(0)$ to its current position at time $t$.
By neglecting rigid-body translations and shear deformations of the tissue, the coordinate system can be appropriately chosen such that there is a reference point
which remains at the origin of the coordinate all the time, that is
\begin{equation}\label{reference-point}
\Gamma(0,t)=0, \quad t\ge0.
\end{equation}
From the definition of $\Gamma$, the initial condition is
\begin{equation}\label{initial-condition}
\mathcal{X}=\Gamma(\mathcal{X},0), \quad \mathcal{X}\in\Omega_0.
\end{equation}
In the spatial description, $a(x,t)$ gives the velocity at time $t$ of the particle which was at $\mathcal{X}=\Gamma^{-1}(x,t)$ initially, that is,
\begin{equation}\label{spatial-description-of-velocity}
a(x,t)=\dfrac{\partial\Gamma(\mathcal{X},t)}{\partial t}\Big|_{\mathcal{X}=\Gamma^{-1}(x,t)},
\end{equation}
The velocity gradient tensor $\nabla a$ describes the relative velocity of each particle with
respect to its neighbors.
By neglecting the antisymmetric component of the velocity gradient tensor $\nabla a$ associated with rigid body rotation,
we focus solely on its symmetric part, denoted by
$\mathfrak{D}$, referred to as the strain-rate tensor in continuum mechanics.
The off-diagonal elements of $\mathfrak{D}$ correspond to the shear rate between axes $x_i$ and $x_j$ with no accompanying volume change.
The diagonal components represent the rates of stretching along each coordinate direction. Accordingly, the trace ${\rm tr}\,\mathfrak{D}=\nabla\cdot a$
characterizes the local rate of volume expansion or contraction. By applying the chain rule, the deformation mapping $\Gamma$ is directly related to the strain-rate tensor
$\mathfrak{D}$ as follows
\begin{equation}\label{strain-rate-tensor}
\partial_t\nabla_{\mathcal X}\Gamma=\mathfrak{D}\nabla_{\mathcal X}\Gamma.
\end{equation}
We assume that the strain-rate tensor $\mathfrak{D}$ is determined locally due to factors that will not be explicitly
modelled, but which might include the effects of growth factors and cellular or subcellular structures giving rise to
mechanical anisotropy, influencing the direction of growth. Then the growth along one direction may be in a quite
different manner or rate to another direction. For simplicity, we restrict our attention to the following anisotropic
domain growth, of which the strain-rate tensor is specified by
$\mathfrak{D}_{i,j}=\delta_{i,j}\sigma_i(t),\, \,  i,j\in\{1,2\}$,
where $\delta_{i,j}$ denotes the Kronecker delta, and $\sigma_i$ is a continuous function defined on $[0,\infty)$.
When $\sigma_1(t)=\sigma_2(t)$ for any $t\ge0$, we have isotropic evolution of $\Omega(t)$.
Combining \eqref{reference-point}--\eqref{strain-rate-tensor}, we obtain the deformation mapping
\begin{align}\label{transformation}
x=\Gamma(\mathcal X,t)
  &=\big(\rho_1(t)\mathcal X_1,\rho_2(t)\mathcal X_2\big)^\top,
\end{align}
and hence
\[
a(\Gamma(\mathcal X,t),t)=\partial_t\Gamma(\mathcal X,t)
=\big(\dot\rho_1(t)\mathcal X_1,\dot\rho_2(t)\mathcal X_2\big)^\top,
\]
where $\rho_i(t)=\exp\left(\int_0^t\sigma_i(s)\,\mathrm{d}s\right)$  for $t\ge0$ and $i=1,2$.
Accordingly, we have
$\Omega(t)=\big\{(\rho_1(t)\mathcal X_1, \rho_2(t)\mathcal X_2);\ \mathcal X\in\Omega_0\big\}$.
For example, if $\Omega_0$ is the unit disk, namely $\Omega_0=\{(\mathcal X_1,\mathcal X_2)\in\mathbb R^2:\ \mathcal X_1^2+\mathcal X_2^2<1\}$,
then the evolving domain is given by
\[
\Omega(t)=\left\{(x_1,x_2)\in\mathbb R^2:\ \frac{x_1^2}{\rho_1^2(t)}+\frac{x_2^2}{\rho_2^2(t)}<1\right\},
\]
which is an ellipse centered at the origin for each fixed $t>0$.

Let $u(x,t)=(u_1(x,t), u_2(x,t),\dots,u_m(x,t))^\top$ denote the vector of densities of $m$ interacting species at
position $x\in\Omega(t)$ and time $t>0$. Here $u$ may also describe the concentration of $m$ chemical species or morphogens,
while the same mathematical formulation can be used for other interpretations, such as cell densities or effective genetic circuits.
The conservation law for $u$ over an elemental volume $V_t$ evolving in time is given by
\begin{equation}\label{conservation law}
\dfrac{\mathrm{d}}{\mathrm{d}t}\int_{V_t}u(x,t)\,\mathrm{d}x=\int_{V_t}\left(-\nabla\cdot J_{\rm diff}-\nabla\cdot J_{\rm adv}+f(u)\right)\,\mathrm{d}x,
\end{equation}
where $J_{\rm diff}$ denotes the diffusive flux associated with $u$, and $J_{\rm adv}$ denotes the environmentally induced advective flux,
which may be generated, for instance, by water currents, wind, or directional transport in the medium.
For simplicity, we take the advective velocity to be constant, $b=(r_1,r_2)^\top$,
where $r_1,r_2\in\mathbb{R}$. Then
$$
J_{\rm adv}=ub^\top, \qquad \nabla\cdot J_{\rm adv}=b\cdot\nabla u=r_1\frac{\partial u}{\partial x_1}+r_2\frac{\partial u}{\partial x_2}.
$$
Applying the Reynolds transport theorem on the left-hand side of \eqref{conservation law},
and using Fick's law of diffusion, we have the reaction-diffusion-advection system
\begin{equation}\label{reaction-diffusion-advection system-1}
\dfrac{\partial u}{\partial t}+(a\cdot \nabla) u+(\nabla\cdot a)u
=\sum_{j=1}^2d_j\frac{\partial^2 u}{\partial x_j^2}-b\cdot\nabla u+f(u).
\end{equation}
Here $d_j>0$ is the diffusivity in the $x_j$ direction, common to all species. The spatial
derivatives act componentwise. The velocity field $a$ induced by domain evolution is
given by \eqref{spatial-description-of-velocity}.
Incorporating the homogeneous Dirichlet boundary conditions
\begin{equation}\label{reaction-diffusion-advection system-2}
u(x,t)=0,\quad x\in\partial\Omega(t), \ t>0,
\end{equation}
and the initial condition
\begin{equation}\label{reaction-diffusion-advection system-3}
u(x,0)=u_0(x), \quad x\in\Omega_0,
\end{equation}
we formulate the initial-boundary value problem \eqref{reaction-diffusion-advection system-1}--\eqref{reaction-diffusion-advection system-3} on $\Omega(t)$.

Now we transform the spatial coordinates into the reference domain $\Omega_0$.
To simplify notation, we write $y\in\Omega_0$ in place of $\mathcal X$, and set $v(y,t)=u(\Gamma(y,t),t)$. By \eqref{transformation} and the chain rule,
problem \eqref{reaction-diffusion-advection system-1}--\eqref{reaction-diffusion-advection system-3} can be reformulated as
\begin{equation}\label{reaction-diffusion-advection system-4}
\left\{
\begin{aligned}
& \frac{\partial v}{\partial t}-A(t)v=F(t,v) \ \ \ &{\rm in}&\ \ \ \Omega_0\times(0,\infty),\\
& v =0 \ \ \ \ &{\rm on}&\ \ \ \partial\Omega_0\times(0,\infty),\\
& v(\cdot,0)=u_0 \ \ \ &{\rm in}&\ \ \ \Omega_0,
\end{aligned}
\right.
\end{equation}
where $A(t)v=\sum_{j=1}^2\frac{d_j}{\rho_j^2(t)}\frac{\partial^2 v}{\partial y_j^2}-\sum_{j=1}^2\frac{r_j}{\rho_j(t)}\frac{\partial v}{\partial y_j}$,
and $F(t,v)=-\sum_{j=1}^2\dfrac{\dot{\rho_j}(t)}{\rho_j(t)}v+f(v)$.

We impose the following assumptions on the reaction term.
\begin{enumerate}[label=(D\arabic*),ref=(D\arabic*),leftmargin=*]
\item\label{ass:d1}
$f\in C^{1}(\mathbb{R}^{m},\mathbb{R}^{m})$; $f_u(u)$ is cooperative and irreducible for $u\in\mathbb{R}^{m}_{+}$.
\item\label{ass:d2}
$f(0)=0$; $\alpha f(u)< f(\alpha u)$ for all $u\in\mathrm{Int}(\mathbb{R}^{m}_{+})$ and $\alpha\in (0,1)$.
\item\label{ass:d3}
There exists $\xi\gg0$ such that
$f_i(\xi)\le c(t)\xi_i$ for every $t\ge0$ and $1\le i\le m$,
where $c(t)=\sum_{j=1}^{2}\dot\rho_j(t)/\rho_j(t)$.
\end{enumerate}

To address  the case of an asymptotically periodic domain, we assume that
\begin{enumerate}[label=(E\arabic*),ref=(E\arabic*),leftmargin=*]
\item\label{ass:e1}
There exist two positive $\omega$-periodic functions $\rho_{1,\omega}$ and $\rho_{2,\omega}$ such that
\begin{align*}
\lim\limits_{t\to +\infty}  (\rho_j(t) -\rho_{j,\omega}(t))&= \lim\limits_{t\to +\infty} (\dot\rho_j(t)-\dot\rho_{j,\omega}(t))=0, \quad j=1,2.
\end{align*}
\end{enumerate}
In this case, the $\Omega(t)$ asymptotically exhibits periodic oscillations as $t\to\infty$.
The corresponding asymptotic limiting system is then given by
\begin{equation}\label{reaction-diffusion-advection system-case1-1}
\left\{
\begin{aligned}
& \frac{\partial v}{\partial t}-\tilde A(t)v=\tilde F(t,v) \ \ \ &{\rm in}&\ \ \ \Omega_0\times(0,\infty),\\
& v = 0 \ \ \ \ &{\rm on}&\ \ \ \partial\Omega_0\times(0,\infty),
\end{aligned}
\right.
\end{equation}
where $\tilde A(t)v=\sum_{j=1}^2\frac{d_j}{\rho_{j,\omega}^2(t)}\frac{\partial^2 v}{\partial y_j^2}-\sum_{j=1}^2\frac{r_j}{\rho_{j,\omega}(t)}\frac{\partial v}{\partial y_j}$,
and $\tilde F(t,v)=-\sum_{j=1}^2\dfrac{\dot{\rho_{j,\omega}}(t)}{\rho_{j,\omega}(t)}v+f(v)$.

Let $\lambda^*_\omega$ be the principal eigenvalue of the $\omega$-periodic parabolic eigenvalue problem:
\begin{equation*}
\left\{
\begin{aligned}
& \frac{\partial v}{\partial t}-\tilde A(t)v=\tilde F_u(t,0)v+\lambda v \ \ \ &{\rm in}&\ \ \ \Omega_0\times\mathbb{R},\\
& v=0 \ \ \ \ &{\rm on}&\ \ \ \partial\Omega_0\times\mathbb{R},\\
& v(y,\omega+t)=v(y,t) \ \ \ &{\rm in}&\ \ \ \Omega_0\times\mathbb{R}.
\end{aligned}
\right.
\end{equation*}
According to Theorem \ref{general threshold dynamics}, if $\lambda^*_\omega<0$, then \eqref{reaction-diffusion-advection system-case1-1} admits a unique
positive $\omega$-periodic solution, which we denote by $v_\omega^*$. As a consequence of Theorem \ref{general threshold dynamics},  we  have the following threshold type result.
\begin{theorem}\label{reaction-diffusion-advection system-case1-2}
Let \ref{ass:e1} hold. Then the following statements are valid:
\begin{itemize}
\item [(a)] If $\lambda_\omega^{\ast}\ge0$, then the zero solution is globally asymptotically stable
		for system \eqref{reaction-diffusion-advection system-4}
		in $C_{0}^1(\overline\Omega_0, \mathbb{M}_{\xi})$.
		\item [(b)] If $\lambda_\omega^{\ast}<0$, then every nontrivial solution $v$ of  system \eqref{reaction-diffusion-advection system-4} with initial data in $C_{0}^1(\overline\Omega_0,\mathbb{M}_{\xi})$
                    satisfies $\lim\limits_{t\to+ \infty} \|v(\cdot,t)-v_\omega^*(\cdot,t)\|_{C^1(\overline\Omega_0)}=0$.
\end{itemize}
\end{theorem}

\begin{remark}\label{reaction-diffusion-advection system-7}
Note that $\rho_1(0)=\rho_2(0)=1$ and $|\Omega(t)|=\rho_1(t)\rho_2(t)|\Omega_0|$ for all $t>0$.
Then Theorem \ref{reaction-diffusion-advection system-case1-2} with $\rho_2(t)\equiv\frac{1}{\rho_1(t)}$ also
gives rise to threshold dynamics for the area-preserving case addressed in \cite{Gorder-Klika-Krause2021}.
\end{remark}

Next, we consider the asymptotically partially degenerate case, restricting attention to rectangular evolving domains $\Omega(t)=\Omega_1(t)\times\Omega_2(t)\subset\mathbb{R}^{2}$, where each $\Omega_i(t)\subset\mathbb{R}$ ($i=1,2$) is a bounded interval evolving from an initial interval $\Omega_{0,i}$ containing the origin. Note that, by an approximation argument, the $C^{2+\theta}$-regularity assumption on $\partial\Omega(t)$ in this section can be relaxed to Lipschitz regularity. We further assume that
\begin{enumerate}[label=(E\arabic*),ref=(E\arabic*),leftmargin=*,start=2]
\item\label{ass:e2}
$\rho_1(t)=(1+t)^q$ with $q\in(0,1)$,  and  $ \rho_2(t)=1+e^{-t}\sin t$, \, \,   $\forall  t\ge 0$.
\end{enumerate}
Then \(\Omega(t)\) can be viewed as a strip-like evolving domain whose length becomes asymptotically unbounded, while its cross-section remains asymptotically bounded.

It is easy to see that the limiting system associated with \eqref{reaction-diffusion-advection system-4} takes the form:
\begin{equation}\label{reaction-diffusion-advection system-case2-1}
\left\{
\begin{aligned}
& \frac{\partial v}{\partial t}-d_2\frac{\partial^2 v}{\partial y_2^2}+r_2\frac{\partial v}{\partial y_2}=f(v) \ \ \ &{\rm in}&\ \ \ \Omega_{0,2}\times(0,\infty),\\
& v = 0 \ \ \ \ &{\rm on}&\ \ \ \partial\Omega_{0,2}\times(0,\infty).
\end{aligned}
\right.
\end{equation}
Define the operator $[Lu](y_2):= d_2\frac{\partial^2 u}{\partial y_2^2}-r_2\frac{\partial u}{\partial y_2}+f_v(0)u$
for $u\in C_0^{2}(\overline{\Omega_{0,2}},\mathbb{R}^{m})$, and let \(\lambda^*\) denote the principal eigenvalue of $L$.
By Theorem \ref{threshold dynamics-PD1}, if \(\lambda^*>0\), then system
\eqref{reaction-diffusion-advection system-case2-1} admits a unique positive $\omega$‑periodic solution
$v^*$. Thus, Theorem \ref{threshold dynamics-PD2} leads to  the following result.

\begin{theorem}\label{reaction-diffusion-advection system-case2-2}
Let \ref{ass:e2} hold. Then the following statements are valid for system \eqref{reaction-diffusion-advection system-4}:
	\begin{itemize}
		\item [(a)] If $\lambda^{\ast}\le0$, then every solution $v$ with $v(\cdot,0)\in C_0(\overline\Omega_0,\mathbb{M}_{\xi})$ satisfies $\lim\limits_{t\to+ \infty}v(\cdot,t)=0$ uniformly on $\overline\Omega_0$.
		\item [(b)] If $\lambda^{\ast}>0$ and $|r_1|<2\sqrt{d_1\lambda^*}$,
                    then every nontrivial solution $v$ with $v(\cdot,0)\in C_0(\overline\Omega_0,\mathbb{M}_{\xi})$ satisfies
                    $\lim\limits_{t\to+ \infty}\sup_{(y_1,y_2)\in K\times\overline\Omega_{0,2}}|v(y_1,y_2,t)-v^*(y_2)|=0$
                   for every compact subset $K$ of $\Omega_{0,1}$.
        \item[(c)]  If $\lambda^{\ast}>0$ and $|r_1|>2\sqrt{d_1\lambda^*}$, then every nonnegative solution converges to zero uniformly on $\overline\Omega_0$ as $t\to\infty$.
	\end{itemize}
\end{theorem}

\section*{Acknowledgments}
Ye's research was supported in part by NSFC (12271186, 12271178).
Zhao's research was supported in part by the NSERC of Canada (RGPIN-2025-04963).

\end{document}